\documentclass[letterpaper, 11pt]{article}
\usepackage[left=1in, right=1in, top=1in, bottom=1in]{geometry}
\usepackage[T1]{fontenc}
\usepackage{amsmath, amsthm, amssymb, mathtools}
\usepackage{natbib}
\usepackage{graphicx}
\usepackage{xcolor}
\usepackage{caption}
\usepackage{subcaption}
\usepackage{setspace}
\usepackage{authblk}
\usepackage{enumitem}
\usepackage{float}
\usepackage{algorithm}
\usepackage{algpseudocode}
\usepackage{multirow}
\usepackage{tabularx} % Add this in the preamble
\usepackage{booktabs}
\usepackage{tikz}
\usetikzlibrary{shapes,arrows.meta,fit}
\newtheorem{remark}{Remark}

\newtheorem{example}{Example}

\graphicspath{{images/}}

\title{Improving Full Strong Branching Decisions by Incorporating Additional Information}
\author[1]{Prachi Shah\thanks{prachi.shah6795@gmail.com}}
\author[2]{Santanu S. Dey\thanks{santanu.dey@isye.gatech.edu}}
\affil[1]{Amazon}
\affil[2]{School of Industrial and Systems Engineering, Georgia Institute of Technology}

\date{\vspace{-5ex}}

\begin{document}

\maketitle              % typeset the header of the contribution

\begin{abstract}
The full strong branching (FSB) rule is well known to produce extremely small branch-and-bound trees. 
This rule guides branching decisions based exclusively on the information regarding local gains in the linear programming (LP) bounds. 
We identify and correct two key shortcomings in FSB. First, the LP gains may be overestimations of the improvement in global dual bounds whenever pruning is possible. We propose a modification to address this issue, that incorporates primal bounds and readjusts the relative importance of the larger and smaller LP gains.  
Second, \textcolor{black}{FSB decisions may be myopic as they consider only local LP gains and cannot foresee the impact of branching decisions on feasibility or integrality beyond immediate children. To address this weakness, we present an approach that detects global asymmetry trends in infeasibility and integrality due to 0 and 1 assignments and incorporates them into the FSB score function. We further extend this approach to achieve more balanced trees even when the branch-and-bound tree prunes primarily by bounds.}

Using randomly generated problem instances with known structures, we derive insights and fine-tune our modified scores. 
Evaluation on MIPLIB 2017 Benchmark instances shows a 22-35\% reduction in mean tree sizes for solved cases and a 3.6-5.6\% decrease in the remaining gap for unsolved ones. Our approach extends to reliability branching (RB), where improved scores reduce mean tree sizes by 5-13\% on solved instances and lower the mean gap by 2.6-4.3\% on unsolved instances, depending on primal bound quality.

%Modern state-of-the-art solvers often find high-quality solutions quickly, providing good primal bounds for the problem early on. 
%We highlight how this information can be incorporated into FSB to possibly improve this rule. 
%However, this information is not considered while evaluating branching scores. 
%Furthermore, in many problems, there is an inherent asymmetry between the size of tree corresponding to the child nodes for the 0 and 1 assignments for binary variables. 
%But, strong-branching scores are symmetric with respect to the dual bound improvements on both sides. 
%$In this work, we show how to estimate this asymmetry and use it to possibly reduce the size of the branch-and-bound tree generated by FSB. 
%incorporate the information on primal bounds and the asymmetry of problem structure into strong-branching decisions by refining the score function. 
%Based on the information regarding the problem structure known to the user, we propose 4 alternate functions for the strong-branching score. % No structural asymmetry => eff sb only counted as the 4th function
%We provide extensive computational results on randomly generated instances and instances from the MIPLIB 2017 Benchmark Set. <Modify to state results>


\end{abstract}
\section{Introduction}

% \begin{enumerate}
%     \item Describe vanilla branch-and-bound
%     \item Branching decisions and importance
%     \item Strong branching
%         \begin{itemize}
%             \item Description, different scores 
%             \item Approximations - most work towards improving time (reliability etc)
%         \end{itemize}
%     \item Motivation
%         \begin{itemize}
%             \item Difference from Optimal trees
%             \item How to make better decisions fundamentally - what is important
%             \item With the advent of machine learning, better quality decisions for offline training
%             \item For both of the above reasons, we are interested in smaller trees and not quicker solution time
%         \end{itemize}
%     \item Contributions
% \end{enumerate}

The branch-and-bound~\citep{land60a} algorithm is one of the most popular frameworks for solving mixed-integer programs (MIPs) and forms the backbone of most modern MIP solvers. 
Linear programming (LP) relaxations are solved ignoring the integrality requirements to get a dual bound. The feasible space is then recursively partitioned to create two subproblems. This operation is called \textit{branching}. This process continues recursively until all leaf nodes are \textit{pruned} either due to infeasibility, having an integral LP solution, or due to the LP bound being worse than a known integer solution~\citep{nemhauser1988integer}. %\red{references for more details.}

In this work, we consider mixed-binary linear programs. Branching is typically performed based on the disjunction $x_j \leq 0 \vee x_j \geq 1$, where $x_j$, called the \textit{branching variable}, is selected from the set of variables violating the integrality requirement in the LP solution $x^*$. 
The choice of branching variable is critical, as it greatly impacts the size of the branch-and-bound tree and, consequently, the solution time. 
This has been demonstrated both theoretically and computationally in the literature~\citep{achterberg2007constraint,le2015important,dey2024theoretical}.
The strategy used to choose the branching variable is called the \textit{branching rule}. Typically, branching rules assign a score to each variable, and the variable with the highest score is preferred for branching.

The branching rule that is empirically known to produce smaller trees than those by other heuristics~\citep{achterberg2005branching} is the \textit{full strong branching rule} (FSB)~\citep{applegate1995finding}.  
It computes the improvement in LP bounds due to branching on every fractional variable, henceforth referred to as the LP gains, by solving (to optimality) $2$ LPs per variable. Let $\Delta^0_i$ and $\Delta^1_i$ be the LP gains, that is, the difference in the LP value of a node and its children when %tentatively branching on variable $i$.
fixing $x_i$ to $0$ and $1$ respectively.
% \footnote{If a child node is infeasible, the corresponding LP gain is infinite. If multiple branching candidates lead to infeasible children, we branch on the variable with two infeasible nodes if such exists, else, maximize the gain on the feasible branch.}. 
Then $\Delta^0_i$ and $\Delta^1_i$ are aggregated to obtain a score for variable $i$.
The product score function is the most commonly used function for aggregation due to its robust performance~\citep{achterberg2005branching,achterberg2007constraint,linderoth1999computational}. The product score function is given by, 
\begin{equation} \label{eq:default_sb}
\text{score}(i) = \max\{\Delta^0_i, \epsilon\} \cdot \max\{\Delta^1_i, \epsilon\},
\end{equation}
where $\epsilon > 0$ is small. %chosen to be close to 0 

A 
%\red{(our?)} 
recent study~\citep{dey2024theoretical} compared the size of FSB trees with that of the \textit{optimal trees}, that is, the smallest branch-and-bound tree that can solve a given instance. It was 
%They \red{(or we?)} 
demonstrated on small instances of many classical problems that although FSB produces smaller trees than all other heuristics, the FSB trees can still be 25-75\% larger than optimal trees on average. In this work, our goal is to improve FSB decisions in order to bridge this difference in tree size from optimal  trees. 

In this paper, we are interested in the tree size rather than the solution time for two major reasons. First,
note that FSB is never implemented in the textbook form described above since it is very expensive to take every branching decision.
However, with the advent of machine learning-based approaches, 
%that can imitate good decisions quickly, better decisions can provide superior experts to learn from.
it is possible to imitate FSB decisions~\citep{alvarez2017machine,khalil2016learning,gasse2019exact,nair2020solving,gupta2020hybrid, zarpellon2021parameterizing}. 
%Therefore, several important advances in using machine learning to improve branching decisions~\citep{alvarez2017machine,khalil2016learning,gasse2019exact,nair2020solving,gupta2020hybrid, zarpellon2021parameterizing} have employed FSB as the expert method to be imitated.
Moreover, while \citep{dey2024theoretical} provides an algorithm to compute optimal trees, this algorithm can only be run for extremely small instances, and learning from such small instances is not practical.
%optimal trees yield the optimal branching decisions, they can be constructed only for instances too small to be interesting in practice. 
Therefore, one of the motivations for this paper is to discover simple rules that can improve over FSB, which can then be used as expert for machine learning algorithms to imitate. Second, we are interested in understanding the science 
behind making 
good branching decisions. %fundamentally.
%how to make good branching decisions fundamentally. 
The success of FSB indicates that LP gains are critical information in making good decisions. However, relying solely on LP gains may result in myopic branching decisions. In this work, we aim to identify and utilize additional information that is readily available in solvers to enhance the effectiveness of FSB.
%However, only considering LP gains can lead to myopic decisions, and we are interested in understanding what other information among those that are readily available in solvers, can be used to improve FSB decisions.

% \red{
% Important with the advent of machine learning based approaches that can imitate good decisions quickly and need for better decisions to learn from. How to make better decisions fundamentally. Success of SB indicates objective gains are important, but what else is? Theory on branching is weak. This line of questioning also identifies features important for machine learning approaches. }
% \red{Our approach - starting from SB, improve it by incorporating more information and integrating it into the strong branching score.}

Researchers in the past have developed several heuristics that collect and leverage different types of information. 
%A popular class of rules consists of derivatives of FSB, attempts to maximize the changes in the objective bounds. 
A popular class of branching rules is based on modified versions of FSB that also seek to maximize LP gains.
The challenge with FSB is the high computational burden of solving two LPs per variable to obtain the value of LP gains. 
%Processing each tree node with FSB is very slow due to the computational burden of solving the additional LPs.
Therefore, these rules use approximations of LP gains, prominent among them being pseudo-cost branching~\cite{benichou1971experiments, linderoth1999computational}, and reliability branching~\citep{achterberg2007constraint}. These rules are quicker at solving MIPs in practice, where they use objective gains per unit change in variable assignment from past branching operations to estimate current LP gains. 
However, they result in much larger trees than FSB. In contrast, look-ahead branching considers objective gains two levels deeper into the branch-and-bound tree~\citep{glankwamdee2006lookahead} resulting in smaller trees at a higher computational cost.

% \margin{maybe remove the trade-off and complexity angle - needless} There is typically a trade-off between the computational cost of obtaining the complex information and the quality of the branching decision. 
% At one end of the spectrum, 
In the context of constraint satisfaction problems, such as the boolean satisfiability (SAT) problem, different branching rules have been developed since the objective is no longer relevant here. One such rule takes into account the reductions in domains of other variables due to fixing a candidate variable~\citep{li1997look}. Another considers the frequency with which a variable appears in conflict clauses leading to infeasibility \citep{moskewicz2001chaff}. In \citep{kilincc2009information}, the notion of conflict clauses has been extended to the MIP setting to pruned nodes and a branching rule is proposed based on the length of clauses a variable is present in. Hybrid branching \citep{achterberg2009hybrid} combines the above-mentioned rules with reliability pseudo-cost branching by taking a convex combination of the scores for each variable.

Other types of information have also been used to make branching decisions. The most-fractional rule uses the fractionality of the variables in $x^*$ to make quick decisions but leads to relatively large tree sizes. In \citep{gilpin2006information}, the authors treat fractionality as uncertainty and design a branching rule to minimize entropy after branching.
In \citep{patel2007active}, the authors develop a branching rule that is based on active constraints in the optimal LP solutions. In \citep{fischetti2011backdoor}, a subset of important branching candidates called backdoors is identified by sampling near optimal fractional solutions and solving a set cover problem over the fractional variables. Similarly, a Monte Carlo Tree Search framework for finding backdoors is presented in \citep{khalil2022finding}.

In this work, we highlight and fix two shortcomings of FSB. Our contributions are summarized below:
\begin{itemize}
    \item We demonstrate that LP gains may be overestimations of improvement in global dual bounds whenever pruning is possible. We present an approach that uses primal bounds to correct the overestimations and a modification of the product score that is tuned for the corrected gains.
    \item We observe that FSB greedily maximizes improvements in bounds and cannot foresee the impact of branching decisions on nodes beyond its immediate children. We leverage global trends in infeasibility to counteract this shortcoming through a further generalization of the product score.
    \item We use randomly generated instances of problems with known structure to derive insights and tune the modified scores. On these problems, the new rules, on average, lead to 10-15\% smaller trees than those based on the FSB, where the range of improvements depend on the quality of primal bound available at the time of initializing the branch-and-bound tree.
    \item We evaluate the proposed rules on instances from MIPLIB 2017 Benchmark Set~\citep{miplib2017} and demonstrate around 15-26\% smaller mean trees on solved instances and 4-6\% decrease in mean gap remaining on unsolved instances, depending on the quality of primal bounds. We show that our approach can be directly extended to reliability branching (RB). Similar experiments on MIPLIB instances show that the improved RB scores decrease the mean tree sizes by 3-15\% on solved instances, depending on the quality of primal bounds, and reduce the gap remaining by 2-3\% for unsolved instances.
\end{itemize}

The rest of this paper is organized as follows. The remainder of this section describes the experimental setup used throughout this paper. In Section~\ref{sec:primal}, we discuss the issue of  overestimation of LP gains and present methods to correct for it using primal bounds. In Section~\ref{sec:la}, we propose a method to identify global trends in infeasibility and incorporate them into a generalization of the product score. In Section~\ref{sec:RLA}, we generalize the results of Section~\ref{sec:la} to account for general asymmetry trends in the branch-and-bound tree. In Section~\ref{sec:MIPLIB}, we evaluate the newly designed scores as well as their extension to RB on MIPLIB instances. We discuss the key findings and conclude in Section~\ref{sec:conclusion}.

\subsection{Experimental set-up} \label{sec:intro_expt_setup}
Experiments are performed on the Red Hat Enterprise Linux 8 workstation. 
% powered by \blue{Dual Intel Xeon Gold 6226 processors (2.7 GHz, 64-bit)}.
We implement a Python 3.10-based textbook branch-and-bound algorithm to isolate the impact of branching decisions, unaffected by other sophisticated processes of modern solvers. Gurobi 11.0 is used as the LP solver. Nodes were processed in the best-bound first order, which produces the smallest trees~\cite{dey2023lower, shah2025non}. %\textcolor{brown}{Moreover, the optimal solution value is used for pruning to ensure that only nodes necessary for certifying optimality are explored.} 
Finally, neither presolve nor cutting planes are employed unless otherwise stated.

FSB is implemented purely as an oracle returning a branching variable, without influencing any other process within branch-and-bound algorithm, thus allowing fair node counting \citep{gamrath2018measuring}. All fractional variables in the LP solution are considered for branching, and objective gains are computed by solving the LPs to optimality. 
{If a child node is infeasible, the corresponding LP gain is infinite. If multiple branching candidates lead to infeasible children, we branch on the variable with two infeasible nodes if such exists, else, maximize the gain on the feasible branch.}

In order to test various ideas for improving FSB, we created a library of MILP instances. 
In particular, we randomly generate 20 instances of $14$ problems from the literature to create a test bed for the experiments in Sections~\ref{sec:primal} and~\ref{sec:la}. These problems include, 
\begin{itemize}
    \item Multi-dimensional Knapsack: The ratio of right-hand side vector to sum of constraint coefficients is 25\% for those labelled ``small'', 50\% for ``medium'' and 75\% for ``large''.
    \item 0-1 coefficient Packing and Covering problems: Including set-packing, set-covering, combinatorial auctions, and independent set and weighted matching on graphs.
    \item Classical mixed-integer problems: Including capacitated lot-sizing, big-bucket lot sizing, capacitated facility location, fixed charge network flow, portfolio optimization with chance constraints and cardinality constrained weighted set coverage.
\end{itemize}
A detailed description of these problems and the distribution of instances is presented in Appendix~\ref{sec:random_instances}. 
Note that these problems have varying numbers of binary variables, since we wanted to have a test bed of instances with tree sizes that were neither too large nor too small. We will also test on Benchmark Set of MIPLIB 2017~\citep{miplib2017} instances in Section~\ref{sec:MIPLIB}. 

Tree sizes are defined as the total number of nodes in the branch-and-bound tree, and a shift of 100 nodes was used to compute shifted geometric means.

\section{Overestimation in LP gains}\label{sec:primal}

Consider the case where, upon branching on a variable, the LP gain exceeds the integrality gap at the node. As long as the LP gain is at least the integrality gap, the child node will be pruned by bound, and any additional LP gain overestimates the possible improvement in global dual bound and may be regarded as redundant. The example below shows how such overestimation may lead to FSB taking poor branching decisions. 
\begin{example}\label{ex:1}
    Consider a node where the additive integrality gap is 1, and there are two candidate variables, $x_i$ with LP gains $(\Delta^0_i, \Delta^1_i)=(0.6, 1.8)$ and $x_j$ with LP gains $(\Delta^0_j, \Delta^1_j)=(0.9, 1)$. Since the second child closes the integrality gap for both variables, we prefer branching on the $x_j$ variable, where the LP gain at the remaining node is larger. However, FSB picks the variable $x_i$ to branch on, as it overestimates the usefulness of the LP gain 1.8. 
\end{example}
% Another scenario worth considering is that strong branching always picks a variable that creates exactly one infeasible child over another, generating feasible children even when both new nodes may be pruned by bound.  

In this section, motivated by the aforementioned example, we use primal bounds to identify and correct such overestimation in LP gains and improve overall branching decisions.
%The above example shows that we should incorporate estimates of gap needed to be closed to improve the performance of FSB. 
The estimate of the gap to be closed is directly tied to the quality of the primal solution. Modern state-of-the-art solvers typically find good quality, near-optimal primal solutions very quickly, providing a good primal bound for the problem early on~\citep{berthold2006primal,berthold2013measuring,berthold1945primal}. 

Solutions found by heuristics are not guaranteed to be optimal, so we use the notion of the primal gap %is used 
to measure their difference from the optimal value.
Let $z^*$ be the optimal value and $\hat{z}$ the primal bound that may be obtained, for example, from a known sub-optimal primal feasible solution. We define the primal gap as, 
\begin{equation*}
    \text{primal gap } = \dfrac{|z^* - \hat{z}|}{|z^*|}.
\end{equation*}
In practice, the primal gap is often positive at the start of the branch-and-bound process. Therefore, an important consideration is to obtain a branching rule that is robust for reasonable primal gaps.

% \red{
% To discuss - 
% \begin{itemize}
%     \item Grid search on small random and experiments on miplib. Skip large random here.
%     \item Remove Dom SB in results and replace by Default SB
%     \item Remove discussion comparing to Dom SB which considers a subset of variables where some node is pruned. Or should this be briefly discussed with detailed results in appendix?
%     \item update tables to show \% difference in stead?
% \end{itemize}
% }
% With the observation mentioned above in mind, 
\subsection{Efficacious gains}
We here propose a modification of the FSB score that factors in the primal bounds.
We apply a threshold on the LP gains $(\Delta^1_i, \Delta^0_i)$ and define \emph{efficacious gains}, $(q^1_i, q^0_i)$ as follows. Let $\hat{z}$ be the incumbent primal bound of the integer program and let $z^{LP}_N$ be the LP value at the node being processed for branching. Then the additive primal-dual gap is $\Delta_{p-d} = |\hat{z} - z^{LP}_N|$ . Note that this may be greater than the integrality gap at the node if the primal bound is not optimal. Then,
% \begin{equation*}
%     q^0_i = \max\big\{ \min (\Delta^0_i, \Delta_{p-d}), \ \epsilon\big\}, \ q^1_i = \max\big\{ \min (\Delta^1_i, \Delta_{p-d}), \ \epsilon\big\}
% \end{equation*}
\begin{equation*}
    q^0_i =  \min \big\{ \max\{\Delta^0_i, \epsilon\}, \ \Delta_{p-d} \big\}, \ q^1_i =\min \big\{ \max\{\Delta^1_i, \epsilon\}, \ \Delta_{p-d} \big\}, 
\end{equation*}
This is also applied to infeasible child nodes, where an infinite LP gain leads to an efficacious gain of $\Delta_{p-d}$. We propose computing strong branching scores using the efficacious gains instead of the actual LP gains. We call this the \textit{efficacious} strong branching (Eff-SB) implementation. 
% \margin{Rename everything to Eff-FSB, Def-FSB?}

For the first set of experiments, we evaluate the performance of Eff-SB (i.e., where $\textup{score}(i) = q^0_i\cdot q^1_i$) against the default product rule in Eq.\eqref{eq:default_sb}, henceforth referred to as Def-SB. 
% First, we assume the optimal IP value is known and use the exact integrality gap at a node to compute efficacious gains. 
To simulate primal heuristics in solvers, we initialize the branch-and-bound process with the primal bound corresponding to 4 different primal gaps in our experiments - 0\%, 2\%, 5\% and 10\%.
The bound may be updated whenever improving integral solutions are discovered during the course of the branch-and-bound tree.

Table \ref{tab:effSB} presents the shifted geometric means of branch-and-bound tree sizes with Def-SB and Eff-SB rules on all the problems in the test bed. 
The same information is shown in terms of percentage reduction in tree sizes relative to Def-SB in Table~\ref{tab:effSB_pct}.
The first row denotes the branching rule, and the second denotes the initial primal gap whenever efficacious gains are used. The first column presents the problem with the numbers in parentheses indicating the number of binary variables, and the second column presents the percentage of leaf nodes that were pruned due to infeasibility or integrality, averaged across all instances.
% While the improvement in tree size is small, it is consistently between 1.7-2.8\% across all primal gaps, indicating a robustness to the quality of primal bounds that was absent in the previous approaches.
% All three scoring functions show some improvement on average across tree sizes, though the tree sizes for some problem classes increase.
\begin{table}
\begin{center}	
\caption{Shifted geometric mean of tree sizes for Def-SB and Eff-SB rounded to integer values. Tree sizes smaller than Def-SB are highlighted in bold. }\label{tab:effSB}
\resizebox{\textwidth}{!}{
\setlength{\tabcolsep}{1.em} 
  {\renewcommand{\arraystretch}{1.1}
\begin{tabular}{|l|c|r|rrrr|}
\hline
Rule                                &      & \multicolumn{1}{c|}{Def-SB} & \multicolumn{4}{c|}{Eff SB} \\ \hline
Primal Gap &
   &
  \multicolumn{1}{c|}{} &
  \multicolumn{1}{c}{0\%} &
  \multicolumn{1}{c}{2\%} &
  \multicolumn{1}{c}{5\%} &
  \multicolumn{1}{c|}{10\%} \\ \hline
\multicolumn{1}{|c|}{Problem} &
  Inf. Leaves &
   &
  \multicolumn{1}{r}{} &
  \multicolumn{1}{r}{} &
  \multicolumn{1}{r}{} &
  \multicolumn{1}{r|}{} \\ \hline
Big-Bucket Lot-Sizing (30)          & 6\%  & 3099    & 3151  & \textbf{3080}  & 3103  & 3109  \\
Capacitated   Lot-Sizing (30)       & 1\%  & 5302    & \textbf{5265}  & \textbf{5164}  & \textbf{5283}  & \textbf{5292}  \\
Fixed Charge   Flow (150)           & 8\%  & 7658    & \textbf{7646}  & \textbf{7613}  & 7675  & 7697  \\
Capacited   Facility Location (100) & 0\%  & 3733    & \textbf{3725}  & \textbf{3722}  & \textbf{3722}  & \textbf{3722}  \\
Portfolio   Optimization (100)      & 0\%  & 19965   & \textbf{19945} & \textbf{19776} & \textbf{19793} & \textbf{19794} \\
Weighted Set   Coverage (200)       & 0\%  & 1274    & \textbf{1252}  & \textbf{1255}  & \textbf{1255}  & \textbf{1255}  \\
Set Covering   (3000)               & 0\%  & 2232    & \textbf{2083}  & \textbf{2088}  & \textbf{2115}  & \textbf{2125}  \\
Set Packing   (1000)                & 0\%  & 5362    & \textbf{4828}  & \textbf{4815}  & \textbf{4805}  & \textbf{4802}  \\
Combinatorial   Auctions (200)      & 0\%  & 1493    & \textbf{1245}  & \textbf{1335}  & \textbf{1337}  & \textbf{1337}  \\
Independent   Set (500)             & 0\%  & 1701    & 1769  & \textbf{1670}  & \textbf{1694}  & \textbf{1688}  \\
Weighted   Matching (1000)          & 0\%  & 3545    & 3679  & 3587  & 3587  & 3587  \\
Multi-dim   Knapsack - Small (100)  & 67\% & 1014    & 1397  & 1339  & 1271  & 1204  \\
Multi-dim   Knapsack - Med (100)    & 40\% & 16196   & 16395 & 16206 & 16423 & 16482 \\
Multi-dim   Knapsack - Large (100)  & 24\% & 18171   & 20128 & 19248 & 19241 & 19240 \\ \hline
Geo Mean                            &      & 4149    & 4190  & 4141  & 4149  & \textbf{4137}  \\ \hline
\end{tabular}
}}
% \end{center}
% \end{table}
\vspace{24pt}
% \begin{table}
% \begin{center}	
\caption{Comparison of shifted geometric mean of tree sizes for the  Def-SB and Eff-SB rules. Tree sizes are presented for Def-SB and percentage reduction in tree size relative to Def-SB are presented for others. }\label{tab:effSB_pct}
\resizebox{\textwidth}{!}{
\setlength{\tabcolsep}{1.em} 
  {\renewcommand{\arraystretch}{1.1}
\begin{tabular}{|l|c|r|rrrr|}
\hline
Rule                                &      & \multicolumn{1}{c|}{Def-SB} & \multicolumn{4}{c|}{Eff SB} \\ \hline
Primal Gap &
   &
  \multicolumn{1}{c|}{} &
  \multicolumn{1}{c}{0\%} &
  \multicolumn{1}{c}{2\%} &
  \multicolumn{1}{c}{5\%} &
  \multicolumn{1}{c|}{10\%} \\ \hline
\multicolumn{1}{|c|}{Problem} &
  Inf. Leaves &
   &
  \multicolumn{1}{r}{} &
  \multicolumn{1}{r}{} &
  \multicolumn{1}{r}{} &
  \multicolumn{1}{r|}{} \\ \hline
Big-Bucket Lot-Sizing (30)        & 6\%  & 3099                        & -2\%          & \textbf{1\%}  & 0\%           & 0\%            \\
Capacitated Lot-Sizing (30)       & 1\%  & 5302                        & \textbf{1\%}  & \textbf{3\%}  & 0\%  & 0\%   \\
Fixed Charge Flow (150)           & 8\%  & 7658                        & 0\%  & \textbf{1\%}  & 0\%           & -1\%           \\
Capacited Facility Location (100) & 0\%  & 3733                        & 0\%  & 0\%  & 0\%  & 0\%   \\
Portfolio Optimization (100)      & 0\%  & 19965                       & 0\%  & \textbf{1\%}  & \textbf{1\%}  & \textbf{1\%}   \\
Weighted Set Coverage (200)       & 0\%  & 1274                        & \textbf{2\%}  & \textbf{1\%}  & \textbf{1\%}  & \textbf{1\%}   \\
Set Covering (3000)               & 0\%  & 2232                        & \textbf{7\%}  & \textbf{6\%}  & \textbf{5\%}  & \textbf{5\%}   \\
Set Packing (1000)                & 0\%  & 5362                        & \textbf{10\%} & \textbf{10\%} & \textbf{10\%} & \textbf{10\%}  \\
Combinatorial Auctions (200)      & 0\%  & 1493                        & \textbf{17\%} & \textbf{11\%} & \textbf{10\%} & \textbf{10\%}  \\
Independent Set (500)             & 0\%  & 1701                        & -4\%          & \textbf{2\%}  & 0\%  & \textbf{1\%}   \\
Weighted Matching (1000)          & 0\%  & 3545                        & -4\%          & -1\%          & -1\%          & -1\%           \\
Multi-dim Knapsack - Small (100)  & 67\% & 1014                        & -38\%         & -32\%         & -25\%         & -19\%          \\
Multi-dim Knapsack - Med (100)    & 40\% & 16196                       & -1\%          & 0\%           & -1\%          & -2\%           \\
Multi-dim Knapsack - Large (100)  & 24\% & 18171                       & -11\%         & -6\%          & -6\%          & -6\%           \\ \hline
Geo Mean                          &      & 4149                        & -1.0\%        & 0.2\%         & 0.0\%         & \textbf{0.3\%} \\ \hline
\end{tabular}
}}
\end{center}
\end{table}

The following are the main insights:

\begin{enumerate}
\item The results show that Eff-SB leads to smaller tree sizes on some problems such as combinatorial auctions, set packing and set covering.
\item On problems with more infeasible leaves, such as multi-dimensional knapsack variants, big-bucket lot-sizing and fixed charge flow, tree sizes increase with efficacious gains.
\item Note that the trees get especially large on the multi-dimensional knapsack with small right-hand sides. %strict versions of 
%random packing and covering IPs 
On these instances, we observe that on average, 67\% of the leaf nodes are pruned by infeasibility. In other words, the gap closed may not be a very good measure for the size of trees, and we should try and take branching decisions which lead to more infeasibility. 

\end{enumerate}

%This exposes the following
\begin{example}\label{ex:2}
%    Now consider the following scenario, 
Under the assumption that most nodes are pruned by infeasibility, we illustrate how using Eff-SB may pose a disadvantage. Consider two branching candidates $i$ and $j$ with $(\Delta^0_i, \Delta^1_i) = (0.7, 0.9)$ and $(\Delta^0_j, \Delta^1_j) = (0.6, \infty)$, where the infinite gain is due to an infeasible node. Here, Def-SB chooses to branch on $j$, giving it an infinite score. Now suppose that $\Delta_{p-d} = 1$, then 
\begin{align}
(q^0_i, q^1_i) = (0.7, 0.9) &\implies \text{Eff-SB score}(i) = 0.63, \\
(q^0_j, q^1_j) = (0.6, 1.0) &\implies \text{Eff-SB score}(j) = 0.6, 
\end{align}
and Eff-SB chooses to branch on $i$, 
%\red{
placing a higher value on closing the additional 0.1 gap on the child with smaller LP gain, at the cost of not pruning either of the child nodes.
%than on the other even that would have pruned a node. 
As discussed above,  pruning by bound is less common for this instance, and therefore, such decisions by Eff-SB may potentially lead to larger trees.
\end{example}

%the scenario described to motivate efficacious gains is less likely, but the above may be more frequent.
% }

In the above example, suppose we instead had a worse primal bound and $\Delta_{p-d} = 1.1$. Then the score($i$) is unaffected, but that score($j$) increases to 0.66 and becomes the preferred branching choice. This may explain the trend observed in Table~\ref{tab:effSB}, where Eff-SB produces smaller trees as the primal gap increases %why we observe in Table~\ref{tab:effSB} that Eff-SB trees of these instances get smaller as the primal gap increases. 

\subsection{Parameterizing the product score}

Our goal is to maintain the advantage that efficacious gains offers in correcting overestimation in situations similar to the 
%example described at the beginning of this section,
Example~\ref{ex:1},
while mitigating the likelihood of poor decisions in the scenario described in Example~\ref{ex:2}.
%above.

Next, we illustrate a simple strategy to improve over  Eff-SB.
\begin{example}
Instead of multiplying the efficacious gains for the two child nodes with equal exponents, we consider increasing the exponent for the node with the larger gain. As an example, %That is if the following score function is used:
consider the following score function:
\begin{eqnarray*}
\text{score}(i) =  \textup{min}\{q^0_i, q^1_i\}^{0.3}\cdot \textup{max}\{q_i^0, q_i^1\}^{0.7},
\end{eqnarray*}
Then in the case of Example~\ref{ex:2}, the score for variable $i$ is $0.83$ and that for variable $j$ is $0.86$. In other words, branching on $j$ will be selected, overcoming  
the scenario described in Example~\ref{ex:2}. At the same time, in the case of Example~\ref{ex:1}, the score of $x_j$ remains larger than $x_i$ and is still preferred for branching.
\end{example}

%With the hypothesis that this can be achieved by increasing the importance of the larger LP gain in the score function, we next endeavour to find better weights for the gains on the two branches. 

%The results of the previous experiments suggest that the product score,
Based on the above example,  we hypothesize that the product score which has been designed through rigorous testing on the actual LP gains~\citep{achterberg2007constraint}, may not be ideal for efficacious gains, and a higher exponent on the larger gain may improve the performance. 
% Due to the nature of the threshold applied to obtain efficacious gains, we may expect a higher weight on the maximum of the gains to improve the performance. \red{Why? Explain the last statement.}
% Hereafter, we focus on the product score since it is more widely used. 
We therefore modify the product score via the following parameterization, 
\begin{equation*}
    \text{score}(i) = \big(q_{i, \min}\big)^{a_{\min}} \cdot \big(q_{i, \max}\big)^{a_{\max}}
\end{equation*}
where, $q_{i, \min} = \min\{q^1_i, q^0_i\}$ and $q_{i, \max} = \max\{q^1_i, q^0_i\}$.
In order to determine good choices for $a_{\min}$ and $a_{\max}$, we conduct a parameter search with $a_{\min} \in \{0, 0.1, 0.2, \hdots, 1\}$ and $a_{\max} = 1 - a_{\min}$ and compute the tree sizes on the same randomly generated instances. 
% Moreover, we also test the robustness of the branching decisions against the quality of the primal bound by considering the scenarios where the primal gap is 2\%, 5\% and 10\%. 

\begin{figure}
    \centering
    \includegraphics[width=\linewidth]{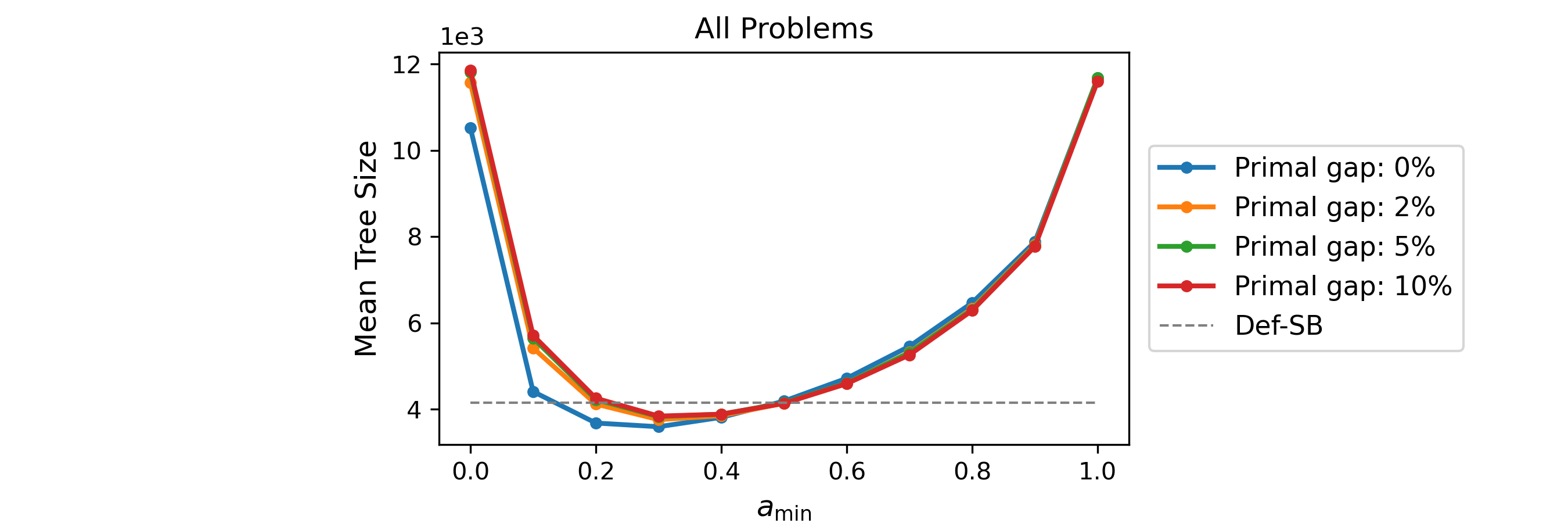}
    \caption{Geometric mean of tree sizes using efficacious gains across problem classes while varying parameter $a_{\min}$ and primal gap.}
    \label{fig:geomean_amin_primalgap}
\end{figure}

\begin{figure}
    \centering
    \includegraphics[width=\linewidth]{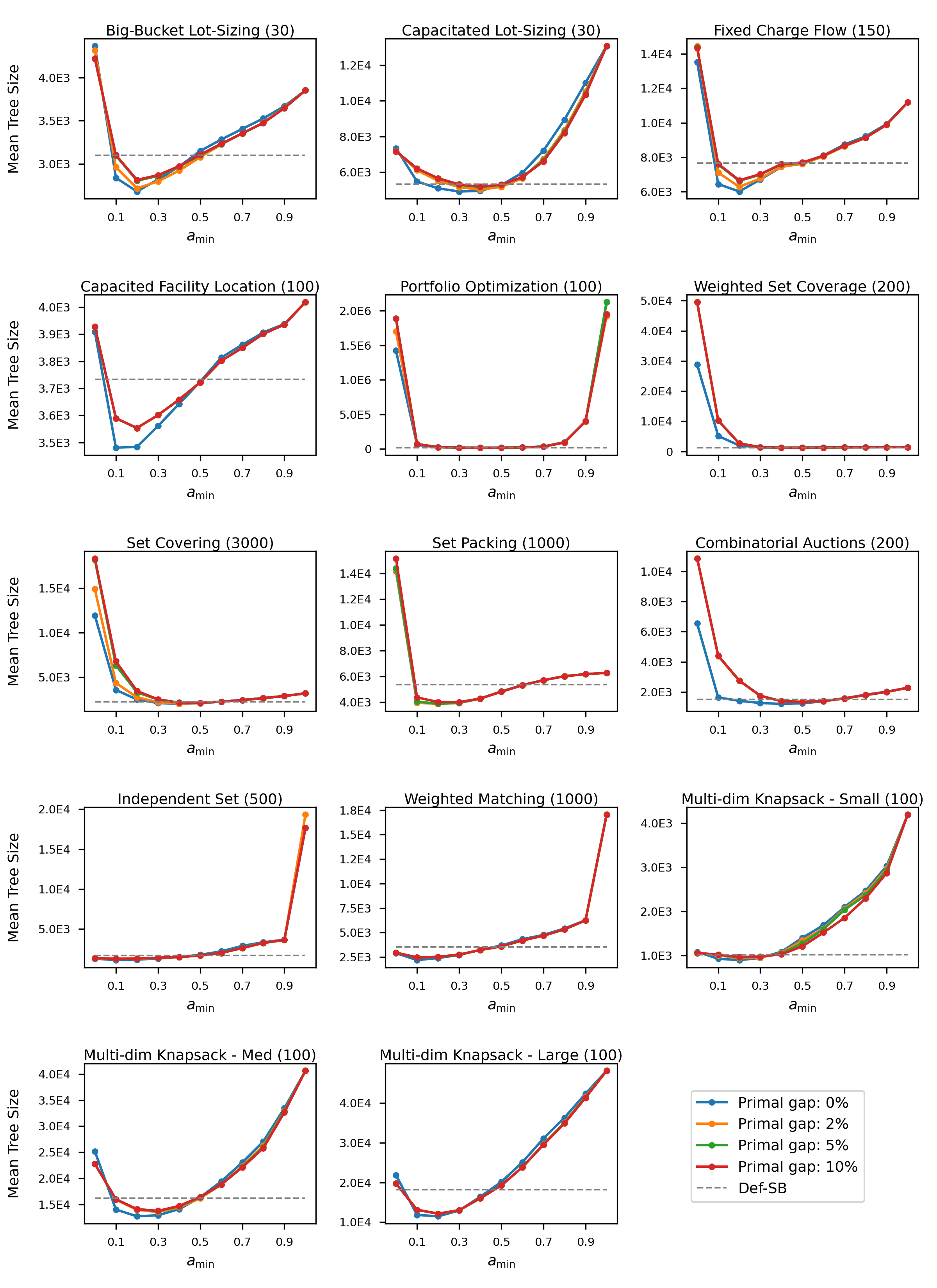}
    \caption{Impact of parameter $a_{\min}$ on the geometric mean of tree sizes using efficacious gains for different primal gaps.}
    \label{fig:amin_plots_quasiconvex}
\end{figure}

The plots in Figure~\ref{fig:amin_plots_quasiconvex} show the change in mean tree size for individual problems as $a_{\min}$ is varied. Here, each curve represents initialization of branch-and-bound tree with a different primal gap, and the dotted lines represent the tree size for Def-SB. A similar plot for geometric means across all the problem classes is shown in Figure~\ref{fig:geomean_amin_primalgap}.
Here are some key observations:
\begin{enumerate} 
\item The tree size increases with an increase in the primal gap for most values of $a_{min}$, indicating the importance of the quality of the primal bound used as the threshold. 

\item Let $a_{\min}^*$ and $a_{\max}^*$ denote the value of the parameters generating the smallest trees for any primal gap, that is, the point of minima for any curve. All plots across problem classes and primal gaps are nearly quasi-convex with $a^*_{\min} < 0.5 < a^*_{\max}$, confirming empirically our hypothesis that placing a higher importance on the larger efficacious gain improves performance.

\item As the primal gap increases, efficacious gains become less effective at eliminating over-estimation in the LP gains. Consequently, $q_{i, \max}$ is more likely to be overestimated than $q_{i, \min}$ for any variable, and increasing weight on $q_{\min}$ is advantageous. This can be observed as an increasing trend in $a^*_{\min}$ with increasing primal gap.

%\red{As hypothesized, this indicates a higher importance on the larger LP gain improves performance. } 

\end{enumerate}
% \red{The overall trend of mean tree sizes across problem types in Figure~\ref{fig:geomean_amin_primalgap} indicate that if the optimal value is known, the preferred choice of parameters should be $(a_{\min}, a_{\max}) = (0.2, 0.8)$ which has a 10.4\% reduction in tree size as compared to the original product score without primal information.}

% Moreover, when the primal gap is 0, in other words, when the optimal value is known, the smallest trees correspond to some $a^*_{\min} \leq 0.2$.  
% Parameters $(a_{\min}, a_{\max}) = (0.3, 0.7)$ are more robust to gaps in the primal bound and the preferred choice when the primal bound is sub-optimal. 
% In practice, until the termination of the branch-and-bound process, we cannot affirm the optimality of any primal bound, therefore, we focus on the latter configuration of parameters. We refer to this rule as the Eff-SB~(3,7) with score, 

Figure~\ref{fig:geomean_amin_primalgap} indicates that parameters $(a_{\min}, a_{\max}) = (0.3, 0.7)$ lead to the smallest tree sizes averaged over problems for all primal gaps. It results in smaller trees than Def-SB on most of the problems and primal gaps as seen in  Figure~\ref{fig:amin_plots_quasiconvex}. Therefore, we henceforth focus on these parameters and refer to the resulting rule as the Eff-SB~(3,7) with score, 
\begin{equation}
    \text{score}(i) = (q_{i, \min})^{0.3} \cdot (q_{i, \max})^{0.7}.
\end{equation}

We would like to highlight that while the changes in weights are applicable for the efficacious gains, using the same parameters
% Consequently, using this parameterized score 
with LP gains instead does not improve performance. That is, consider the full strong branching implementation with the following score,
\begin{equation}
    \text{score}(i) = (\Delta_{i, \min})^{0.3} \cdot (\Delta_{i, \max})^{0.7}
\end{equation}
where, $\Delta_{i, \min} = \max\{ \min\{\Delta^0_i, \Delta^1_i,\}, \epsilon \}$ and $\Delta_{i, \max} = \max\{ \max\{\Delta^0_i, \Delta^1_i,\}, \epsilon \}$. We call this rule Def-SB~(3,7). 

\begin{table}
\begin{center}
		
		\caption{Shifted geometric mean of tree sizes for Def-SB, Def-SB~(3,7) and Eff-SB~(3,7) rounded to integers. Tree sizes smaller than Def-SB are highlighted in bold. }\label{tab:effSB37}
\resizebox{\textwidth}{!}{
\setlength{\tabcolsep}{1.1em} 
  {\renewcommand{\arraystretch}{1.1}
\begin{tabular}{|l|r|r|rrrr|}
\hline
Rule &
  \multicolumn{1}{c|}{Def-SB} &
  \multicolumn{1}{c|}{Def-SB(3,7)} &
  \multicolumn{4}{c|}{Eff-SB~(3,7)} \\ \hline
Primal Gap &
  \multicolumn{1}{c|}{} &
  \multicolumn{1}{c|}{} &
  \multicolumn{1}{c}{0\%} &
  \multicolumn{1}{c}{2\%} &
  \multicolumn{1}{c}{5\%} &
  \multicolumn{1}{c|}{10\%} \\ \hline
\multicolumn{1}{|c|}{Problem}       &       &       &       &       &       &       \\ \hline
Big-Bucket Lot-Sizing (30)          & 3099  & \textbf{2909}  & \textbf{2823}  & \textbf{2801}  & \textbf{2864}  & \textbf{2872}  \\
Capacitated   Lot-Sizing (30)       & 5302  & 5371  & \textbf{4902}  & \textbf{5146}  & \textbf{5280}  & 5302  \\
Fixed Charge   Flow (150)           & 7658  & \textbf{7353}  & \textbf{6711}  & \textbf{6772}  & \textbf{6993}  & \textbf{7019}  \\
Capacited   Facility Location (100) & 3733  & \textbf{3673}  & \textbf{3561}  & \textbf{3602}  & \textbf{3602}  & \textbf{3602}  \\
Portfolio   Optimization (100)      & 19965 & 21679 & 20179 & 20619 & 20653 & 20654 \\
Weighted Set   Coverage (200)       & 1274  & 1731  & 1309  & 1415  & 1415  & 1415  \\
Set Covering   (3000)               & 2232  & 3088  & \textbf{2097}  & \textbf{2197}  & 2490  & 2509  \\
Set Packing   (1000)                & 5362  & 7302  & \textbf{3963}  & \textbf{3941}  & \textbf{3935}  & \textbf{3995}  \\
Combinatorial   Auctions (200)      & 1493  & 2607  & \textbf{1262}  & 1732  & 1744  & 1744  \\
Independent   Set (500)             & 1701  & \textbf{1430}  & \textbf{1306}  & \textbf{1377}  & \textbf{1398}  & \textbf{1401}  \\
Weighted   Matching (1000)          & 3545  & \textbf{2740}  & \textbf{2715}  & \textbf{2749}  & \textbf{2749}  & \textbf{2749}  \\
Multi-dim   Knapsack - Small (100)  & 1014  & \textbf{975}   & \textbf{940}   & \textbf{949}   & \textbf{967}   & \textbf{973}   \\
Multi-dim   Knapsack - Med (100)    & 16196 & \textbf{14283} & \textbf{12971} & \textbf{13584} & \textbf{13760} & \textbf{13805} \\
Multi-dim   Knapsack - Large (100)  & 18171 & \textbf{13082} & \textbf{12953} & \textbf{12996} & \textbf{13019} & \textbf{13023} \\ \hline
Geo Mean                            & 4149  & 4308  & \textbf{3596}  & \textbf{3761}  & \textbf{3831}  & \textbf{3843}  \\ \hline
\end{tabular}
}}
% \end{center}
% \end{table}
\vspace{24pt}
% \begin{table}
% \begin{center}
		
		\caption{Comparison of shifted geometric mean of tree sizes to evaluate Def-SB~(3,7) and Eff-SB~(3,7) rules. Tree sizes are presented for Def-SB, and the percentage reduction in tree size relative to Def-SB is presented for others.}\label{tab:effSB37_pct}
\resizebox{\textwidth}{!}{
\setlength{\tabcolsep}{1.1em} 
  {\renewcommand{\arraystretch}{1.1}
\begin{tabular}{|l|r|r|rrrr|}
\hline
Rule &
  \multicolumn{1}{c|}{Def-SB} &
  \multicolumn{1}{c|}{Def-SB(3,7)} &
  \multicolumn{4}{c|}{Eff-SB~(3,7)} \\ \hline
Primal Gap &
  \multicolumn{1}{c|}{} &
  \multicolumn{1}{c|}{} &
  \multicolumn{1}{c}{0\%} &
  \multicolumn{1}{c}{2\%} &
  \multicolumn{1}{c}{5\%} &
  \multicolumn{1}{c|}{10\%} \\ \hline
\multicolumn{1}{|c|}{Problem}       &       &       &       &       &       &       \\ \hline
Big-Bucket Lot-Sizing (30)        & 3099                        & \textbf{6\%}                     & \textbf{9\%}    & \textbf{10\%}  & \textbf{8\%}   & \textbf{7\%}   \\
Capacitated Lot-Sizing (30)       & 5302                        & -1\%                             & \textbf{8\%}    & \textbf{3\%}   & \textbf{0\%}   & 0\%            \\
Fixed Charge Flow (150)           & 7658                        & \textbf{4\%}                     & \textbf{12\%}   & \textbf{12\%}  & \textbf{9\%}   & \textbf{8\%}   \\
Capacited Facility Location (100) & 3733                        & \textbf{2\%}                     & \textbf{5\%}    & \textbf{4\%}   & \textbf{4\%}   & \textbf{4\%}   \\
Portfolio Optimization (100)      & 19965                       & -9\%                             & -1\%            & -3\%           & -3\%           & -3\%           \\
Weighted Set Coverage (200)       & 1274                        & -36\%                            & -3\%            & -11\%          & -11\%          & -11\%          \\
Set Covering (3000)               & 2232                        & -38\%                            & \textbf{6\%}    & \textbf{2\%}   & -12\%          & -12\%          \\
Set Packing (1000)                & 5362                        & -36\%                            & \textbf{26\%}   & \textbf{27\%}  & \textbf{27\%}  & \textbf{26\%}  \\
Combinatorial Auctions (200)      & 1493                        & -75\%                            & \textbf{15\%}   & -16\%          & -17\%          & -17\%          \\
Independent Set (500)             & 1701                        & \textbf{16\%}                    & \textbf{23\%}   & \textbf{19\%}  & \textbf{18\%}  & \textbf{18\%}  \\
Weighted Matching (1000)          & 3545                        & \textbf{23\%}                    & \textbf{23\%}   & \textbf{22\%}  & \textbf{22\%}  & \textbf{22\%}  \\
Multi-dim Knapsack - Small (100)  & 1014                        & \textbf{4\%}                     & \textbf{7\%}    & \textbf{6\%}   & \textbf{5\%}   & \textbf{4\%}   \\
Multi-dim Knapsack - Med (100)    & 16196                       & \textbf{12\%}                    & \textbf{20\%}   & \textbf{16\%}  & \textbf{15\%}  & \textbf{15\%}  \\
Multi-dim Knapsack - Large (100)  & 18171                       & \textbf{28\%}                    & \textbf{29\%}   & \textbf{28\%}  & \textbf{28\%}  & \textbf{28\%}  \\ \hline
Geo Mean                          & 4149                        & -3.8\%                           & \textbf{13.3\%} & \textbf{9.4\%} & \textbf{7.7\%} & \textbf{7.4\%} \\ \hline

\end{tabular}
}}
\end{center}
\end{table}

Table \ref{tab:effSB37} compares the performance of Eff-SB~(3,7) and Def-SB~(3,7) against Def-SB in terms of tree sizes, and Table~\ref{tab:effSB37_pct} presents the same results in terms of percentage reduction in tree sizes relative to Def-SB. 
\begin{enumerate}
    \item Eff-SB~(3,7) demonstrates a 13.3\% reduction in mean tree sizes over problem classes when provided with the optimal primal value and a 9.4\%, 7.7\% and 7.4\% reduction when the primal gap is 2\%, 5\% and 10\%, respectively. 
    \item On the other hand, the mean tree size of Def-SB~(3,7) is 3.8\% larger than that of Def-SB.
    This is because a higher emphasis on the larger LP gain can be detrimental without correcting for the overestimation, and may lead to much larger trees as seen in the case of combinatorial auctions, set packing, set covering, and weighted set coverage.
    \item {It is also remarkable that when the primal bound is optimal, Eff-SB~(3, 7) yields smaller trees on average for all problems except portfolio optimization and weighted set coverage. On these two problems, the mean tree size increases by 1\%-3\% only. } With suboptimal bounds Eff~(3,7) leads to larger trees on four out of the fourteen problems.
\end{enumerate}

While parameters $(a_{\min}, a_{\max}) = (0.3, 0.7)$ lead to smaller trees than Def-SB in most cases, the improvements are even greater if the best parameters are known for the problem. The results obtained by picking the best parameter for each problem are shown in Table~\ref{tab:opteffSB} as the shifted geometric mean of tree sizes and in Table~\ref{tab:opteffSB_pct} as the percentage reduction of the same relative to Def-SB. The optimal parameters lead to trees smaller than Def-SB on all problems and on average, further reduce the tree sizes by 3.5\%-5\% than Eff-SB~(3,7).

Note that the best parameters have a strictly larger value for $a^*_{\max}$ than $a^*_{\min}$ for all problems except one, where they are equal.
Another interpretation for this may be that without overestimation it is advantageous to emphasize more on the maximum gain to ensure that sufficient progress is made on at least one branch. But with overestimation, the maximum may be unreliable and therefore popular FSB scores (that product rule or the  $(5\Delta_{\min} + \Delta_{\max})$ rule~\cite{achterberg2007constraint, achterberg2005branching}) are designed to be more sensitive to the minimum.

\begin{table}
\begin{center}	
\caption{Shifted geometric mean of tree sizes for Eff-SB with the best choice of $(a_{min}, a_{max})$ rounded to integer values. }\label{tab:opteffSB}
\resizebox{\textwidth}{!}{
\setlength{\tabcolsep}{1.2em} 
  {\renewcommand{\arraystretch}{1.1}
\begin{tabular}{|l|c|r|rrrr|}
\hline
Rule &  & \multicolumn{1}{c|}{Def-SB} & \multicolumn{4}{c|}{Eff SB - $(a^*_{\min}, a^*_{\max})$} \\ \hline
Primal Gap &  & \multicolumn{1}{c|}{} & \multicolumn{1}{c}{0\%} & \multicolumn{1}{c}{2\%} & \multicolumn{1}{c}{5\%} & \multicolumn{1}{c|}{10\%} \\ \hline
\multicolumn{1}{|c|}{Problem} & $(a^*_{\min}, a^*_{\max})$ & \multicolumn{1}{r|}{} & \multicolumn{1}{r}{} & \multicolumn{1}{r}{} & \multicolumn{1}{r}{} & \multicolumn{1}{r|}{} \\ \hline
Big-Bucket Lot-Sizing (30) & (0.2, 0.8)& 3099 & \textbf{2683} & \textbf{2720} & \textbf{2808} & \textbf{2818} \\
Capacitated Lot-Sizing (30) & (0.4, 0.6)& 5302 & \textbf{4938} & \textbf{5008} & \textbf{5147} & \textbf{5171} \\
Fixed Charge Flow (150) & (0.2, 0.8)& 7658 & \textbf{6014} & \textbf{6284} & \textbf{6626} & \textbf{6666} \\
Capacited Facility Location (100) & (0.2, 0.8)& 3733 & \textbf{3484} & \textbf{3554} & \textbf{3554} & \textbf{3554} \\
Portfolio Optimization (100) & (0.4, 0.6)& 19965 & \textbf{19380} & \textbf{19427} & \textbf{19448} & \textbf{19449} \\
Weighted Set Coverage (200) & (0.4, 0.6)& 1274 & \textbf{1227} & \textbf{1246} & \textbf{1246} & \textbf{1246} \\
Set Covering (3000) & (0.4, 0.6)& 2232 & \textbf{2001} & \textbf{2034} & \textbf{2135} & \textbf{2141} \\
Set Packing (1000) & (0.2, 0.8)& 5362 & \textbf{3883} & \textbf{3856} & \textbf{3863} & \textbf{4001} \\
Combinatorial Auctions (200) & (0.5, 0.5)& 1493 & \textbf{1245} & \textbf{1335} & \textbf{1337} & \textbf{1337} \\
Independent Set (500) & (0.1, 0.9)& 1701 & \textbf{1105} & \textbf{1272} & \textbf{1282} & \textbf{1279} \\
Weighted Matching (1000) & (0.1, 0.9)& 3545 & \textbf{2185} & \textbf{2464} & \textbf{2464} & \textbf{2464} \\
Multi-dim Knapsack - Small (100) & (0.2, 0.8)& 1014 & \textbf{890} & \textbf{934} & \textbf{953} & \textbf{959} \\
Multi-dim Knapsack - Med (100) & (0.3, 0.7)& 16196 & \textbf{12971} & \textbf{13584} & \textbf{13760} & \textbf{13805} \\
Multi-dim Knapsack - Large (100) & (0.2, 0.8)& 18171 & \textbf{11467} & \textbf{12126} & \textbf{12137} & \textbf{12137} \\ \hline
Geo Mean & \multicolumn{1}{l|}{} & 4149 & \textbf{3371} & \textbf{3515} & \textbf{3567} & \textbf{3582} \\ \hline
\end{tabular}
}}
% \end{center}
% \end{table}

\vspace{24pt}

% \begin{table}
% \begin{center}	
\caption{Shifted geometric mean of tree sizes for Eff-SB with the best choice of $(a_{min}, a_{max})$ presented as percentage reduction relative to Def-SB.}\label{tab:opteffSB_pct}
\resizebox{\textwidth}{!}{
\setlength{\tabcolsep}{1.2em} 
  {\renewcommand{\arraystretch}{1.1}
\begin{tabular}{|l|c|r|rrrr|}
\hline
Rule &  & \multicolumn{1}{c|}{Def-SB} & \multicolumn{4}{c|}{Eff SB - $(a^*_{\min}, a^*_{\max})$} \\ \hline
Primal Gap &  & \multicolumn{1}{c|}{} & \multicolumn{1}{c}{0\%} & \multicolumn{1}{c}{2\%} & \multicolumn{1}{c}{5\%} & \multicolumn{1}{c|}{10\%} \\ \hline
\multicolumn{1}{|c|}{Problem} & $(a^*_{\min}, a^*_{\max})$ & \multicolumn{1}{r|}{} & \multicolumn{1}{r}{} & \multicolumn{1}{r}{} & \multicolumn{1}{r}{} & \multicolumn{1}{r|}{} \\ \hline
Big-Bucket Lot-Sizing (30) & (0.2, 0.8)& 3099 & \textbf{13\%} & \textbf{12\%} & \textbf{9\%} & \textbf{9\%} \\
Capacitated Lot-Sizing (30) & (0.4, 0.6)& 5302 & \textbf{7\%} & \textbf{6\%} & \textbf{3\%} & \textbf{2\%} \\
Fixed Charge Flow (150) & (0.2, 0.8)& 7658 & \textbf{21\%} & \textbf{18\%} & \textbf{13\%} & \textbf{13\%} \\
Capacited Facility Location (100) & (0.2, 0.8)& 3733 & \textbf{7\%} & \textbf{5\%} & \textbf{5\%} & \textbf{5\%} \\
Portfolio Optimization (100) & (0.4, 0.6)& 19965 & \textbf{3\%} & \textbf{3\%} & \textbf{3\%} & \textbf{3\%} \\
Weighted Set Coverage (200) & (0.4, 0.6)& 1274 & \textbf{4\%} & \textbf{2\%} & \textbf{2\%} & \textbf{2\%} \\
Set Covering (3000) & (0.4, 0.6)& 2232 & \textbf{10\%} & \textbf{9\%} & \textbf{4\%} & \textbf{4\%} \\
Set Packing (1000) & (0.2, 0.8)& 5362 & \textbf{28\%} & \textbf{28\%} & \textbf{28\%} & \textbf{25\%} \\
Combinatorial Auctions (200) & (0.5, 0.5)& 1493 & \textbf{17\%} & \textbf{11\%} & \textbf{10\%} & \textbf{10\%} \\
Independent Set (500) & (0.1, 0.9)& 1701 & \textbf{35\%} & \textbf{25\%} & \textbf{25\%} & \textbf{25\%} \\
Weighted Matching (1000) & (0.1, 0.9)& 3545 & \textbf{38\%} & \textbf{30\%} & \textbf{30\%} & \textbf{30\%} \\
Multi-dim Knapsack - Small (100) & (0.2, 0.8)& 1014 & \textbf{12\%} & \textbf{8\%} & \textbf{6\%} & \textbf{5\%} \\
Multi-dim Knapsack - Med (100) & (0.3, 0.7)& 16196 & \textbf{20\%} & \textbf{16\%} & \textbf{15\%} & \textbf{15\%} \\
Multi-dim Knapsack - Large (100) & (0.2, 0.8)& 18171 & \textbf{37\%} & \textbf{33\%} & \textbf{33\%} & \textbf{33\%} \\ \hline
Geo Mean & \multicolumn{1}{l|}{} & 4149 & \textbf{18.8\%} & \textbf{15.3\%} & \textbf{14.0\%} & \textbf{13.7\%} \\ \hline
\end{tabular}
}}
\end{center}
\end{table}

\begin{remark}

    Some solvers might be implicitly using efficacious gains since LPs at nodes are re-solved using dual simplex. However, our goal here is to understand the properties of FSB, and to the best of our knowledge, 
     this work is the first systematic study of the impact of using efficacious gain on the size of branch-and-bound trees. The result of this section show that if a solver is implementing efficacious gains for FSB, then it may be advisable to consider reducing $a_{\min}$ to $0.3$.

\end{remark}

\begin{remark}
    We also consider an alternative approach to correcting overestimation error that treats pruning by bound as equivalent to pruning by infeasibility, thereby preferring to branch on any variable that leads to a pruned child by assigning it an infinite score. While this rule performs very well when the incumbent primal bound is optimal, its performance quickly deteriorates with suboptimal solutions. A detailed examination of this rule is presented in Appendix~\ref{sec:Greedy pruning}.
\end{remark}

\section{Factoring infeasibility and integrality}\label{sec:la}

Full strong branching selects a branching variable only based on the improvement in local LP gains, greedily racing towards pruning by bounds. 
% It cannot foresee the impact of a branching decision on nodes beyond the immediate children. 
In particular, it does not take into account the changes in feasible region or the changes in fractionality of the LP solution and therefore cannot predict the occurrence of integrality or infeasibility beyond the immediate children. 
Such long-range impact of branching decisions is difficult to model or predict. Moreover, its effects are often asymmetric on the 0 and 1 branches depending on the problem structure. For example, consider a problem that is dominated by packing constraints. Here, fixing too many variables to 1 would render the problem infeasible, however, fixing variables to 0 will maintain feasibility.
% \margin{Should we show box whisker plots for random problems to substantiate asymmetry?}
Thus, at any node, we may expect the sub-tree size on the branch with the variable fixed to 1 to be limited to some extent due to pruning by infeasibility and find it advantageous to give greater importance to the LP gain on the branch fixing the variable 0. Similarly, for problems with covering constraints, a higher importance on the LP gain for the branch fixing 1 may be more beneficial. 

We propose modifying the strong branching score to distinctly consider the LP gains on 0 and 1 branches to account for such differences. 
Building on the optimized score from the previous section, we consider the following generalization of the score, 
\begin{equation}
    \text{score}(i) = \big(q_{i, 0}\big)^{a_{0}} \cdot \big(q_{i, 1}\big)^{a_{1}} \cdot \big(q_{i, \min}\big)^{0.3} \cdot \big(q_{i, \max}\big)^{0.7}. \label{eq:asymmetrical_score}
\end{equation}
where only one of $a_0$ and $a_1$ may be positive, corresponding to the branching with a higher need of on closing the gap. 
\subsection{Inferring asymmetry from partially solved trees}
We propose to choose values of $a_0$ and $a_1$ by considering global information from partially solved branch-and-bound trees.
The key idea behind this approach is the following. To infer the direction of asymmetry, for every infeasible leaf node, we would like to % need to 
identify the subset of branching constraints that produces infeasibility.
%essential for infeasibility from others that may not affect it. 
\textcolor{black}{This is important as several constraints may be incidental due to the path of the node. Considering all constraints may not indicate the correct asymmetrical trend in infeasibility as illustrated by the example in Appendix~\ref{sec:necc_branching_constrs}.}

Identifying any minimally infeasible subsystem of constraints involves solving a related LP~\citep{gleeson1990identifying}. 
However, note that the last branching constraint was necessary for rendering the LP at the node to be infeasible and, therefore, must be present in any minimally infeasible subsystem. 
Moreover, as we are only interested in ascertaining the direction of asymmetry rather than identifying the contributing variables, it is sufficient to consider only the last branch on the path to each infeasible node.
% we only consider the last branch on the path to each infeasible node for a quicker estimate.  
Note that if all variables have the same coefficient sign, that is, the problem is either a packing-type or a covering-type, then all branching constraints in any minimal infeasible subsystem must fix variables to the same value as the last branch. Therefore, to choose a value of $a_0$ and $a_1$, we track the number of times the last fixing before an infeasible node was $0$ or $1$ in a partial tree. Henceforth, we refer to collecting this information as the \textit{last assignment} rule. 

% Besides infeasible leaf nodes, we would also like to consider asymmetry due to integral nodes 
%for two reasons.
Besides infeasible leaf nodes, we would also like to consider and correct for the asymmetry due to integral leaf nodes. 
%integral leaf nodes because integrality may also be an asymmetric trend due to the problem structure similar to infeasibility.
%in this approach for two primary reasons. 
%First, 
Observe that in set-based, cardinality-constrained or similar problems where the constraints only have 0-1 coefficients and integral right-hand side, infeasibility never occurs. For instance, in set-packing, if a variable is fixed to 1, all other variables appearing in the same constraints have a value of 0 in any feasible LP solution. Since they are never fractional, they are never branched on and therefore cannot take a positive value to cause infeasibility. Instead, fixing enough variables to 1 may quickly lead to an integral solution. 
The last assignment rule when applied to integral nodes is also helpful in detecting the asymmetry caused by integrality. To see this, consider two packing instances one with very small right-hand-sides and another with large right-hand-sides. For the former, we expect a large number of infeasible nodes due to a few variables being fixed to $1$ -- this information is captured by the last assignment rule as discussed in previous paragraph -- and we adjust $a_0$ to put more emphasis on the $0$ branching side. On the other hand, for the second type of instance, fixing a few variables to $0$ may lead to an integral node, again an information captured by the last assignment rule, leading to an adjustment of $a_1$ to put more emphasis on the $1$ branching side. 

%integrality may also be an asymmetric trend due to the problem structure similar to infeasibility. In a knapsack problem if the right-hand side is large enough, fixing only a few variables to 0 may make it possible to assign all others to 1, quickly leading to an integral solution.
% In a knapsack problem with a small right-hand side, fixing a few variables to 1 may soon violate the constraints and lead to infeasibility. On the other hand, if the right-hand side was large enough, fixing only a few variables to 0 may make it possible to assign all others to 1, leading to an integral solution.

\begin{algorithm}
\caption{Updating strong branching score function based on the last assignment rule}  \label{algo:last_assignment}
\begin{algorithmic}[1]
\State \textbf{Input:} $\mathcal{L}_I$, $l(N)$ for all $N \in \mathcal{L}_I$
\State \textbf{Output:} Exponents $(a_0, a_1)$ for score function in Eq.~\eqref{eq:asymmetrical_score}
\State \textbf{Initialize:} $a_0 \gets 0$, $a_1 \gets 0$
\If{$|\mathcal{L}_I| \geq k_I$ }
    \State $n_0 \gets \big|\{N \in \mathcal{L}_I : l(N) = 0\}\big|$,   $\;n_1 \gets \big|\{N \in \mathcal{L}_I : l(N) = 1\}\big|$
    \State $a \gets \frac{n_0 - n_1}{n_0 + n_1}$
    \If{$a > 0$}
        \State $a_1 \gets \eta \times a$
    \Else
        \State $a_0 \gets -\eta \times a$
    \EndIf
\EndIf
\State \Return $(a_0, a_1)$
\end{algorithmic}
\end{algorithm}

The detailed procedure to update the score function based on this idea is presented in Algorithm~\ref{algo:last_assignment} and is applied at a regular frequency within the branch-and-bound framework.
Let $\mathcal{L}_I$ be the set of all integral or infeasible leaf nodes. 
For every node $N \in \mathcal{L}_I$, we track whether $N$ was created by fixing the branching variable to 0 or 1, denoted by the function $l(\cdot): \mathcal{L}_I \longrightarrow \{0, 1\}$.
% Intuitively, $l(N) = 0$ for most nodes indicates a higher occurrence of pruning due to infeasibility and integrality on the 0-branch and an increased importance of closing the gap on the 1-branch. 
The number of nodes with $l(N)$ equal to 0 and 1 are represented by $n_0$ and $n_1$ respectively. Their difference is normalized by the total number of nodes, multiplied by a constant parameter $\eta$ for appropriate scaling and attributed as an exponent to the branch with fewer pruning due to integrality and infeasibility. Finally, note that this algorithm is employed only if there is a large enough sample of nodes in $\mathcal{L}_I$ to ascertain any asymmetry in pruning between the two branches and is indicated by the threshold $k_I$.

Other approaches in literature have also leveraged branching constraints causing infeasibility to improve branching decisions. Hybrid branching~\citep{achterberg2009hybrid} combines conflict clauses used to select branching variables in SAT solvers, and the fathoming-based approach in \citep{kilincc2009information} is designed to obtain information from nodes pruned due to infeasibility as well as bounds. 
These methods use this information collected at the variable level to generate the relative importance of branching on a variable. 
%compare the importance of individual branching variable candidates.
On the other hand, our approach fundamentally differs in that we use this information to infer global asymmetry trends due to infeasibility and integrality in the branch-and-bound tree, and do not use it to make distinctions amongst variables.

\subsection{Experiments on last assignment rule}

We now want to test the score generated by the last assignment rule and understand the impact of parameter $\eta$ on tree sizes. We conduct computational experiments similar to Section~\ref{sec:primal}, using the same test bed of random problems and initialize the branch-and-bound process with feasible solutions at primal gaps of  0\%, 5\% and 10\%. We use $k_I = 10$ and vary $\eta \in \{ 0.05, 0.1, 0.15, 0.2, 0.25\}$. \textcolor{black}{This range is selected so that the exponent on the asymmetrical term is smaller than the exponent for $q_{min}$}. \textcolor{black}{Algorithm~\ref{algo:last_assignment} is called to update the score after every 50 nodes within the branch-and-bound process. }

\begin{figure}
    \centering
    \includegraphics[width=\linewidth]{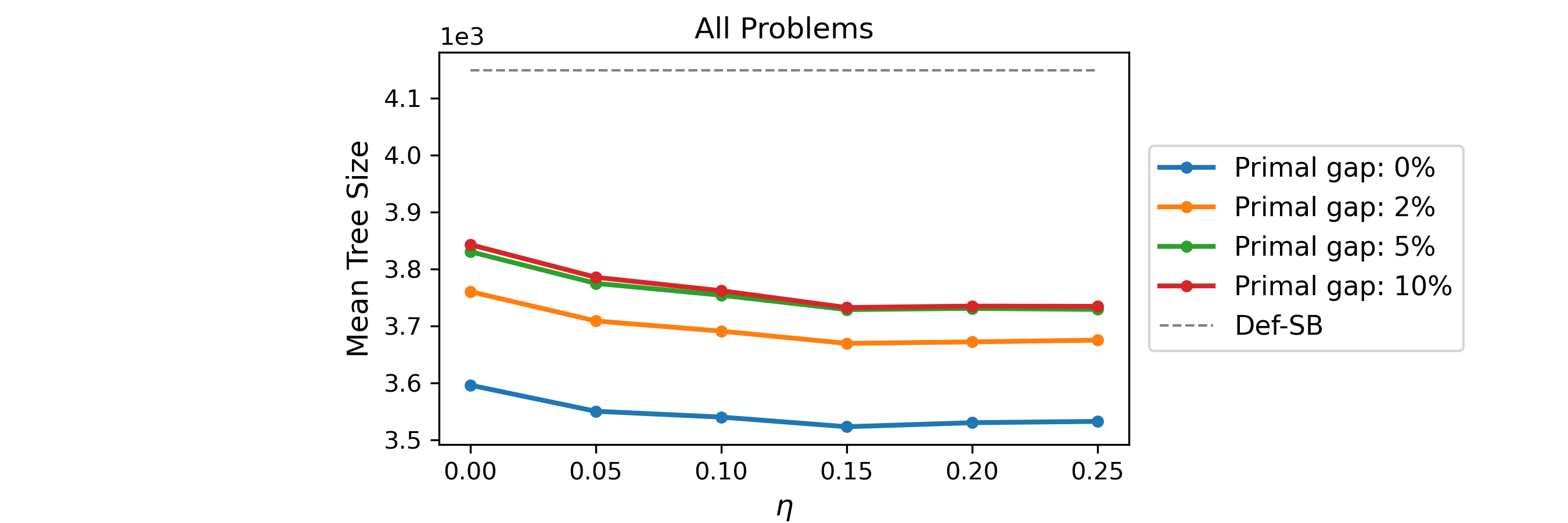}
    \caption{Geometric mean of tree sizes across problem classes when using the last assignment rule and varying parameter $\eta$ and primal gap.}
    \label{fig:geomean_eta_primalgap}
\end{figure}

\begin{figure}
    \centering
    \includegraphics[width=\linewidth]{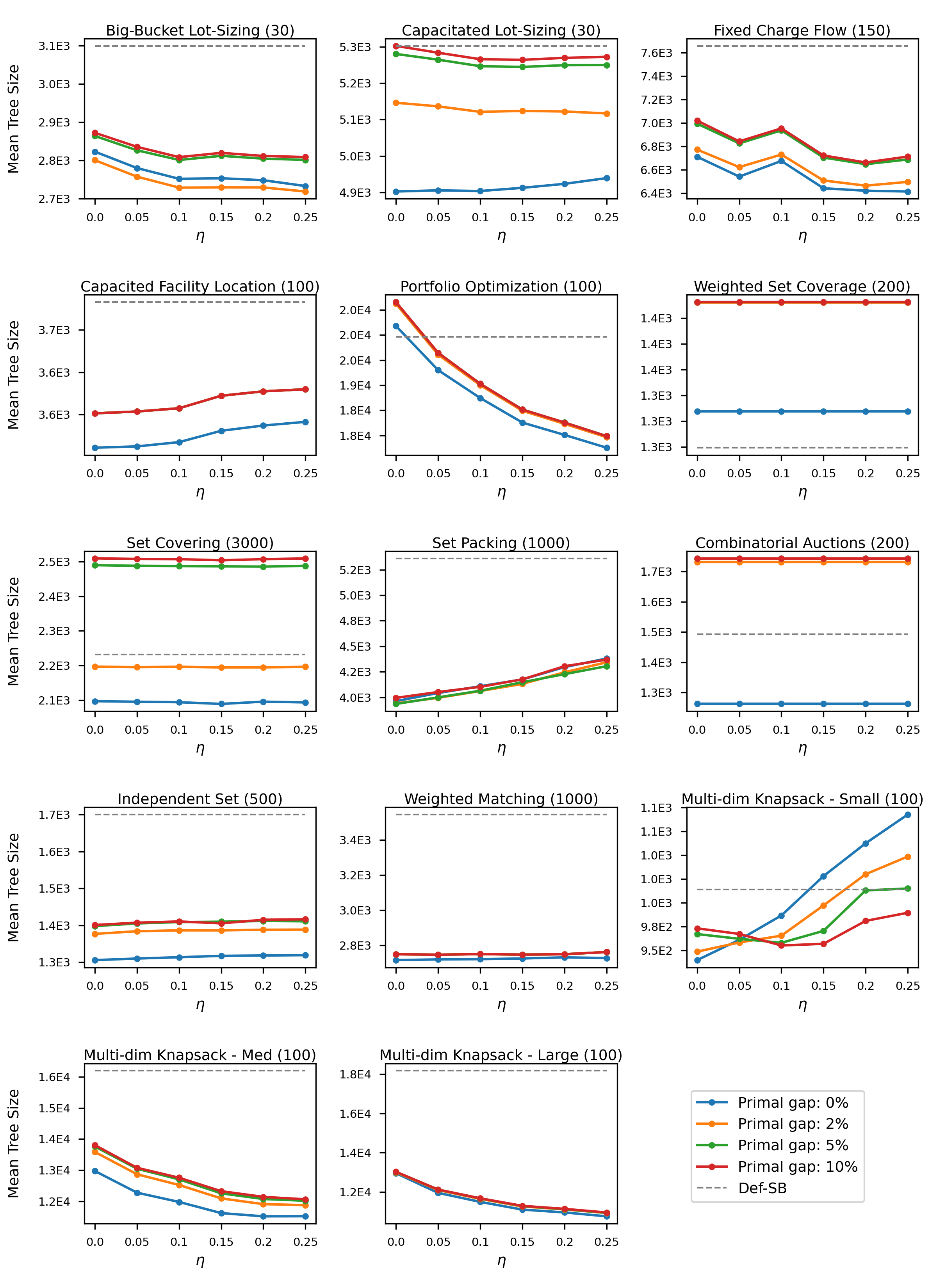}
    \caption{Impact of parameter $\eta$ on the geometric mean of tree sizes using the last assignment rule for different primal gaps.}
    \label{fig:eta_plots_quasiconvex}
\end{figure}

Figure~\ref{fig:eta_plots_quasiconvex} presents the geometric mean of tree sizes for individual problems against different choices of $\eta$, and Figure~\ref{fig:geomean_eta_primalgap} presents the geometric mean of tree sizes aggregated across all problems. Note that $\eta=0$ implies that the last assignment rule was not used and corresponds to the Eff-SB~(3,7) score.

The results in Figure~\ref{fig:eta_plots_quasiconvex} indicate very different trends for different problems. In particular:
\begin{enumerate}
    \item Several problems show an overall decreasing trend with an increase in $\eta$ throughout the range of values considered in the experiment. These include big-bucket lot-sizing, portfolio optimization, and multi-dimensional knapsack with medium and large right-hand sides.
    % random packing and covering problems, except for the ``strict'' variants.
    \item Multi-dimensional knapsack problem with small right-hand side and capacitated lot-sizing have an initial decrease and then an increasing trend.
\item  On the other hand, set-packing and capacitated facility location problems see an almost steady increase.
\item Some problems, such as combinatorial auctions, independent set, weighted matching and weighted set coverage, do not have a substantial number of integer or infeasible leaf nodes. Therefore, the last assignment rule is less relevant, and these problems are not significantly impacted by changes in $\eta$.
\end{enumerate}

\begin{table}
\begin{center}
		
		\caption{Evaluating LA-SB~(3,7,0.15) based on shifted geometric mean of tree sizes rounded to integers. Smaller tree sizes between Eff-SB(3,7) and LA-SB(3,7,0.15) are highlighted in bold.}\label{tab:la_0.15}
\resizebox{1.\textwidth}{!}{
\setlength{\tabcolsep}{0.5em} 
  {\renewcommand{\arraystretch}{1.1}
\begin{tabular}{|l|c|r|rrrr|rrrr|}
\hline
Rule                              &                                       & \multicolumn{1}{c|}{Def-SB} & \multicolumn{4}{c|}{Eff-SB (3,7)}                                                                       & \multicolumn{4}{c|}{LA-SB (3,7,0.15)}                                                                   \\ \hline
Primal Gap                        &                                       & \multicolumn{1}{c|}{}       & \multicolumn{1}{c}{0\%} & \multicolumn{1}{c}{2\%} & \multicolumn{1}{c}{5\%} & \multicolumn{1}{c|}{10\%} & \multicolumn{1}{c}{0\%} & \multicolumn{1}{c}{2\%} & \multicolumn{1}{c}{5\%} & \multicolumn{1}{c|}{10\%} \\ \hline
\multicolumn{1}{|c|}{Problem}     & $|\mathcal{L}_I|/|\mathcal{L}|$ &                             &                         &                         &                         &                           &                         &                         &                         &                           \\ \hline
Big-Bucket Lot-Sizing (30)        & 7.3\%                                 & 3099                        & 2823                    & 2801                    & 2864                    & 2872                      & \textbf{2753}           & \textbf{2730}           & \textbf{2812}           & \textbf{2820}             \\
Capacitated Lot-Sizing (30)       & 1.9\%                                 & 5302                        & \textbf{4902}           & 5146                    & 5280                    & 5302                      & 4913                    & \textbf{5124}           & \textbf{5245}           & \textbf{5264}             \\
Fixed Charge Flow (150)           & 8.0\%                                 & 7658                        & 6711                    & 6772                    & 6993                    & 7019                      & \textbf{6442}           & \textbf{6508}           & \textbf{6704}           & \textbf{6722}             \\
Capacited Facility Location (100) & 1.2\%                                 & 3733                        & \textbf{3561}           & \textbf{3602}           & \textbf{3602}           & \textbf{3602}             & 3581                    & 3622                    & 3622                    & 3622                      \\
Portfolio Optimization (100)      & 19.6\%                                & 19965                       & 20179                   & 20619                   & 20653                   & 20654                     & \textbf{18257}          & \textbf{18491}          & \textbf{18516}          & \textbf{18516}            \\
Weighted Set Coverage (200)       & 0.5\%                                 & 1274                        & 1309                    & 1415                    & 1415                    & 1415                      & 1309                    & 1415                    & 1415                    & 1415                      \\
Set Covering (3000)               & 0.6\%                                 & 2232                        & 2097                    & 2197                    & 2490                    & 2509                      & \textbf{2090}           & \textbf{2195}           & \textbf{2486}           & \textbf{2504}             \\
Set Packing (1000)                & 6.1\%                                 & 5362                        & \textbf{3963}           & \textbf{3941}           & \textbf{3935}           & \textbf{3995}             & 4175                    & 4131                    & 4148                    & 4176                      \\
Combinatorial Auctions (200)      & 0.3\%                                 & 1493                        & 1262                    & 1732                    & 1744                    & 1744                      & 1262                    & 1732                    & 1744                    & 1744                      \\
Independent Set (500)             & 3.5\%                                 & 1701                        & \textbf{1306}           & \textbf{1377}           & \textbf{1398}           & \textbf{1401}             & 1317                    & 1386                    & 1410                    & 1405                      \\
Weighted Matching (1000)          & 1.4\%                                 & 3545                        & \textbf{2715}           & 2749                    & 2749                    & 2749                      & 2725                    & \textbf{2747}           & \textbf{2747}           & \textbf{2747}             \\
Multi-dim Knapsack - Small (100)  & 67\%                                  & 1014                        & \textbf{940}            & \textbf{949}            & \textbf{967}            & 973                       & 1028                    & 997                     & 971                     & \textbf{957}              \\
Multi-dim Knapsack - Med (100)    & 40.1\%                                & 16196                       & 12971                   & 13584                   & 13760                   & 13805                     & \textbf{11617}          & \textbf{12090}          & \textbf{12253}          & \textbf{12322}            \\
Multi-dim Knapsack - Large (100)  & 23.8\%                                & 18171                       & 12953                   & 12996                   & 13019                   & 13023                     & \textbf{11098}          & \textbf{11257}          & \textbf{11277}          & \textbf{11289}            \\ \hline
Geo Mean                          &                                       & 4149                        & 3596                    & 3761                    & 3831                    & 3843                      & \textbf{3524}           & \textbf{3670}           & \textbf{3729}           & \textbf{3733}             \\ \hline
\end{tabular}
}}
% \end{center}
% \end{table}
\vspace{24pt}

% \begin{table}
% \begin{center}		
		\caption{
        Comparison of the shifted geometric mean of tree sizes to evaluate LA-SB~(3,7,0.15). Tree sizes are presented for Def-SB, and the percentage reduction in tree size relative to Def-SB is presented for others.}\label{tab:la_0.15_pct}
\resizebox{1.\textwidth}{!}{
\setlength{\tabcolsep}{0.5em} 
  {\renewcommand{\arraystretch}{1.2}
% Please add the following required packages to your document preamble:
% \usepackage[table,xcdraw]{xcolor}
% Beamer presentation requires \usepackage{colortbl} instead of \usepackage[table,xcdraw]{xcolor}
\begin{tabular}{|l|c|r|rrrr|rrrr|}
\hline
Rule                              &        & \multicolumn{1}{c|}{Def-SB} & \multicolumn{4}{c|}{Eff-SB (3,7)}                                                                                                                 & \multicolumn{4}{c|}{LA-SB (3,7,0.15)}                                                                   \\ \hline
Primal Gap                        &        & \multicolumn{1}{c|}{}       & \multicolumn{1}{c}{0\%}            & \multicolumn{1}{c}{2\%}            & \multicolumn{1}{c}{5\%}            & \multicolumn{1}{c|}{10\%}          & \multicolumn{1}{c}{0\%} & \multicolumn{1}{c}{2\%} & \multicolumn{1}{c}{5\%} & \multicolumn{1}{c|}{10\%} \\ \hline
\multicolumn{1}{|c|}{Problem}     &    $|\mathcal{L}_I|/|\mathcal{L}|$   &                             &                                    &                                    &                                    &                                    &                         &                         &                         &                           \\ \hline
Big-Bucket Lot-Sizing (30)        & 7.3\%  & 3099                        & 9\%                                & 10\%                               & 8\%                                & 7\%                                & \textbf{11\%}           & \textbf{12\%}           & \textbf{9\%}            & \textbf{9\%}              \\
Capacitated Lot-Sizing (30)       & 1.9\%  & 5302                        & \textbf{8\%}                       & 3\%                                & 0\%                                & 0\%                                & 7\%                     & \textbf{3\%}            & \textbf{1\%}            & \textbf{1\%}              \\
Fixed Charge Flow (150)           & 8.0\%  & 7658                        & 12\%                               & 12\%                               & 9\%                                & 8\%                                & \textbf{16\%}           & \textbf{15\%}           & \textbf{12\%}           & \textbf{12\%}             \\
Capacited Facility Location (100) & 1.2\%  & 3733                        & \textbf{5\%}                       & \textbf{4\%}                       & \textbf{4\%}                       & \textbf{4\%}                       & 4\%                     & 3\%                     & 3\%                     & 3\%                       \\
Portfolio Optimization (100)      & 19.6\% & 19965                       & -1\%                               & -3\%                               & -3\%                               & -3\%                               & \textbf{9\%}            & \textbf{7\%}            & \textbf{7\%}            & \textbf{7\%}              \\
Weighted Set Coverage (200)       & 0.5\%  & 1274                        & -3\%                               & -11\%                              & -11\%                              & -11\%                              & -3\%                    & -11\%                   & -11\%                   & -11\%                     \\
Set Covering (3000)               & 0.6\%  & 2232                        & 6\%                                & 2\%                                & -12\%                              & -12\%                              & \textbf{6\%}            & \textbf{2\%}            & \textbf{-11\%}          & \textbf{-12\%}            \\
Set Packing (1000)                & 6.1\%  & 5362                        & $\,$\textbf{26\%} & $\,$\textbf{27\%} & $\,$\textbf{27\%} & $\,$\textbf{26\%} & 22\%                    & 23\%                    & 23\%                    & 22\%                      \\
Combinatorial Auctions (200)      & 0.3\%  & 1493                        & 15\%                               & -16\%                              & -17\%                              & -17\%                              & 15\%                    & -16\%                   & -17\%                   & -17\%                     \\
Independent Set (500)             & 3.5\%  & 1701                        & \textbf{23\%}                      & \textbf{19\%}                      & \textbf{18\%}                      & \textbf{18\%}                      & 23\%                    & 18\%                    & 17\%                    & 17\%                      \\
Weighted Matching (1000)          & 1.4\%  & 3545                        & \textbf{23\%}                      & 22\%                               & 22\%                               & 22\%                               & 23\%                    & \textbf{23\%}           & \textbf{23\%}           & \textbf{23\%}             \\
Multi-dim Knapsack - Small (100)  & 67\%   & 1014                        & \textbf{7\%}                       & \textbf{6\%}                       & \textbf{5\%}                       & 4\%                                & -1\%                    & 2\%                     & 4\%                     & \textbf{6\%}              \\
Multi-dim Knapsack - Med (100)    & 40.1\% & 16196                       & 20\%                               & 16\%                               & 15\%                               & 15\%                               & \textbf{28\%}           & \textbf{25\%}           & \textbf{24\%}           & \textbf{24\%}             \\
Multi-dim Knapsack - Large (100)  & 23.8\% & 18171                       & 29\%                               & 28\%                               & 28\%                               & 28\%                               & \textbf{39\%}           & \textbf{38\%}           & \textbf{38\%}           & \textbf{38\%}             \\ \hline
Geo Mean                          &        & 4149                        & $\;$ 13.3\%                        & $\;$9.4\%                          & $\;$7.7\%                          & $\;$7.4\%                          & \textbf{15.1\%}         & \textbf{11.6\%}         & \textbf{10.1\%}         & \textbf{10.0\%}           \\ \hline
\end{tabular}
}}
\end{center}
\end{table}

On average, Figure~\ref{fig:geomean_eta_primalgap} demonstrates a steady decrease in the mean tree size as $\eta$ is increased from 0 to 0.15, and no significant change is seen 
with increasing the value of $\eta $ thereafter. We therefore select $\eta=0.15$ for further computations and refer to the resulting last assignment based rule as LA-SB~(3,7,0.15).
We present a detailed comparison of the mean tree sizes of LA-SB~(3,7,0.15) against Def-SB and Eff-SB~(3,7) in Table~\ref{tab:la_0.15}, and Table~\ref{tab:la_0.15_pct} presents the same results in terms of percentage reduction in tree sizes relative to Def-SB. 
\textcolor{black}{Let $\mathcal{L}$ be the set of all leaf nodes. Then, the column $|\mathcal{L}_I|/|\mathcal{L}|$ represents the average fraction of all leaves that were pruned due to integrality or infeasibility. }
The following are the main observations:
\begin{enumerate}
    \item The overall mean tree size is 15.1\% to 10.0\% smaller than Def-SB as the primal gap increases from 0\% to 10\%.   
    \item The impact of the last assignment is more pronounced for larger primal gaps and leads to 1.8\% to 2.6\% improvement over Eff-SB~(3,7) as the primal gap increases from 0\% to 10\%. 
    \item \textcolor{black}{The last assignment rule yields smaller trees than Eff-SB~(3,7) for some primal gap on problems where the average fraction of leaf nodes pruned due to integrality or infeasibility is at least 7\%. On other problems, any change in tree size is small except for set packing.}
\end{enumerate}

\begin{remark}
    \textcolor{black}{In this work, nodes are processed using the best-bound rule. If nodes are processed using the depth-first rule, the effectiveness of the last assignment rule may depend on the quality of the primal bound. For example, in a set-packing problem, diving on branches fixing variables to 0 will eventually result in an integral solution in the absence of any primal bound. This may potentially lead to an incorrect estimate of the direction of asymmetry using the last assignment rule. However, with a good primal bound, pruning by bound is more likely to occur before such an integral solution may be reached in this situation. } 

\end{remark}
\section{General Rebalancing}\label{sec:RLA}

The approach based on the last assignment rule in the previous section can be viewed as rebalancing of the branch-and-bound tree to account for asymmetrical pruning due to integrality and infeasibility. An asymmetrical trend may arise in pruning by bound if the LP gains on one branch are consistently smaller than those on the other branch due to structural properties of a problem. 
In this section, we extend the approach from Section~\ref{sec:la} to rebalance the 
%overall 
trees by 
%additionally 
accounting for any asymmetrical trend in pruning by bounds.

% \margin{why balanced trees are preferred?}
%An asymmetrical trend may arise in pruning by bound if the LP gains on one branch are consistently smaller than those on the other branch due to structural properties of a problem. 
If the global trend suggests that the gap on one branch ($0$ or $1$ side) is consistently harder to close than on the other, leading to larger subtrees, it may be beneficial to place greater emphasis on the closed gap on 
%the harder side. 
that branch.
To this end, we modify the Algorithm~\ref{algo:last_assignment} by replacing the set of integral and infeasible leaves, $\mathcal{L}_I$, with the set of all leaf nodes $\mathcal{L}$. We then use this modified algorithm to determine the exponents $a_0$ and $a_1$ in Eq.~\eqref{eq:asymmetrical_score}. We call this the \textit{rebalancing last assignment rule}. 

\subsection{Experiments with rebalancing last assignment rule}

\begin{figure}
    \centering
    \includegraphics[width=\linewidth]{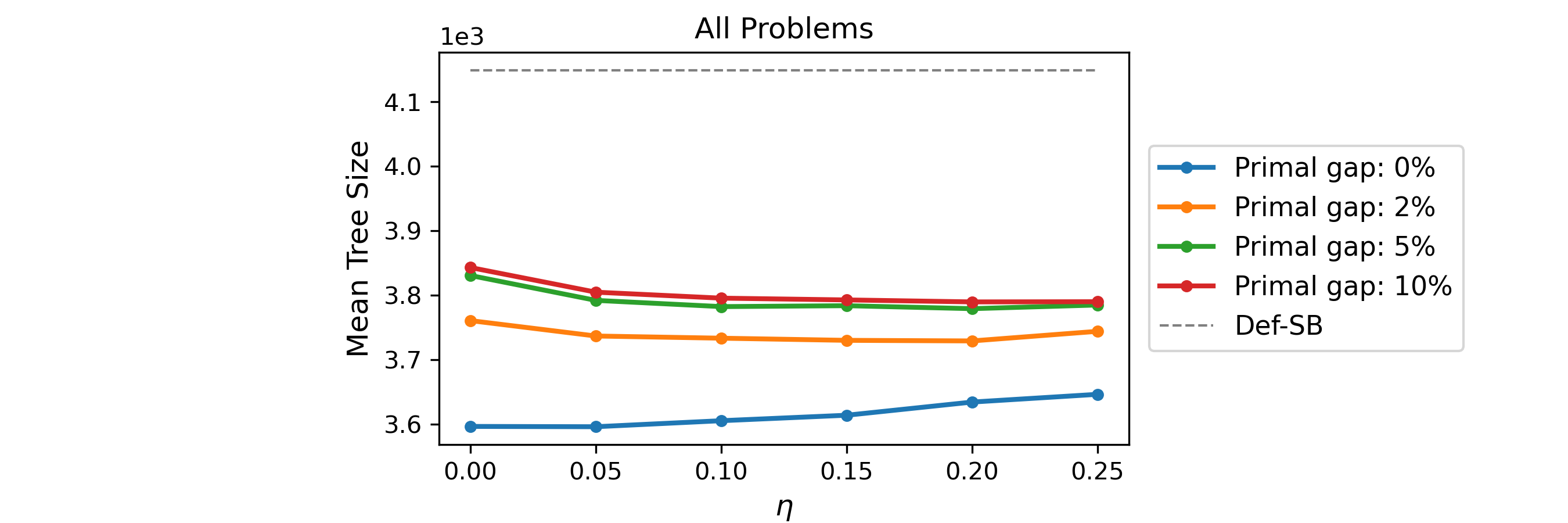}
    \caption{Geometric mean of tree sizes across problem classes when using the rebalancing last assignment rule and varying parameter $\eta$ and primal gap.}
    \label{fig:geomean_etabound_primalgap}
\end{figure}

\begin{figure}
    \centering
    \includegraphics[width=\linewidth]{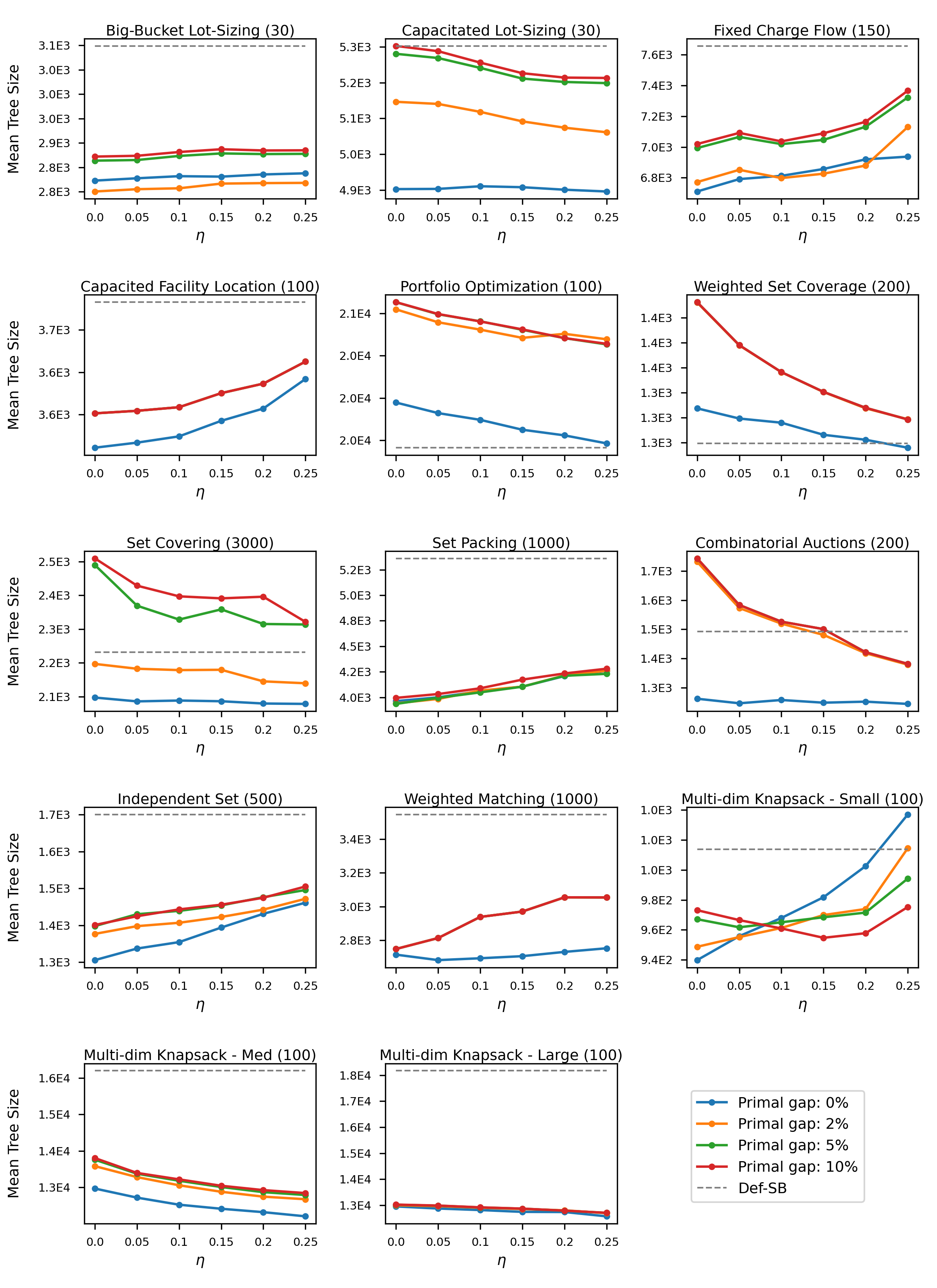}
    \caption{Impact of parameter $\eta$ on the geometric mean of tree sizes using the rebalancing last assignment rule for different primal gaps.}
    \label{fig:etabound_plots_quasiconvex}
\end{figure}

We conduct a similar search on $\eta \in \{ 0.05, 0.1, 0.15, 0.2, 0.25\}$ with the same values of other parameters. Figure~\ref{fig:etabound_plots_quasiconvex} presents the geometric mean of tree sizes for individual problems against different choices of $\eta$, and Figure~\ref{fig:geomean_etabound_primalgap} presents the geometric mean of tree sizes aggregated across all problems. Note that $\eta=0$ implies that the rebalancing last assignment rule was not used and corresponds to the Eff-SB~(3,7) score.

Figure~\ref{fig:etabound_plots_quasiconvex} indicates that here again the trends with respect to $\eta$ differ significantly across different problems. We make the following key observations.
\begin{enumerate}
    \item Problems with a trend of decreasing tree sizes with increasing $\eta$ include capacitated lot-sizing, portfolio optimization, weighted set-coverage, set covering, combinatorial auctions and multi-dimensional knapsack with medium right-hand side. These are precisely the problems where the slopes in Figure~\ref{fig:amin_plots_quasiconvex} are negative at the point 0.3, that is, increasing $a_{min}$ leads to smaller trees. On these problems, the ideal parameters in Table~\ref{tab:opteffSB} correspond to a value of $a^*_{min} \geq 0.3$ (and typically $>0.3$). {In other words, for these classes of problems, the choice of 
    $a^*_{min} = 0.3$, placed far too little emphasis on the smaller gains, possibly leading to larger trees. Therefore, rebalancing last assignment is able to compensate for the fixed value of $a_{min}$. 
    
%    This means that a higher emphasis on the smaller gain was ideally needed for these problems. Rebalancing might be reversing this by increasing the importance on the larger subtree with smaller lp gain
}

    \item Problems with a trend of increasing tree sizes with increasing $\eta$ include big-bucket lot-sizing, fixed charge flow, capacitated facility location, set packing, independent set and weighted matching and multi-dimensional knapsack with small right-hand-side. On these problems, the slopes in Figure~\ref{fig:amin_plots_quasiconvex} are positive at the point 0.3 and the ideal parameters in Table~\ref{tab:opteffSB} correspond to a smaller value of $a_{min} \leq 0.2$.
\end{enumerate}

Figure~\ref{fig:geomean_etabound_primalgap} indicates that all primal gaps except $0\%$, have a decreasing trend of average tree sizes (across all 14 problem classes) with increasing values of $\eta$. Whereas, with 0\% primal gap, the average tree sizes increase with increasing values of $\eta$. We select $\eta=0.15$ as further increasing $\eta$ leads to a significant increase in the average tree sizes for the 0\% primal gap case,  
%for optimal primal bounds 
without a significant decrease in tree sizes for 
%suboptimal bounds. 
other primal gap values.
We refer to the rebalancing last assignment rule with this parameter as R-LA-SB~(3,7,0.15). Table~\ref{tab:rla_0.15} and~\ref{tab:rla_0.15_pct} compare the rebalancing last assignment rule with the last assignment rule for integrality and infeasibility.

\begin{table}
\begin{center}	
		\caption{Comparing performance of R-LA-SB~(3,7,0.15) and LA-SB~(3,7,0.15) based on shifted geometric mean of tree sizes rounded to integers. Smaller tree sizes are highlighted in bold.}\label{tab:rla_0.15}
\resizebox{1.\textwidth}{!}{
\setlength{\tabcolsep}{0.5em} 
  {\renewcommand{\arraystretch}{1.1}
\begin{tabular}{|l|c|r|rrrr|rrrr|}
\hline
Rule                              &        & \multicolumn{1}{c|}{Def-SB} & \multicolumn{4}{c|}{LA-SB (3,7,0.15)}                                                                   & \multicolumn{4}{c|}{R-LA-SB (3,7,0.15)}                                                                 \\ \hline
Primal Gap                        &        & \multicolumn{1}{c|}{}       & \multicolumn{1}{c}{0\%} & \multicolumn{1}{c}{2\%} & \multicolumn{1}{c}{5\%} & \multicolumn{1}{c|}{10\%} & \multicolumn{1}{c}{0\%} & \multicolumn{1}{c}{2\%} & \multicolumn{1}{c}{5\%} & \multicolumn{1}{c|}{10\%} \\ \hline
\multicolumn{1}{|c|}{Problem}     &        &                             &                         &                         &                         &                           &                         &                         &                         &                           \\ \hline
Big-Bucket Lot-Sizing (30)        & 7.3\%  & 3099                        & \textbf{2753}           & \textbf{2730}           & \textbf{2812}           & \textbf{2820}             & 2831                    & 2817                    & 2879                    & 2887                      \\
Capacitated Lot-Sizing (30)       & 1.9\%  & 5302                        & 4913                    & 5124                    & 5245                    & 5264                      & \textbf{4908}           & \textbf{5092}           & \textbf{5211}           & \textbf{5226}             \\
Fixed Charge Flow (150)           & 8.0\%  & 7658                        & \textbf{6442}           & \textbf{6508}           & \textbf{6704}           & \textbf{6722}             & 6857                    & 6827                    & 7046                    & 7089                      \\
Capacited Facility Location (100) & 1.2\%  & 3733                        & \textbf{3581}           & \textbf{3622}           & \textbf{3622}           & \textbf{3622}             & 3593                    & 3625                    & 3625                    & 3625                      \\
Portfolio Optimization (100)      & 19.6\% & 19965                       & \textbf{18257}          & \textbf{18491}          & \textbf{18516}          & \textbf{18516}            & 20050                   & 20484                   & 20523                   & 20525                     \\
Weighted Set Coverage (200)       & 0.5\%  & 1274                        & 1309                    & 1415                    & 1415                    & 1415                      & \textbf{1283}           & \textbf{1326}           & \textbf{1326}           & \textbf{1326}             \\
Set Covering (3000)               & 0.6\%  & 2232                        & 2090                    & 2195                    & 2486                    & 2504                      & \textbf{2086}           & \textbf{2180}           & \textbf{2358}           & \textbf{2391}             \\
Set Packing (1000)                & 6.1\%  & 5362                        & 4175                    & 4131                    & 4148                    & 4176                      & \textbf{4104}           & \textbf{4104}           & \textbf{4104}           & \textbf{4173}             \\
Combinatorial Auctions (200)      & 0.3\%  & 1493                        & 1262                    & 1732                    & 1744                    & 1744                      & \textbf{1248}           & \textbf{1481}           & \textbf{1501}           & \textbf{1501}             \\
Independent Set (500)             & 3.5\%  & 1701                        & \textbf{1317}           & \textbf{1386}           & \textbf{1410}           & \textbf{1405}             & 1394                    & 1423                    & 1454                    & 1456                      \\
Weighted Matching (1000)          & 1.4\%  & 3545                        & 2725                    & \textbf{2747}           & \textbf{2747}           & \textbf{2747}             & \textbf{2706}           & 2971                    & 2971                    & 2971                      \\
Multi-dim Knapsack - Small (100)  & 67\%   & 1014                        & 1028                    & 997                     & 971                     & 957                       & \textbf{982}            & \textbf{970}            & \textbf{968}            & \textbf{955}              \\
Multi-dim Knapsack - Med (100)    & 40.1\% & 16196                       & \textbf{11617}          & \textbf{12090}          & \textbf{12253}          & \textbf{12322}            & 12417                   & 12882                   & 13013                   & 13047                     \\
Multi-dim Knapsack - Large (100)  & 23.8\% & 18171                       & \textbf{11098}          & \textbf{11257}          & \textbf{11277}          & \textbf{11289}            & 12744                   & 12853                   & 12867                   & 12869                     \\ \hline
Geo Mean                          &        & 4149                        & \textbf{3524}           & \textbf{3670}           & \textbf{3729}           & \textbf{3733}             & 3614                    & 3730                    & 3784                    & 3793                      \\ \hline
\end{tabular}
}}
% \end{center}
% \end{table}
\vspace{24pt}

% \begin{table}
% \begin{center}		
		\caption{
        Comparing performance of R-LA-SB~(3,7,0.15) and LA-SB~(3,7,0.15). Tree sizes are presented for Def-SB, and the percentage reduction in tree size relative to Def-SB is presented for others.}\label{tab:rla_0.15_pct}
\resizebox{1.\textwidth}{!}{
\setlength{\tabcolsep}{0.5em} 
  {\renewcommand{\arraystretch}{1.2}
% Please add the following required packages to your document preamble:
% \usepackage[table,xcdraw]{xcolor}
% Beamer presentation requires \usepackage{colortbl} instead of \usepackage[table,xcdraw]{xcolor}
\begin{tabular}{|l|c|r|rrrr|rrrr|}
\hline
Rule                              &        & \multicolumn{1}{c|}{Def-SB} & \multicolumn{4}{c|}{LA-SB (3,7,0.15)}                                                                   & \multicolumn{4}{c|}{R-LA-SB (3,7,0.15)}                                                                 \\ \hline
Primal Gap                        &        & \multicolumn{1}{c|}{}       & \multicolumn{1}{c}{0\%} & \multicolumn{1}{c}{2\%} & \multicolumn{1}{c}{5\%} & \multicolumn{1}{c|}{10\%} & \multicolumn{1}{c}{0\%} & \multicolumn{1}{c}{2\%} & \multicolumn{1}{c}{5\%} & \multicolumn{1}{c|}{10\%} \\ \hline
\multicolumn{1}{|c|}{Problem}     &        &                             &                         &                         &                         &                           &                         &                         &                         &                           \\ \hline
Big-Bucket Lot-Sizing (30)        & 7.3\%  & 3099                        & \textbf{11\%}           & \textbf{12\%}           & \textbf{9\%}            & \textbf{9\%}              & 9\%                     & 9\%                     & 7\%                     & 7\%                       \\
Capacitated Lot-Sizing (30)       & 1.9\%  & 5302                        & 7\%                     & 3\%                     & 1\%                     & 1\%                       & \textbf{7\%}            & \textbf{4\%}            & \textbf{2\%}            & \textbf{1\%}              \\
Fixed Charge Flow (150)           & 8.0\%  & 7658                        & \textbf{16\%}           & \textbf{15\%}           & \textbf{12\%}           & \textbf{12\%}             & 10\%                    & 11\%                    & 8\%                     & 7\%                       \\
Capacited Facility Location (100) & 1.2\%  & 3733                        & \textbf{4\%}            & \textbf{3\%}            & \textbf{3\%}            & \textbf{3\%}              & 4\%                     & 3\%                     & 3\%                     & 3\%                       \\
Portfolio Optimization (100)      & 19.6\% & 19965                       & \textbf{9\%}            & \textbf{7\%}            & \textbf{7\%}            & \textbf{7\%}              & 0\%                     & -3\%                    & -3\%                    & -3\%                      \\
Weighted Set Coverage (200)       & 0.5\%  & 1274                        & -3\%                    & -11\%                   & -11\%                   & -11\%                     & \textbf{-1\%}           & \textbf{-4\%}           & \textbf{-4\%}           & \textbf{-4\%}             \\
Set Covering (3000)               & 0.6\%  & 2232                        & 6\%                     & 2\%                     & -11\%                   & -12\%                     & \textbf{7\%}            & \textbf{2\%}            & \textbf{-6\%}           & \textbf{-7\%}             \\
Set Packing (1000)                & 6.1\%  & 5362                        & 22\%                    & 23\%                    & 23\%                    & 22\%                      & \textbf{23\%}           & \textbf{23\%}           & \textbf{23\%}           & \textbf{22\%}             \\
Combinatorial Auctions (200)      & 0.3\%  & 1493                        & 15\%                    & -16\%                   & -17\%                   & -17\%                     & \textbf{16\%}           & \textbf{1\%}            & \textbf{-1\%}           & \textbf{-1\%}             \\
Independent Set (500)             & 3.5\%  & 1701                        & \textbf{23\%}           & \textbf{18\%}           & \textbf{17\%}           & \textbf{17\%}             & 18\%                    & 16\%                    & 15\%                    & 14\%                      \\
Weighted Matching (1000)          & 1.4\%  & 3545                        & 23\%                    & \textbf{23\%}           & \textbf{23\%}           & \textbf{23\%}             & \textbf{24\%}           & 16\%                    & 16\%                    & 16\%                      \\
Multi-dim Knapsack - Small (100)  & 67\%   & 1014                        & -1\%                    & 2\%                     & 4\%                     & 6\%                       & \textbf{3\%}            & \textbf{4\%}            & \textbf{4\%}            & \textbf{6\%}              \\
Multi-dim Knapsack - Med (100)    & 40.1\% & 16196                       & \textbf{28\%}           & \textbf{25\%}           & \textbf{24\%}           & \textbf{24\%}             & 23\%                    & 20\%                    & 20\%                    & 19\%                      \\
Multi-dim Knapsack - Large (100)  & 23.8\% & 18171                       & \textbf{39\%}           & \textbf{38\%}           & \textbf{38\%}           & \textbf{38\%}             & 30\%                    & 29\%                    & 29\%                    & 29\%                      \\ \hline
Geo Mean                          &        & 4149                        & \textbf{15.1\%}         & \textbf{11.6\%}         & \textbf{10.1\%}         & \textbf{10.0\%}           & 12.9\%                  & 10.1\%                  & 8.8\%                   & 8.6\%                     \\ \hline
\end{tabular}
}}
\end{center}
\end{table}

The following are the key observations,
\begin{enumerate}
    \item The problems on which R-LA-SB~(3,7,0.15) yields significantly smaller trees than LA-SB~(3,7,0.15) are capacitated lot-sizing, weighted set coverage, set-covering, combinatorial auctions and multi-dimensional knapsack with small right-hand sides. Except for the last problem, all others have 
    % negative slopes in Figure~\ref{fig:etabound_plots_quasiconvex} and 
    less than 2\% of leaf nodes pruned due to integrality or infeasibility. 
    
    \item The problems on which R-LA-SB~(3,7,0.15) yields significantly larger trees than LA-SB~(3,7,0.15) are big-bucket lot-sizing, fixed-charge flow, portfolio optimization, independent set, weighted matching and multi-knapsack with medium and large right-hand sides. Except for independent set and weighted matching, all others have at least 7\% leaf nodes pruned due to integrality or infeasibility.
    % \item The problems on which the rebalancing modification yields smaller trees are capacitated lot-sizing, weighted set coverage, set covering and combinatorial auctions. These problems have less than 2\% of leaf nodes pruned due to integrality or infeasibility. 
    % \item Other problems with less than 2\% of integral or infeasible leaves include capacitated facility location, where the rebalancing modification has a similar performance as the original last assignment, and weighted matching, where the rebalancing modification leads to larger tree sizes.
    % \item On all other problems, the tree sizes with the rebalancing modification are similar or larger.
\end{enumerate}
The above results are a validation for the motivation for designing these rules, since LA-SB~(3,7,0.15) typically improves the tree sizes when the fraction of integer and infeasible leaf nodes is high, and R-LA-SB~(3,7,0.15) improves tree sizes when the fraction of integer and infeasible leaf nodes is low. 

Lastly, while LA-SB~(3,7,0.15) has smaller trees on average across problems, R-LA-SB~(3,7,0.15) is more robust in the sense that any increase in the resulting tree sizes as compared to Def-SB is no larger than 7\%.

\subsection{Pruning aware last assignment rule}
Motivated by the results in Tables~\ref{tab:rla_0.15} and~\ref{tab:rla_0.15_pct},
we now apply a simple policy to switch between LA-SB~(3,7,0.15) and R-LA-SB~(3,7,0.15) in order to obtain the benefits of both. Since in our experiments less than $2\%$ integer and infeasible nodes indicates the use of LA-SB~(3,7,0.15) and more than $7\%$ integer and infeasible nodes indicates the use of R-LA-SB~(3,7,0.15), we choose to use a threshold policy of $5\%$, i.e., if the integer and infeasible nodes are less than $5\%$, we use R-LA-SB~(3,7,0.15), otherwise we use LA-SB~(3,7,0.15). We call this rule pruning-aware last assignment and refer to it as PA-LA~(3,7,0.15). 

Tables~\ref{tab:pala_0.15} and~\ref{tab:pala_0.15_pct} present the results of the experiments with the PA-LA~(3,7,0.15) rule and compare its performance with Def-SB and Eff-SB~(3,7). 
These results indicate that PA-LA-SB~(3,7,0.15) indeed strikes the balance between small trees and robustness, thereby achieving the benefits of both the constituting rules. 
On average across problem classes, it leads to smaller trees than all previously considered rules for all initial primal gaps except $0\%$ where LA-SB~(3,7,0.15) tree sizes are smaller. Its performance for $0\%$ primal gap still remains superior compared to higher primal gaps. It is also robust as it leads to smaller trees than Def-SB on more combinations of problems and primal gaps than previously considered rules. Any increase in tree sizes compared to Def-SB is again at most 7\%.

\begin{table}
\begin{center}
		
		\caption{Evaluating PA-LA-SB~(3,7,0.15) based on shifted geometric mean of tree sizes rounded to integers. Smaller tree sizes between Eff-SB(3,7) and PA-LA-SB(3,7,0.15) are highlighted in bold.}\label{tab:pala_0.15}
\resizebox{1.\textwidth}{!}{
\setlength{\tabcolsep}{0.5em} 
  {\renewcommand{\arraystretch}{1.1}
\begin{tabular}{|l|c|r|rrrr|rrrr|}
\hline
Rule                              &                                       & \multicolumn{1}{c|}{Def-SB} & \multicolumn{4}{c|}{Eff-SB (3,7)}                                                                       & \multicolumn{4}{c|}{PA-LA-SB (3,7,0.15)}                                                                   \\ \hline
Primal Gap                        &                                       & \multicolumn{1}{c|}{}       & \multicolumn{1}{c}{0\%} & \multicolumn{1}{c}{2\%} & \multicolumn{1}{c}{5\%} & \multicolumn{1}{c|}{10\%} & \multicolumn{1}{c}{0\%} & \multicolumn{1}{c}{2\%} & \multicolumn{1}{c}{5\%} & \multicolumn{1}{c|}{10\%} \\ \hline
\multicolumn{1}{|c|}{Problem}     & $|\mathcal{L}_I|/|\mathcal{L}|$ &                             &                         &                         &                         &                           &                         &                         &                         &                           \\ \hline
Big-Bucket Lot-Sizing (30)        & 7.3\%                           & 3099                        & 2823                    & 2801                    & 2864                    & 2872                      & \textbf{2780}           & \textbf{2751}           & \textbf{2827}           & \textbf{2835}             \\
Capacitated Lot-Sizing (30)       & 1.9\%                           & 5302                        & 4902                    & 5146                    & 5280                    & 5302                      & \textbf{4900}           & \textbf{5085}           & \textbf{5206}           & \textbf{5221}             \\
Fixed Charge Flow (150)           & 8.0\%                           & 7658                        & 6711                    & 6772                    & 6993                    & 7019                      & \textbf{6706}           & \textbf{6745}           & \textbf{6913}           & \textbf{6931}             \\
Capacited Facility Location (100) & 1.2\%                           & 3733                        & \textbf{3561}           & \textbf{3602}           & \textbf{3602}           & \textbf{3602}             & 3593                    & 3625                    & 3625                    & 3625                      \\
Portfolio Optimization (100)      & 19.6\%                          & 19965                       & 20179                   & 20619                   & 20653                   & 20654                     & \textbf{18435}          & \textbf{18829}          & \textbf{18876}          & \textbf{18876}            \\
Weighted Set Coverage (200)       & 0.5\%                           & 1274                        & 1309                    & 1415                    & 1415                    & 1415                      & \textbf{1283}           & \textbf{1326}           & \textbf{1326}           & \textbf{1326}             \\
Set Covering (3000)               & 0.6\%                           & 2232                        & 2097                    & 2197                    & 2490                    & 2509                      & \textbf{2086}           & \textbf{2180}           & \textbf{2358}           & \textbf{2391}             \\
Set Packing (1000)                & 6.1\%                           & 5362                        & \textbf{3963}           & \textbf{3941}           & \textbf{3935}           & \textbf{3995}             & 4159                    & 4118                    & 4129                    & 4197                      \\
Combinatorial Auctions (200)      & 0.3\%                           & 1493                        & 1262                    & 1732                    & 1744                    & 1744                      & \textbf{1248}           & \textbf{1481}           & \textbf{1501}           & \textbf{1501}             \\
Independent Set (500)             & 3.5\%                           & 1701                        & \textbf{1306}           & \textbf{1377}           & \textbf{1398}           & \textbf{1401}             & 1394                    & 1424                    & 1455                    & 1457                      \\
Weighted Matching (1000)          & 1.4\%                           & 3545                        & 2715                    & \textbf{2749}                    & \textbf{2749}                    & \textbf{2749}                      & \textbf{2705}           & 2967           & 2967           & 2967             \\
Multi-dim Knapsack - Small (100)  & 67\%                            & 1014                        & \textbf{940}            & \textbf{949}            & \textbf{967}            & 973                       & 1028                    & 997                     & 971                     & \textbf{957}              \\
Multi-dim Knapsack - Med (100)    & 40.1\%                          & 16196                       & 12971                   & 13584                   & 13760                   & 13805                     & \textbf{11617}          & \textbf{12090}          & \textbf{12253}          & \textbf{12322}            \\
Multi-dim Knapsack - Large (100)  & 23.8\%                          & 18171                       & 12953                   & 12996                   & 13019                   & 13023                     & \textbf{11098}          & \textbf{11257}          & \textbf{11277}          & \textbf{11289}            \\ \hline
Geo Mean                          &                                 & 4149                        & 3596                    & 3761                    & 3831                    & 3843                      & \textbf{3542}           & \textbf{3652}           & \textbf{3700}           & \textbf{3709}             \\ \hline
\end{tabular}
}}
% \end{center}
% \end{table}
\vspace{24pt}

% \begin{table}
% \begin{center}		
		\caption{
        Comparison of the shifted geometric mean of tree sizes to evaluate PA-LA-SB~(3,7,0.15). Tree sizes are presented for Def-SB, and the percentage reduction in tree size relative to Def-SB is presented for others.}\label{tab:pala_0.15_pct}
\resizebox{1.\textwidth}{!}{
\setlength{\tabcolsep}{0.5em} 
  {\renewcommand{\arraystretch}{1.2}
% Please add the following required packages to your document preamble:
% \usepackage[table,xcdraw]{xcolor}
% Beamer presentation requires \usepackage{colortbl} instead of \usepackage[table,xcdraw]{xcolor}
\begin{tabular}{|l|c|r|rrrr|rrrr|}
\hline
Rule                              &        & \multicolumn{1}{c|}{Def-SB} & \multicolumn{4}{c|}{Eff-SB (3,7)}                                                                                                                 & \multicolumn{4}{c|}{PA-LA-SB (3,7,0.15)}                                                                   \\ \hline
Primal Gap                        &        & \multicolumn{1}{c|}{}       & \multicolumn{1}{c}{0\%}            & \multicolumn{1}{c}{2\%}            & \multicolumn{1}{c}{5\%}            & \multicolumn{1}{c|}{10\%}          & \multicolumn{1}{c}{0\%} & \multicolumn{1}{c}{2\%} & \multicolumn{1}{c}{5\%} & \multicolumn{1}{c|}{10\%} \\ \hline
\multicolumn{1}{|c|}{Problem}     &    $|\mathcal{L}_I|/|\mathcal{L}|$   &                             &                                    &                                    &                                    &                                    &                         &                         &                         &                           \\ \hline
Big-Bucket Lot-Sizing (30)        & 7.3\%                           & 3099                        & 9\%                     & 10\%                    & 8\%                     & 7\%                       & \textbf{10\%}           & \textbf{11\%}           & \textbf{9\%}            & \textbf{9\%}              \\
Capacitated Lot-Sizing (30)       & 1.9\%                           & 5302                        & 8\%                     & 3\%                     & 0\%                     & 0\%                       & \textbf{8\%}            & \textbf{4\%}            & \textbf{2\%}            & \textbf{2\%}              \\
Fixed Charge Flow (150)           & 8.0\%                           & 7658                        & 12\%                    & 12\%                    & 9\%                     & 8\%                       & \textbf{12\%}           & \textbf{12\%}           & \textbf{10\%}           & \textbf{9\%}              \\
Capacited Facility Location (100) & 1.2\%                           & 3733                        & \textbf{5\%}            & \textbf{4\%}            & \textbf{4\%}            & \textbf{4\%}              & 4\%                     & 3\%                     & 3\%                     & 3\%                       \\
Portfolio Optimization (100)      & 19.6\%                          & 19965                       & -1\%                    & -3\%                    & -3\%                    & -3\%                      & \textbf{8\%}            & \textbf{6\%}            & \textbf{5\%}            & \textbf{5\%}              \\
Weighted Set Coverage (200)       & 0.5\%                           & 1274                        & -3\%                    & -11\%                   & -11\%                   & -11\%                     & \textbf{-1\%}           & \textbf{-4\%}           & \textbf{-4\%}           & \textbf{-4\%}             \\
Set Covering (3000)               & 0.6\%                           & 2232                        & 6\%                     & 2\%                     & -12\%                   & -12\%                     & \textbf{7\%}            & \textbf{2\%}            & \textbf{-6\%}           & \textbf{-7\%}             \\
Set Packing (1000)                & 6.1\%                           & 5362                        & \textbf{26\%}           & \textbf{27\%}           & \textbf{27\%}           & \textbf{26\%}             & 22\%                    & 23\%                    & 23\%                    & 22\%                      \\
Combinatorial Auctions (200)      & 0.3\%                           & 1493                        & 15\%                    & -16\%                   & -17\%                   & -17\%                     & \textbf{16\%}           & \textbf{1\%}            & \textbf{-1\%}           & \textbf{-1\%}             \\
Independent Set (500)             & 3.5\%                           & 1701                        & \textbf{23\%}           & \textbf{19\%}           & \textbf{18\%}           & \textbf{18\%}             & 18\%                    & 16\%                    & 14\%                    & 14\%                      \\
Weighted Matching (1000)          & 1.4\%                           & 3545                        & 23\%                    & 22\%                    & 22\%                    & 22\%                      & \textbf{24\%}           & \textbf{16\%}           & \textbf{16\%}           & \textbf{16\%}             \\
Multi-dim Knapsack - Small (100)  & 67\%                            & 1014                        & \textbf{7\%}            & \textbf{6\%}            & \textbf{5\%}            & 4\%                       & -1\%                    & 2\%                     & 4\%                     & \textbf{6\%}              \\
Multi-dim Knapsack - Med (100)    & 40.1\%                          & 16196                       & 20\%                    & 16\%                    & 15\%                    & 15\%                      & \textbf{28\%}           & \textbf{25\%}           & \textbf{24\%}           & \textbf{24\%}             \\
Multi-dim Knapsack - Large (100)  & 23.8\%                          & 18171                       & 29\%                    & 28\%                    & 28\%                    & 28\%                      & \textbf{39\%}           & \textbf{38\%}           & \textbf{38\%}           & \textbf{38\%}             \\ \hline
Geo Mean                          &                                 & 4149                        & 13.3\%                  & 9.4\%                   & 7.7\%                   & 7.4\%                     & \textbf{14.6\%}         & \textbf{12.0\%}         & \textbf{10.8\%}         & \textbf{10.6\%}           \\ \hline
\end{tabular}
}}
\end{center}
\end{table}

\section{Computational Experiments on MIPLIB instances}\label{sec:MIPLIB}

In previous sections, we have considered randomly generated instances of problems with known structure to enable us to draw insights and design better rules. This section uses instances from the Benchmark Set of MIPLIB 2017~\cite{miplib2017} to test the \textcolor{black}{efficacious strong branching (Eff-SB~(3,7)) and pruning aware last assignment rule (PA-LA-SB~(3,7,0.15))} on more general problems.

We exclude instances marked hard or infeasible from our test set. We use Gurobi~\cite{gurobi} to presolve the remaining instances and only include instances with at most 2000 binary variables and no general integer variables after presolve. Finally, we use SCIP~\cite{BestuzhevaEtal2021OO} to add \textcolor{black}{up to 5 rounds of } cuts at the root node 
% {with a time limit of 30 minutes,} 
and save the resulting models. These models are then used in the computational experiments with the same set-up as described in Section~\ref{sec:intro_expt_setup}. The final test set contains 51 instances after excluding 2 that have an integer LP solution at the root node.

We initialize the branch-and-bound process with primal solutions at a gap of 0\%, 2\%, 5\% and 10\% from optimal. Every combination of instance, rule and primal gap is solved with 3 different random seeds to account for possible dual degeneracy in the problems, even though our textbook implementations of branch-and-bound algorithm and full strong branching do not use any random processes. A node limit of 20,000 is imposed, and the remaining gap at the time of termination relative to the integrality gap at the root node is presented for incomplete runs. The geometric mean of the remaining gap is calculated with a 1\% shift.

\subsection{Results}
% \margin{Not sure if separating LA makes sense with hybrid.}

\begin{table}
\begin{center}
		\caption{Evaluation of Eff-SB (3,7) and PA-LA-SB~(3,7,0.15) on MIPLIB instances; presenting shifted geometric mean of tree sizes (numbers) if solved within node limit, else gap remaining (percentages) at termination.
        }\label{tab:miplib_sb}
\resizebox{1.\linewidth}{!}{
\setlength{\tabcolsep}{1.2em} 
  {\renewcommand{\arraystretch}{1.}
\begin{tabular}{|l|r|rrrr|rrrr|}
\hline
Rule                                                                       & \multicolumn{1}{c|}{Def-SB} & \multicolumn{4}{c|}{Eff-SB (3,7)}                                                                       & \multicolumn{4}{c|}{PA-LA-SB (3,7,0.15)}                                                                \\ \hline
Primal Gap                                                                 & \multicolumn{1}{c|}{0\%}    & \multicolumn{1}{c}{0\%} & \multicolumn{1}{c}{2\%} & \multicolumn{1}{c}{5\%} & \multicolumn{1}{c|}{10\%} & \multicolumn{1}{c}{0\%} & \multicolumn{1}{c}{2\%} & \multicolumn{1}{c}{5\%} & \multicolumn{1}{c|}{10\%} \\ \hline
\multicolumn{1}{|c|}{Instance}                                             &                             &                         &                         &                         &                           &                         &                         &                         &                           \\ \hline
app1-1                                                                     & 1679                        & \textbf{42}             & \textbf{42}             & \textbf{42}             & \textbf{42}               & \textbf{42}             & \textbf{42}             & \textbf{42}             & \textbf{42}               \\
assign1-5-8                                                                & 61.0\%                      & \textbf{60.9\%}         & \textbf{60.9\%}         & \textbf{60.9\%}         & \textbf{60.9\%}           & \textbf{60.9\%}         & \textbf{60.9\%}         & \textbf{60.9\%}         & \textbf{60.9\%}           \\
b1c1s1                                                                     & 43.4\%                      & \textbf{42.6\%}         & \textbf{42.6\%}         & \textbf{42.6\%}         & \textbf{42.6\%}           & \textbf{42.6\%}         & \textbf{42.6\%}         & \textbf{42.6\%}         & \textbf{42.6\%}           \\
beasleyC3                                                                  & 78.7\%                      & \textbf{77.4\%}         & \textbf{77.4\%}         & \textbf{77.4\%}         & \textbf{77.4\%}           & \textbf{77.4\%}         & \textbf{77.4\%}         & \textbf{77.4\%}         & \textbf{77.4\%}           \\
binkar10\_1                                                                & 25.1\%                      & 24.9\%                  & 25.0\%                  & 25.0\%                  & 25.0\%                    & \textbf{24.7\%}         & \textbf{24.8\%}         & \textbf{24.8\%}         & \textbf{24.8\%}           \\
cmflsp50-24-8-8                                                            & 2.9\%                       & 2.6\%                   & 5.1\%                   & 5.2\%                   & 5.1\%                     & \textbf{2.4\%}          & \textbf{4.7\%}          & \textbf{4.7\%}          & \textbf{4.7\%}            \\
cost266-UUE                                                                & 28.8\%                      & \textbf{27.7\%}         & \textbf{27.7\%}         & \textbf{27.7\%}         & \textbf{27.8\%}           & 28.2\%                  & 28.4\%                  & 28.3\%                  & 28.4\%                    \\
dano3\_3                                                                   & \textbf{28}                 & \textbf{18}             & 30                      & 30                      & 30                        & \textbf{18}             & 30                      & 30                      & 30                        \\
dano3\_5                                                                   & \textbf{144}                & \textbf{134}            & 244                     & 244                     & 244                       & \textbf{134}            & 232                     & 232                     & 232                       \\
exp-1-500-5-5                                                              & \textbf{57.4\%}             & 55.0\%                  & 55.0\%                  & 55.0\%                  & 55.0\%                    & 55.0\%                  & 55.0\%                  & 55.0\%                  & 55.0\%                    \\
glass-sc                                                                   & \textbf{26.3\%}             & 26.7\%                  & 26.7\%                  & 26.7\%                  & 26.7\%                    & 26.7\%                  & 26.7\%                  & 26.7\%                  & 26.7\%                    \\
glass4                                                                     & 100.0\%                     & \textbf{50.0\%}         & \textbf{50.0\%}         & \textbf{50.0\%}         & \textbf{50.0\%}           & \textbf{50.0\%}         & \textbf{50.0\%}         & \textbf{50.0\%}         & \textbf{50.0\%}           \\
gmu-35-40                                                                  & 100.0\%                     & 46.1\%                  & 46.9\%                  & 47.4\%                  & 73.4\%                    & \textbf{45.9\%}         & \textbf{45.9\%}         & \textbf{47.2\%}         & \textbf{71.5\%}           \\
gmu-35-50                                                                  & 100.0\%                     & 100.0\%                 & 99.5\%                  & 100.0\%                 & 100.0\%                   & 100.0\%                 & \textbf{99.4\%}         & 100.0\%                 & 100.0\%                   \\
lotsize                                                                    & 92.0\%                      & \textbf{91.7\%}         & \textbf{91.7\%}         & \textbf{91.7\%}         & \textbf{91.7\%}           & 91.8\%                  & 91.8\%                  & 91.8\%                  & 91.8\%                    \\
map10                                                                      & \textbf{732}                & 570                     & 608                     & 780                     & 982                       & \textbf{560}            & \textbf{592}            & \textbf{714}            & 854                       \\
mc11                                                                       & 79.8\%                      & 79.8\%                  & 79.8\%                  & 79.8\%                  & 79.8\%                    & 79.8\%                  & 79.8\%                  & 79.8\%                  & 79.8\%                    \\
mcsched                                                                    & 7.8\%                       & \textbf{1.4\%}          & 4.2\%                   & 4.2\%                   & 4.2\%                     & \textbf{1.4\%}          & \textbf{2.5\%}          & \textbf{2.6\%}          & \textbf{2.6\%}            \\
milo-v12-6-r2-40-1                                                         & 40.0\%                      & \textbf{37.8\%}         & \textbf{37.9\%}         & \textbf{37.9\%}         & \textbf{37.9\%}           & 38.5\%                  & 38.6\%                  & 38.6\%                  & 38.6\%                    \\
neos-3627168-kasai                                                         & 19.2\%                      & \textbf{19.1\%}         & \textbf{19.1\%}         & \textbf{19.1\%}         & \textbf{19.1\%}           & \textbf{19.1\%}         & \textbf{19.1\%}         & \textbf{19.1\%}         & \textbf{19.1\%}           \\
neos-5188808-nattai                                                        & 17.3\%                      & \textbf{3.4\%}          & \textbf{3.8\%}          & \textbf{3.8\%}          & \textbf{4.0\%}            & 6.9\%                   & 7.3\%                   & 5.8\%                   & 7.8\%                     \\
neos-848589                                                                & 95.6\%                      & 95.6\%                  & 95.6\%                  & 95.6\%                  & 95.6\%                    & 95.6\%                  & 95.6\%                  & 95.6\%                  & 95.6\%                    \\
neos-860300                                                                & 270                         & \textbf{188}            & \textbf{190}            & \textbf{210}            & 228                       & \textbf{188}            & \textbf{190}            & \textbf{210}            & \textbf{226}              \\
neos-911970                                                                & 8.3\%                       & 8.3\%                   & 8.3\%                   & 8.3\%                   & 8.3\%                     & 8.3\%                   & 8.3\%                   & 8.3\%                   & 8.3\%                     \\
neos17                                                                     & 7.8\%                       & 6.8\%                   & 6.8\%                   & 6.8\%                   & 6.8\%                     & \textbf{4.7\%}          & \textbf{4.7\%}          & \textbf{4.7\%}          & \textbf{7.7\%}            \\
neos5                                                                      & \textbf{48.0\%}             & 48.4\%                  & 48.4\%                  & 48.4\%                  & 48.4\%                    & 48.4\%                  & 48.4\%                  & 48.4\%                  & 48.4\%                    \\
ns1830653                                                                  & 13402                       & 10069                   & 10053                   & 10408                   & 10619                     & \textbf{9753}           & \textbf{9791}           & \textbf{9946}           & \textbf{10497}            \\
p200x1188c                                                                 & 55.9\%                      & 50.6\%                  & 50.6\%                  & 50.6\%                  & 50.6\%                    & \textbf{50.3\%}         & \textbf{50.3\%}         & \textbf{50.3\%}         & \textbf{50.3\%}           \\
pg                                                                         & 61.3\%                      & \textbf{60.8\%}         & \textbf{60.8\%}         & \textbf{60.8\%}         & \textbf{60.8\%}           & \textbf{60.8\%}         & \textbf{60.8\%}         & \textbf{60.8\%}         & \textbf{60.8\%}           \\
pg5\_34                                                                    & 49.4\%                      & \textbf{48.8\%}         & \textbf{48.8\%}         & \textbf{48.8\%}         & \textbf{48.8\%}           & \textbf{48.8\%}         & \textbf{48.8\%}         & \textbf{48.8\%}         & \textbf{48.8\%}           \\
pk1                                                                        & \textbf{44.5\%}             & 46.1\%                  & 46.1\%                  & 46.2\%                  & 46.2\%                    & 46.1\%                  & 46.1\%                  & 46.1\%                  & 46.2\%                    \\
ran14x18-disj-8                                                            & 51.6\%                      & 51.5\%                  & 51.5\%                  & 51.5\%                  & 51.5\%                    & \textbf{50.9\%}         & \textbf{50.9\%}         & \textbf{50.9\%}         & \textbf{50.9\%}           \\
rd-rplusc-21                                                               & 69.3\%                      & \textbf{27.4\%}         & 59.2\%                  & 46.4\%                  & 67.6\%                    & 32.4\%                  & \textbf{28.9\%}         & \textbf{34.2\%}         & \textbf{32.9\%}           \\
rmatr100-p10                                                               & \textbf{588}                & 686                     & 754                     & 860                     & 860                       & 686                     & 716                     & 756                     & 756                       \\
seymour                                                                    & 42.7\%                      & \textbf{42.5\%}         & \textbf{42.5\%}         & \textbf{42.5\%}         & \textbf{42.5\%}           & \textbf{42.5\%}         & \textbf{42.5\%}         & \textbf{42.5\%}         & \textbf{42.5\%}           \\
seymour1                                                                   & \textbf{1395}               & 1303                    & 1407                    & 1407                    & 1407                      & \textbf{1285}           & 1449                    & 1449                    & 1449                      \\
sp150x300d                                                                 & \textbf{3.6\%}              & 9.2\%                   & 9.2\%                   & 9.2\%                   & 9.2\%                     & 9.2\%                   & 9.1\%                   & 9.1\%                   & 9.1\%                     \\
supportcase26                                                              & 59.0\%                      & \textbf{27.6\%}         & 27.6\%                  & 27.6\%                  & 27.6\%                    & \textbf{27.6\%}         & \textbf{27.4\%}         & \textbf{27.4\%}         & \textbf{27.4\%}           \\
supportcase40                                                              & \textbf{17442}              & \textbf{14400}          & \textbf{16383}          & 17469                   & 17469                     & 15546                   & 16981                   & 17521                   & 17521                     \\
supportcase7                                                               & 159                         & \textbf{88}             & \textbf{116}            & \textbf{125}            & \textbf{128}              & \textbf{88}             & \textbf{116}            & \textbf{125}            & \textbf{128}              \\
tr12-30                                                                    & \textbf{56.1\%}             & 56.2\%                  & 56.2\%                  & 56.2\%                  & 56.2\%                    & 56.2\%                  & 56.2\%                  & 56.2\%                  & 56.2\%                    \\
uct-subprob                                                                & \textbf{47.5\%}             & 49.0\%                  & 49.0\%                  & 49.0\%                  & 49.0\%                    & 48.9\%                  & 48.9\%                  & 48.9\%                  & 48.9\%                    \\ \hline
\multicolumn{1}{|c|}{Geo mean -- Nodes\footnotemark[1]} & 907                         & \textbf{586}            & 656                     & 697                     & 721                       & 587                     & \textbf{651}            & \textbf{678}            & \textbf{700}              \\ \hline
\multicolumn{1}{|c|}{Geo mean -- Gap\footnotemark[2]}   & 37.2\%                      & \textbf{31.2\%}         & 33.5\%                  & 33.3\%                  & 34.2\%                    & 31.6\%                  & \textbf{32.5\%}         & \textbf{32.5\%}         & \textbf{33.6\%}           \\ \hline
\end{tabular}
}}
\end{center}
\end{table}

\footnotetext[1]{Across instances where some combination of rule, primal bound and seed could be solved within the node limit.} 

\footnotetext[2]{Across instances where some combination of rule, primal bound and seed was unsolved within the node limit.}

The results of these experiments are shown in Table~\ref{tab:miplib_sb}. Five instances showed no change in dual bound after 20,000 nodes with any of the rules and are excluded from the results. In addition, four instances could not be solved in 14 days and are not included in the results. Each value in the table represents the shifted geometric mean of the number of nodes across different seeds, if the instance is solved within the node limit for all seeds. If the instance is unsolved for all seeds, the table reports the shifted geometric mean of the remaining gap as a percentage. Otherwise, we report geometric means of both, number of nodes and gap remaining.

\begin{enumerate}
{
    \item Similar to the trend seen in random instances, Eff-SB~(3, 7) and PA-LA-SB~(3,7,0.15) improve the average performance of the branch-and-bound process. The performance improvement is greater when provided with better primal bounds. 
    \item Over the 10 instances that could be solved within the node limit for some combination of rule, primal gap and seed, Eff-SB~(3,7) leads to a 35.4\%, 27.8\%, 23.1\% and 20.5\%  reduction in mean tree sizes relative to Def-SB with primal gaps 0\%, 2\%, 5\% and 10\%, respectively. On these instances, PA-LA-SB~(3,7,0.15) leads to 35.3\%, 28.3\%, 25.3\% and 22.8\% smaller tree sizes than Def-SB on the four primal gaps, respectively.
    \item Over the 32 instances that remained unsolved within the node limit for some combination of rule, primal gap and seed, the mean gap remaining at termination was 6.0\%, 3.7\%, 4.0\% and 3.1\% less with Eff-SB~(3, 7) than Def-SB with primal gaps 0\%, 2\%, 5\% and 10\%, respectively. The improvement in gap with PA-LA-SB~(3,7,0.15) compared to Def-SB was 5.6\%, 4.7\%, 4.8\% and 3.6\% respectively for the four primal gaps. 
    % \item {Eff~SB(3,7) improved the branch-and-bound algorithm's performance for some primal gaps on 32 instances and for all primal gaps on 27 instances. On 6 instances, its performance worsened for all primal gaps; on 2 instances, it led to no significant change.}
    \item PA-LA-SB~(3,7,0.15) leads to better performance than Eff-SB~(3,7) in terms of number of nodes and gap remaining when initialization is with 2\%, 5\% and 10\% primal gap.
    With 0\% primal gap, the performance of PA-LA-SB~(3,7,0.15) deteriorates only marginally as compared to Eff-SB~(3,7) in terms of nodes as well as gap remaining. Its performance, however, remains better than Def-SB and also better than the performance of PA-LA-SB~(3,7,0.15) on other primal gaps. The above trends indicate improved robustness with respect to primal gaps and is similar to the observations on the random test bed. 
    }
\end{enumerate}

\subsection{Extension to Reliability Branching}
In practice, FSB is computationally very expensive and impractical for larger instances. Solvers therefore use reliability-branching (RB), which attempts to approximate strong branching by predicting LP gains instead of solving LPs to compute their exact values. 
% Instead of solving LPs to compute the gains, reliability-branching estimates the gains from previous branching outcomes. 

We are now interested in assessing if the reduction in sizes of strong branching trees also translates to reliability-branching as we extend the idea of using efficacious gains to it. We derive efficacious gain estimates by applying the same threshold of incumbent additive gap $\Delta_{p-d}$ to the estimates of LP gains derived using pseudocosts. That is, let $(\hat{\Delta}^0_i, \hat{\Delta}^1_i)$ be the estimates of LP gains derived using pseudo-costs or the strong branching evaluations. Then the estimated efficacious gains are computed as follows, 
\begin{equation*}
    \hat{q}^0_i =  \min \big\{ \max\{\hat{\Delta}^0_i, \epsilon\}, \ \Delta_{p-d} \big\}, \ \hat{q}^1_i =\min \big\{ \max\{\hat{\Delta}^1_i, \epsilon\}, \ \Delta_{p-d} \big\}, 
\end{equation*}
These are then used to compute variable scores in the same way as in EFF-SB. We call this efficacious reliability branching.

Reliability-branching is implemented as described in \cite{achterberg2007constraint}. At any node, we allow a maximum of 100 FSB evaluations, limit each LP evaluation to 500 dual simplex iterations and deem pseudocosts reliable if a variable has at least 8 data points on both branches. We do not apply early termination of strong branching evaluations in case the score does not improve, that is, the look ahead parameter is set to infinity.

% We run similar experiments as previously described for strong branching to evaluate the impact of incorporating primal information through efficacious gains estimates instead of the actual LP gain estimates. 

\begin{table}
\begin{center}
		\caption{Evaluating Eff-RB (3,7) and PA-LA-RB(3,7,0.15), on MIPLIB instances; presenting shifted geometric mean of tree sizes (numbers) if solved within node limit, else gap remaining (percentages) at termination.
        }\label{tab:miplib_rb}
\resizebox{1.\linewidth}{!}{
\setlength{\tabcolsep}{.2em} 
  {\renewcommand{\arraystretch}{1.1}
\begin{tabular}{|l|r|rrrr|rrrr|}
\hline
Rule                                      & \multicolumn{1}{c|}{Def-RB} & \multicolumn{4}{c|}{Eff-RB (3,7)}                                                                       & \multicolumn{4}{c|}{PA-LA-RB (3,7,0.15)}                                                                \\ \hline
Primal Gap                                & \multicolumn{1}{c|}{0\%}    & \multicolumn{1}{c}{0\%} & \multicolumn{1}{c}{2\%} & \multicolumn{1}{c}{5\%} & \multicolumn{1}{c|}{10\%} & \multicolumn{1}{c}{0\%} & \multicolumn{1}{c}{2\%} & \multicolumn{1}{c}{5\%} & \multicolumn{1}{c|}{10\%} \\ \hline
\multicolumn{1}{|c|}{Instance}            & \multicolumn{1}{c|}{}       & \multicolumn{1}{c}{}    & \multicolumn{1}{c}{}    & \multicolumn{1}{c}{}    & \multicolumn{1}{c|}{}     & \multicolumn{1}{c}{}    & \multicolumn{1}{c}{}    & \multicolumn{1}{c}{}    & \multicolumn{1}{c|}{}     \\ \hline
app1-1                                                                     & 42                                     & 42                                              & 42                                     & 42                      & 42                                     & 42                                     & 42                      & 42                      & 42                        \\
assign1-5-8                                                                & \textbf{64.0\%}                        & 69.2\%                                          & 69.2\%                                 & 69.2\%                  & 69.2\%                                 & 69.2\%                                 & 69.2\%                  & 69.2\%                  & 69.2\%                    \\
b1c1s1                                                                     & 54.6\%                                 & \textbf{53.1\%}                                 & 53.7\%                                 & 53.7\%                  & 53.7\%                                 & \textbf{53.1\%}                        & \textbf{53.6\%}         & \textbf{53.6\%}         & \textbf{53.6\%}           \\
beasleyC3                                                                  & 79.9\%                                 & \textbf{79.2\%}                                 & \textbf{79.2\%}                        & \textbf{79.2\%}         & \textbf{79.2\%}                        & 79.5\%                                 & 79.5\%                  & 79.5\%                  & 79.5\%                    \\
binkar10\_1                                                                & \textbf{30.0\%}                        & 30.1\%                                          & 30.2\%                                 & 30.2\%                  & 30.2\%                                 & \textbf{29.9\%}                        & 30.2\%                  & 30.2\%                  & 30.2\%                    \\
cmflsp50-24-8-8                                                            & \textbf{22.9\%}                        & 23.9\%                                          & 23.8\%                                 & 23.9\%                  & 23.9\%                                 & 25.8\%                                 & 25.6\%                  & 25.8\%                  & 26.2\%                    \\
cod105                                                                     & \textbf{59.3\%}                        & 59.5\%                                          & 59.5\%                                 & 59.5\%                  & 59.5\%                                 & 59.4\%                                 & 59.4\%                  & 59.4\%                  & 59.4\%                    \\
cost266-UUE                                                                & 37.1\%                                 & \textbf{36.4\%}                                 & \textbf{36.4\%}                        & \textbf{36.4\%}         & \textbf{36.4\%}                        & 38.4\%                                 & 38.4\%                  & 38.5\%                  & 38.5\%                    \\
dano3\_3                                                                   & \textbf{16}                            & 21                                              & 28                                     & 28                      & 28                                     & 21                                     & 28                      & 28                      & 28                        \\
dano3\_5                                                                   & 1705                                   & \textbf{1279}                                   & 2716                                   & 2716                    & 2716                                   & 1328                                   & 1810                    & 1810                    & 1810                      \\
exp-1-500-5-5                                                              & 60.5\%                                 & \textbf{57.9\%}                                 & \textbf{57.9\%}                        & \textbf{57.9\%}         & \textbf{57.9\%}                        & \textbf{57.9\%}                        & \textbf{57.9\%}         & \textbf{57.9\%}         & \textbf{57.9\%}           \\
glass-sc                                                                   & 34.5\%                                 & \textbf{34.1\%}                                 & \textbf{34.1\%}                        & \textbf{34.1\%}         & \textbf{34.1\%}                        & \textbf{34.1\%}                        & \textbf{34.1\%}         & \textbf{34.1\%}         & \textbf{34.1\%}           \\
glass4                                                                     & 85.1\%                                 & \textbf{75.0\%}                                 & \textbf{75.0\%}                        & \textbf{75.0\%}         & \textbf{77.0\%}                        & 85.0\%                                 & 80.0\%                  & \textbf{75.0\%}         & 82.5\%                    \\
gmu-35-40                                                                  & 100.0\%                                & 50.1\%                                          & 50.0\%                                 & 50.0\%                  & 50.0\%                                 & \textbf{48.1\%}                        & \textbf{48.0\%}         & \textbf{48.0\%}         & \textbf{48.0\%}           \\
gmu-35-50                                                                  & 91.8\%                                 & \textbf{86.6\%}                                 & \textbf{86.5\%}                        & \textbf{86.5\%}         & \textbf{86.5\%}                        & 87.8\%                                 & 87.5\%                  & 87.5\%                  & 87.5\%                    \\
graph20-20-1rand                                                           & 93.0\%                                 & \textbf{92.6\%}                                 & \textbf{92.6\%}                        & \textbf{92.6\%}         & \textbf{92.6\%}                        & \textbf{92.6\%}                        & \textbf{92.6\%}         & \textbf{92.6\%}         & \textbf{92.6\%}           \\
lotsize                                                                    & 92.2\%                                 & \textbf{91.9\%}                                 & \textbf{91.9\%}                        & \textbf{91.9\%}         & \textbf{91.9\%}                        & \textbf{91.9\%}                        & \textbf{91.9\%}         & \textbf{91.9\%}         & \textbf{91.9\%}           \\
map10                                                                      & \textbf{53.5\%}                        & 64.0\%                                          & 64.0\%                                 & 64.0\%                  & 64.0\%                                 & 65.7\%                                 & 65.7\%                  & 65.7\%                  & 65.7\%                    \\
mc11                                                                       & 81.2\%                                 & \textbf{81.1\%}                                 & \textbf{81.1\%}                        & \textbf{81.1\%}         & \textbf{81.1\%}                        & \textbf{81.1\%}                        & \textbf{81.1\%}         & \textbf{81.1\%}         & \textbf{81.1\%}           \\
mcsched                                                                    & 22.4\%                                 & 17.8\%                                          & 19.0\%                                 & 19.0\%                  & 19.1\%                                 & \textbf{15.8\%}                        & \textbf{16.3\%}         & \textbf{16.3\%}         & \textbf{16.3\%}           \\
milo-v12-6-r2-40-1                                                         & \textbf{54.5\%}                        & 64.5\%                                          & 64.6\%                                 & 64.6\%                  & 64.4\%                                 & 60.8\%                                 & 60.7\%                  & 60.8\%                  & 60.8\%                    \\
neos-3627168-kasai                                                         & 20.5\%                                 & 19.7\%                                          & 19.8\%                                 & 19.8\%                  & 19.7\%                                 & \textbf{19.6\%}                        & \textbf{19.6\%}         & \textbf{19.6\%}         & \textbf{19.6\%}           \\
neos-5188808-nattai                                                        & \phantom{18584}53.7\% & \phantom{18584}\textbf{34.0\%} & \phantom{18584}39.5\% & \textbf{18845(3.6\%)}   & \phantom{18584}16.9\% & \phantom{18584}39.8\% & \textbf{19230(15.8\%)}  & 18584(8.7\%)            & \textbf{19537(10.6\%)}    \\
neos-848589                                                                & 95.6\%                                 & 95.6\%                                          & 95.6\%                                 & 95.6\%                  & 95.6\%                                 & 95.6\%                                 & 95.6\%                  & 95.6\%                  & 95.6\%                    \\
neos-860300                                                                & 527                                    & \textbf{438}                                    & 473                                    & 455                     & \textbf{377}                           & 443                                    & \textbf{451}            & \textbf{402}            & 381                       \\
neos-911970                                                                & \textbf{11.5\%}                        & 12.4\%                                          & 12.9\%                                 & 12.2\%                  & 11.3\%                                 & 12.1\%                                 & 12.4\%                  & 12.4\%                  & 12.0\%                    \\
neos17                                                                     & 21.1\%                                 & \textbf{17.4\%}                                 & \textbf{17.5\%}                        & \textbf{17.5\%}         & \textbf{17.5\%}                        & 18.5\%                                 & 18.5\%                  & 18.5\%                  & 18.5\%                    \\
neos5                                                                      & \textbf{51.8\%}                        & 54.1\%                                          & 54.1\%                                 & 54.1\%                  & 54.1\%                                 & 53.8\%                                 & 53.8\%                  & 53.8\%                  & 53.8\%                    \\
ns1830653                                                                  & 26.5\%                                 & 23.0\%                                          & 22.5\%                                 & \textbf{22.4\%}         & \textbf{22.4\%}                        & \textbf{22.0\%}                        & \textbf{22.1\%}         & 22.5\%                  & 23.0\%                    \\
p200x1188c                                                                 & 56.1\%                                 & 49.8\%                                          & 49.8\%                                 & 49.8\%                  & 49.8\%                                 & \textbf{49.7\%}                        & \textbf{49.7\%}         & \textbf{49.7\%}         & \textbf{49.7\%}           \\
pg                                                                         & 61.8\%                                 & \textbf{60.6\%}                                 & \textbf{60.6\%}                        & \textbf{60.6\%}         & \textbf{60.6\%}                        & \textbf{60.6\%}                        & \textbf{60.6\%}         & \textbf{60.6\%}         & \textbf{60.6\%}           \\
pg5\_34                                                                    & 48.9\%                                 & \textbf{48.6\%}                                 & \textbf{48.6\%}                        & \textbf{48.6\%}         & \textbf{48.6\%}                        & \textbf{48.6\%}                        & \textbf{48.6\%}         & \textbf{48.6\%}         & \textbf{48.6\%}           \\
pk1                                                                        & 89.6\%                                 & 88.9\%                                          & 88.9\%                                 & 88.9\%                  & 88.9\%                                 & \textbf{88.3\%}                        & \textbf{88.3\%}         & \textbf{88.3\%}         & \textbf{88.3\%}           \\
ran14x18-disj-8                                                            & \textbf{58.0\%}                        & 60.7\%                                          & 60.6\%                                 & 60.6\%                  & 60.6\%                                 & 59.4\%                                 & 59.4\%                  & 59.4\%                  & 59.4\%                    \\
rd-rplusc-21                                                               & 100.0\%                                & 100.0\%                                         & 100.0\%                                & 100.0\%                 & 100.0\%                                & 100.0\%                                & 100.0\%                 & 100.0\%                 & 100.0\%                   \\
rmatr100-p10                                                               & 1484                                   & \textbf{1255}                                   & 1349                                   & 1460                    & 1466                                   & 1261                                   & \textbf{1323}           & \textbf{1412}           & \textbf{1457}             \\
seymour                                                                    & \textbf{49.5\%}                        & 51.8\%                                          & 51.8\%                                 & 51.8\%                  & 51.8\%                                 & 51.8\%                                 & 51.8\%                  & 51.8\%                  & 51.8\%                    \\
seymour1                                                                   & \textbf{3250}                          & 4106                                            & 4346                                   & 4361                    & 4361                                   & 4102                                   & 4237                    & 4237                    & 4237                      \\
sp150x300d                                                                 & \textbf{9.1\%}                         & 9.2\%                                           & 9.2\%                                  & 9.2\%                   & 9.3\%                                  & 9.2\%                                  & 9.2\%                   & 9.2\%                   & 9.2\%                     \\
supportcase26                                                              & 69.4\%                                 & 60.7\%                                          & 60.5\%                                 & 61.0\%                  & 62.1\%                                 & \textbf{60.2\%}                        & \textbf{61.0\%}         & \textbf{60.5\%}         & \textbf{60.6\%}           \\
supportcase40                                                              & 27.7\%                                 & \textbf{15.2\%}                                 & 18.9\%                                 & 18.1\%                  & 18.2\%                                 & 15.7\%                                 & \textbf{17.3\%}         & \textbf{17.3\%}         & \textbf{17.3\%}           \\
supportcase7                                                               & 605                                    & \textbf{231}                                    & 346                                    & 351                     & 331                                    & 233                                    & \textbf{345}            & \textbf{348}            & \textbf{328}              \\
tr12-30                                                                    & 58.5\%                                 & \textbf{57.8\%}                                 & \textbf{57.8\%}                        & \textbf{57.8\%}         & \textbf{57.8\%}                        & \textbf{57.8\%}                        & \textbf{57.8\%}         & \textbf{57.8\%}         & \textbf{57.8\%}           \\
uct-subprob                                                                & \textbf{62.1\%}                        & 63.8\%                                          & 63.8\%                                 & 63.8\%                  & 63.8\%                                 & 63.2\%                                 & 63.2\%                  & 63.2\%                  & 63.1\%                    \\ \hline
\multicolumn{1}{|c|}{Geo mean -- Nodes\footnotemark[1]} & 943                         & \textbf{813}            & 968                     & 968                     & 950                       & 819                     & \textbf{902}            & \textbf{894}            & \textbf{893}              \\ \hline
\multicolumn{1}{|c|}{Geo mean -- Gap\footnotemark[2]}   & 49.1\%                      & \textbf{46.2\%}         & 46.8\%                  & \textbf{44.0\%}         & 45.6\%                    & 46.5\%                  & \textbf{45.5\%}         & 44.7\%                  & \textbf{45.1\%}           \\ \hline
\end{tabular}
}}
\end{center}
\end{table}

We present similar tables for RB as those for FSB. Five instances showed no change in dual bound after 20,000 nodes with any of the rules and are excluded from the results. {In addition, the solution of two instances could not be completed in 14 days and are not included in the results.} Table~\ref{tab:miplib_rb} compares Eff-RB~(3,7) against Def-RB on the remaining 44 instances. 
\begin{enumerate}
    % \item Here again we see improvements in the number of nodes as well as the gap remaining at termination. 
   {
    \item On the 8 instances solved within the node limit in any run, the mean number of nodes of Eff-RB~(3,7) decreases by 13.8\% when using optimal primal bounds, but increases by 2.7\%,  2.6\% and 0.7\% when the primal gap is 2\%, 5\% and 10\%, respectively. On these instances, PA-LA-RB~(3,7,0.15) lead to a decrease in mean tree sizes for all primal gaps. The number of nodes reduces by 13.2\%, 4.4\%, 5.2\% and 5.3\% respectively with the four primal gaps relative to Def-SB.
    \item On the 37 instances that could not be solved within the node limit by some run, the mean gap remaining decreases by 2.9\%, 2.3\%, 5.1\% and 3.5\% with Eff-RB~(3,7) and by 2.6\%, 3.6\%, 4.3\% and 4.0\% with PA-LA-RB~(3,7,0.15) for the different primal gaps. 
    \item These results again indicate a more robust performance of PA-LA-RB~(3,7,0.15) as compared to Eff-RB~(3,7).
    }
\end{enumerate}

\section{Conclusion and Future Directions}\label{sec:conclusion}

% In this work, we have tried to better understand the properties and shortcomings of FSB. We have identified some factors that can significantly impact tree sizes and used them to design new FSB scores. 

% Parameters selected that work well on average, but the results on random instances show that the best parameters differ based on problem structure. 
% For efficacious gains - (3,7) works well on most problems. But the point of minima is often different, eg. (1,9) for independent set and matching and (4,6) for combinatorial auction and set covering. 
% If the optimal parameters could be predicted or estimated for a problem class, that would make the proposed approached even more powerful. 

% Determining the right parameters is even more relevant for LA, where the differences in trends are even more pronounced for different problems. While the chosen 0.15 parameters leads to more wins than losses and better average performance, the gains are not as consistent as seen with Eff due to this reason. Our results indicate that asymmetrical trend in infeasibility is an important consideration, but LA may not be the best approach to detect and correct them. Future research directions can investigate other metrics for capturing asymmetry due to infeasibility. A more general approach towards rebalancing trees that also incorporates possible asymmetry in LP gains can be explored. 
{
In this work, we have attempted to better understand the properties and limitations of FSB. We have identified key factors that significantly impact tree sizes and have leveraged them to design new FSB scores. Through computational experiments on problems with known structures, we have studied how these scores impact branching decisions and have tuned them to optimize the average performance. We finally evaluate the proposed scores on instances from MIPLIB 2017 Benchmark set and demonstrate a remarkable improvement in performance over FSB. 
Our new rules extend to RB and offer significant performance improvement, albeit less than in the case of FSB.  
This performance difference may be explained by the fact that our approach assumes that LP gains are accurate, which is not always true for pseudo-cost-based estimates.

While the selected parameters perform well on average, results on random instances reveal that the optimal parameters vary depending on the problem structure. 
In the case of efficacious gains, $(a_{\min}, a_{\max})=(0.3,0.7)$ works well for most problems, but $(0.1,0.9)$ yields smaller trees on independent set and matching, and $(0.4,0.6)$ is the preferred choice for combinatorial auction and set covering. If the optimal parameters could be predicted or estimated for a given problem class, the proposed approach could be even more powerful. 

\textcolor{black}{
Moreover, the preferred values of the exponents $(a_{\min}, a_{\max})=(0.3,0.7)$ decrease as primal bound improves. Therefore, a more sophisticated implementation of this approach may update the exponents based on the incumbent primal gap; using parameters $(0.5, 0.5)$ when the primal bound is unknown and gradually transitioning up to $(0.3, 0.7)$ as the primal gap decreases. 
}
% \red{change 3,7 -> 4,6 based on primal gap - outside scope of paper}

Determining good parameters is even more crucial for the last assignment based rules, where trends vary more significantly across different problem types. Although the chosen parameter $\eta=0.15$ leads to more wins than losses and improves average performance, its effectiveness is weakened due to these variations. Our results indicate that asymmetry plays a critical role in many problems, yet the last assignment rule may not be the most effective approach for detecting and addressing it. Future research could explore alternative metrics for capturing asymmetrical pruning. 
% \red{A more general strategy for rebalancing trees—one that accounts for asymmetry in infeasibility and potential asymmetry in LP gains—could be a promising direction for further investigation.} 
\textcolor{black}{Furthermore, known problem structure can be exploited to design specialized scores. An example of such a specialized score for cardinality constrained problems is presented in Appendix~\ref{sec:cardinality_constr}. } 

\textcolor{black}{The proposed approaches based on last assignments assume that variables are similar and therefore the global trend in asymmetry can be applied to all candidate variables at the current node. Future research may extend this approach to design variable specific score to account for dissimilar variables in any problem. If information regarding the problem structure is known, different scores may be chosen for different groups of variables.}

% \red{LA with best-bound rule - search the tree in a strange way, it may give different results. But unlikely. }

}
% \begin{itemize}
%     \item We have tried to understand the properties and shortcomings of FSB better 
    
%     % \item Theoretical side: most important is the role of infeasibility vs bounds 

%     \item We have identified some factors that can have a significant impact on tree sizes and have used them to design new FSB scores.
    
%     \item Parameters selected that work well on average, but the results show that best parameters differ based on problem structure. 
%     For efficacious gains - 3,7 works well on most problems. But the point of minima is often different, eg. 1, 9 for independent set and matching and 4,6 for combinatorial auction and set covering. Even more relevant for LA where the differences in trends are even more pronounced for different problems. 
%     If the optimal parameters could be predicted or estimated for a problem class, that would make the proposed approached even more powerful. 

%     LA - more victories than losses. How to 

%     \item Open directions - other metrics for asymmetry, more general rebalancing.

% \end{itemize}

\paragraph{Acknowledgements} 
\textcolor{black}{We are grateful to Nikolaos Sahinidis for his support in performing the computational experiments in this paper. We also thank the support from AFOSR grant \# FA9550-22-1-0052.}

\bibliographystyle{plain}
\bibliography{bibliography}

\begin{thebibliography}{10}

\bibitem{achterberg2007constraint}
Tobias Achterberg.
\newblock {\em Constraint integer programming}.
\newblock PhD thesis, 2007.

\bibitem{achterberg2009hybrid}
Tobias Achterberg and Timo Berthold.
\newblock Hybrid branching.
\newblock In {\em Integration of AI and OR Techniques in Constraint Programming for Combinatorial Optimization Problems: 6th International Conference, CPAIOR 2009 Pittsburgh, PA, USA, May 27-31, 2009 Proceedings 6}, pages 309--311. Springer, 2009.

\bibitem{achterberg2005branching}
Tobias Achterberg, Thorsten Koch, and Alexander Martin.
\newblock Branching rules revisited.
\newblock {\em Operations Research Letters}, 33(1):42--54, 2005.

\bibitem{alvarez2017machine}
Alejandro~Marcos Alvarez, Quentin Louveaux, and Louis Wehenkel.
\newblock A machine learning-based approximation of strong branching.
\newblock {\em INFORMS Journal on Computing}, 29(1):185--195, 2017.

\bibitem{applegate1995finding}
David Applegate, Robert Bixby, Va{\v{s}}ek Chv{\'a}tal, and William Cook.
\newblock {\em Finding cuts in the TSP (A preliminary report)}, volume~95.
\newblock Citeseer, 1995.

\bibitem{benichou1971experiments}
Michel B{\'e}nichou, Jean-Michel Gauthier, Paul Girodet, Gerard Hentges, Gerard Ribi{\`e}re, and Olivier Vincent.
\newblock Experiments in mixed-integer linear programming.
\newblock {\em Mathematical programming}, 1:76--94, 1971.

\bibitem{bergman2016decision}
David Bergman, Andre~A Cire, Willem-Jan Van~Hoeve, and John Hooker.
\newblock {\em Decision diagrams for optimization}, volume~1.
\newblock Springer, 2016.

\bibitem{berthold2006primal}
Timo Berthold.
\newblock {\em Primal heuristics for mixed integer programs}.
\newblock PhD thesis, Zuse Institute Berlin (ZIB), 2006.

\bibitem{berthold2013measuring}
Timo Berthold.
\newblock Measuring the impact of primal heuristics.
\newblock {\em Operations Research Letters}, 41(6):611--614, 2013.

\bibitem{berthold1945primal}
Timo Berthold, Andrea Lodi, and Domenico Salvagnin.
\newblock Primal heuristics in integer programming.
\newblock {\em American history}, 1861(1900), 1945.

\bibitem{BestuzhevaEtal2021OO}
Ksenia Bestuzheva, Mathieu Besan{\c{c}}on, Wei-Kun Chen, Antonia Chmiela, Tim Donkiewicz, Jasper van Doornmalen, Leon Eifler, Oliver Gaul, Gerald Gamrath, Ambros Gleixner, Leona Gottwald, Christoph Graczyk, Katrin Halbig, Alexander Hoen, Christopher Hojny, Rolf van~der Hulst, Thorsten Koch, Marco L{\"u}bbecke, Stephen~J. Maher, Frederic Matter, Erik M{\"u}hmer, Benjamin M{\"u}ller, Marc~E. Pfetsch, Daniel Rehfeldt, Steffan Schlein, Franziska Schl{\"o}sser, Felipe Serrano, Yuji Shinano, Boro Sofranac, Mark Turner, Stefan Vigerske, Fabian Wegscheider, Philipp Wellner, Dieter Weninger, and Jakob Witzig.
\newblock {The SCIP Optimization Suite 8.0}.
\newblock Technical report, Optimization Online, December 2021.

\bibitem{cornuejols1991comparison}
G{\'e}rard Cornu{\'e}jols, Ranjani Sridharan, and Jean-Michel Thizy.
\newblock A comparison of heuristics and relaxations for the capacitated plant location problem.
\newblock {\em European journal of operational research}, 50(3):280--297, 1991.

\bibitem{dey2023lower}
Santanu~S Dey, Yatharth Dubey, and Marco Molinaro.
\newblock Lower bounds on the size of general branch-and-bound trees.
\newblock {\em Mathematical Programming}, 198(1):539--559, 2023.

\bibitem{dey2024theoretical}
Santanu~S Dey, Yatharth Dubey, Marco Molinaro, and Prachi Shah.
\newblock A theoretical and computational analysis of full strong-branching.
\newblock {\em Mathematical Programming}, 205(1):303--336, 2024.

\bibitem{fischetti2011backdoor}
Matteo Fischetti and Michele Monaci.
\newblock Backdoor branching.
\newblock In {\em International Conference on Integer Programming and Combinatorial Optimization}, pages 183--191. Springer, 2011.

\bibitem{gamrath2018measuring}
Gerald Gamrath and Christoph Schubert.
\newblock Measuring the impact of branching rules for mixed-integer programming.
\newblock In {\em Operations Research Proceedings 2017: Selected Papers of the Annual International Conference of the German Operations Research Society (GOR), Freie Universi{\"a}t Berlin, Germany, September 6-8, 2017}, pages 165--170. Springer, 2018.

\bibitem{gasse2019exact}
Maxime Gasse, Didier Ch{\'e}telat, Nicola Ferroni, Laurent Charlin, and Andrea Lodi.
\newblock Exact combinatorial optimization with graph convolutional neural networks.
\newblock {\em Advances in neural information processing systems}, 32, 2019.

\bibitem{gilpin2006information}
Andrew Gilpin and Tuomas Sandholm.
\newblock Information-theoretic approaches to branching in search.
\newblock In {\em Proceedings of the fifth international joint conference on Autonomous agents and multiagent systems}, pages 545--547, 2006.

\bibitem{glankwamdee2006lookahead}
Wasu Glankwamdee and Jeff Linderoth.
\newblock Lookahead branching for mixed integer programming.
\newblock In {\em Twelfth INFORMS Computing Society Meeting}, pages 130--150, 2006.

\bibitem{gleeson1990identifying}
John Gleeson and Jennifer Ryan.
\newblock Identifying minimally infeasible subsystems of inequalities.
\newblock {\em ORSA Journal on Computing}, 2(1):61--63, 1990.

\bibitem{miplib2017}
Ambros Gleixner, Gregor Hendel, Gerald Gamrath, Tobias Achterberg, Michael Bastubbe, Timo Berthold, Philipp~M. Christophel, Kati Jarck, Thorsten Koch, Jeff Linderoth, Marco L\"ubbecke, Hans~D. Mittelmann, Derya Ozyurt, Ted~K. Ralphs, Domenico Salvagnin, and Yuji Shinano.
\newblock {MIPLIB 2017: Data-Driven Compilation of the 6th Mixed-Integer Programming Library}.
\newblock {\em Mathematical Programming Computation}, 2021.

\bibitem{gupta2020hybrid}
Prateek Gupta, Maxime Gasse, Elias Khalil, Pawan Mudigonda, Andrea Lodi, and Yoshua Bengio.
\newblock Hybrid models for learning to branch.
\newblock {\em Advances in neural information processing systems}, 33:18087--18097, 2020.

\bibitem{gurobi}
{Gurobi Optimization, LLC}.
\newblock {Gurobi Optimizer Reference Manual}, 2024.

\bibitem{khalil2016learning}
Elias Khalil, Pierre Le~Bodic, Le~Song, George Nemhauser, and Bistra Dilkina.
\newblock Learning to branch in mixed integer programming.
\newblock In {\em Proceedings of the AAAI conference on artificial intelligence}, volume~30, 2016.

\bibitem{khalil2022finding}
Elias~B Khalil, Pashootan Vaezipoor, and Bistra Dilkina.
\newblock Finding backdoors to integer programs: a monte carlo tree search framework.
\newblock In {\em Proceedings of the AAAI Conference on Artificial Intelligence}, volume~36, pages 3786--3795, 2022.

\bibitem{kilincc2009information}
Fatma K{\i}l{\i}n{\c{c}}~Karzan, George~L Nemhauser, and Martin~WP Savelsbergh.
\newblock Information-based branching schemes for binary linear mixed integer problems.
\newblock {\em Mathematical Programming Computation}, 1:249--293, 2009.

\bibitem{land60a}
Ailsa~H. Land and Alison~G. Doig.
\newblock An automatic method of solving discrete programming problems.
\newblock {\em Econometrica}, 28(3):497--520, 1960.

\bibitem{le2015important}
Pierre Le~Bodic and George~L Nemhauser.
\newblock How important are branching decisions: Fooling mip solvers.
\newblock {\em Operations Research Letters}, 43(3):273--278, 2015.

\bibitem{leyton2000towards}
Kevin Leyton-Brown, Mark Pearson, and Yoav Shoham.
\newblock Towards a universal test suite for combinatorial auction algorithms.
\newblock In {\em Proceedings of the 2nd ACM conference on Electronic commerce}, pages 66--76, 2000.

\bibitem{li1997look}
Chu~Min Li and Anbulagan.
\newblock Look-ahead versus look-back for satisfiability problems.
\newblock In {\em International Conference on Principles and Practice of Constraint Programming}, pages 341--355. Springer, 1997.

\bibitem{linderoth1999computational}
Jeff~T Linderoth and Martin~WP Savelsbergh.
\newblock A computational study of search strategies for mixed integer programming.
\newblock {\em INFORMS Journal on Computing}, 11(2):173--187, 1999.

\bibitem{moskewicz2001chaff}
Matthew~W Moskewicz, Conor~F Madigan, Ying Zhao, Lintao Zhang, and Sharad Malik.
\newblock Chaff: Engineering an efficient sat solver.
\newblock In {\em Proceedings of the 38th annual Design Automation Conference}, pages 530--535, 2001.

\bibitem{nair2020solving}
Vinod Nair, Sergey Bartunov, Felix Gimeno, Ingrid Von~Glehn, Pawel Lichocki, Ivan Lobov, Brendan O'Donoghue, Nicolas Sonnerat, Christian Tjandraatmadja, Pengming Wang, et~al.
\newblock Solving mixed integer programs using neural networks.
\newblock {\em arXiv preprint arXiv:2012.13349}, 2020.

\bibitem{nemhauser1988integer}
George~L. Nemhauser and Laurence~A. Wolsey.
\newblock {\em Integer and Combinatorial Optimization}.
\newblock Wiley-Interscience Series in Discrete Mathematics and Optimization. John Wiley \& Sons, 1988.

\bibitem{pagnoncelli2009computational}
Bernardo~K Pagnoncelli, Shabbir Ahmed, and Alexander Shapiro.
\newblock Computational study of a chance constrained portfolio selection problem.
\newblock {\em Journal of Optimization Theory and Applications}, 142(2):399--416, 2009.

\bibitem{patel2007active}
Jagat Patel and John~W Chinneck.
\newblock Active-constraint variable ordering for faster feasibility of mixed integer linear programs.
\newblock {\em Mathematical Programming}, 110:445--474, 2007.

\bibitem{prouvost2020ecole}
Antoine Prouvost, Justin Dumouchelle, Lara Scavuzzo, Maxime Gasse, Didier Ch{\'e}telat, and Andrea Lodi.
\newblock Ecole: A gym-like library for machine learning in combinatorial optimization solvers.
\newblock In {\em Learning Meets Combinatorial Algorithms at NeurIPS2020}, 2020.

\bibitem{qiu2014covering}
Feng Qiu, Shabbir Ahmed, Santanu~S Dey, and Laurence~A Wolsey.
\newblock Covering linear programming with violations.
\newblock {\em INFORMS Journal on Computing}, 26(3):531--546, 2014.

\bibitem{quadt2008capacitated}
Daniel Quadt and Heinrich Kuhn.
\newblock Capacitated lot-sizing with extensions: a review.
\newblock {\em 4OR}, 6(1):61--83, 2008.

\bibitem{shah2025non}
Prachi Shah, Santanu~S Dey, and Marco Molinaro.
\newblock Non-monotonicity of branching rules with respect to linear relaxations.
\newblock {\em INFORMS Journal on Computing}, 2025.

\bibitem{yang2022learning}
Yu~Yang, Natashia Boland, Bistra Dilkina, and Martin Savelsbergh.
\newblock Learning generalized strong branching for set covering, set packing, and 0--1 knapsack problems.
\newblock {\em European Journal of Operational Research}, 301(3):828--840, 2022.

\bibitem{zarpellon2021parameterizing}
Giulia Zarpellon, Jason Jo, Andrea Lodi, and Yoshua Bengio.
\newblock Parameterizing branch-and-bound search trees to learn branching policies.
\newblock In {\em Proceedings of the aaai conference on artificial intelligence}, volume~35, pages 3931--3939, 2021.

\end{thebibliography}
\newpage

\appendix

\section{Alternative to Efficacious Gains} \label{sec:Greedy pruning}
An alternative approach to solve the issue posed by overestimated LP gains in Example~\ref{ex:1} could be to identify pruning by bounds and treat them as the same as pruning by infeasibility. In other words, any variable that leads to a child that can be pruned either by the incumbent primal bound or infeasibility gets an infinite LP gain. If multiple branching candidates lead to pruned children, branch on the variable with two pruned nodes if such exists, else, maximize the gain on the remaining branch. 

We call this variant of FSB the Pruning Focused rule and evaluate its tree sizes against Def-SB on the test bed of random problems. The results are presented as shifted geometric means of tree sizes in Table~\ref{tab:greedy pruning} and as percentage reduction of the same values relative to Def-SB in Table~\ref{tab:greedy pruning pct}. When the optimal primal bound is known, this rule leads to a 15.7\% reduction in overall mean tree sizes relative to Def-SB, surpassing Eff-SB. However, this rule is not robust to suboptimal primal bounds, and the tree sizes are only 2.4\% smaller than Def-SB with a 2\% gap in primal bounds. This sudden change can be explained by the fact that the variable scores as a function of the primal bound have discontinuous jumps between finite and infinite values, as demonstrated by the following example. 

\begin{example}
Consider the branching decision in Example~\ref{ex:1}. If we have a primal solution such that $\Delta_{p-d} = 1$, both variables have an infinite score, and $x_j$ is preferred for branching since it has a higher gain on the unpruned branch. But with a worse primal solution resulting in $\Delta_{p-d} = 1 + \epsilon$ for any small $\epsilon > 0$, $x_j$ get a finite score of 0.9, while the score of $x_i$ remains infinite and is chosen for branching.
\end{example}

Moreover, the results in Table~\ref{tab:greedy pruning} indicate that improving the primal bound can often lead to worse decisions, as on 9 of the 14 problems, the tree sizes are smaller for 10\% primal gap than for 2\% primal gap. 

\begin{table}
\begin{center}
		
		\caption{Shifted geometric mean of tree sizes for the  Def-SB and Pruning Focused rules rounded to integers. Smaller tree sizes are highlighted in bold to compare the differences between different rules. }\label{tab:greedy pruning}
\resizebox{1.\textwidth}{!}{
\setlength{\tabcolsep}{0.5em} 
  {\renewcommand{\arraystretch}{1.1}
\begin{tabular}{|l|r|rrrr|rrrr|}
\hline
Rule & \multicolumn{1}{c|}{Def-SB} & \multicolumn{4}{c|}{Pruning focussed} & \multicolumn{4}{c|}{Pruning focussed (3,7)} \\ \hline
Primalgap & \multicolumn{1}{c|}{} & \multicolumn{1}{c}{0\%} & \multicolumn{1}{c}{2\%} & \multicolumn{1}{c}{5\%} & \multicolumn{1}{c|}{10\%} & \multicolumn{1}{r}{0\%} & \multicolumn{1}{r}{2\%} & \multicolumn{1}{r}{5\%} & \multicolumn{1}{r|}{10\%} \\ \hline
\multicolumn{1}{|c|}{Problem} &  & \multicolumn{1}{r}{} & \multicolumn{1}{r}{} & \multicolumn{1}{r}{} & \multicolumn{1}{r|}{} & \multicolumn{1}{r}{} & \multicolumn{1}{r}{} & \multicolumn{1}{r}{} & \multicolumn{1}{r|}{} \\ \hline
Big-Bucket Lot-Sizing (30) & 3099 & \textbf{2798} & 3141 & 3121 & \textbf{3013} & \textbf{2653} & \textbf{2989} & \textbf{2920} & \textbf{2848} \\
Capacitated Lot-Sizing (30) & 5302 & \textbf{4745} & 5792 & 5442 & \textbf{5264} & \textbf{4860} & 5653 & 5416 & 5359 \\
Fixed Charge Flow (150) & 7658 & \textbf{6389} & 7865 & \textbf{7546} & \textbf{7219} & \textbf{6137} & \textbf{7249} & \textbf{6982} & \textbf{6793} \\
Capacited Facility Location (100) & 3733 & \textbf{3566} & \textbf{3648} & \textbf{3648} & \textbf{3648} & \textbf{3499} & \textbf{3586} & \textbf{3586} & \textbf{3586} \\
Portfolio Optimization (100) & 19965 & 20877 & 22591 & 21735 & 21719 & 21711 & 23604 & 22620 & 22584 \\
Weighted Set Coverage (200) & 1274 & 1335 & 1368 & 1342 & 1342 & 1597 & 1768 & 1732 & 1732 \\
Set Covering (3000) & 2232 & 2261 & 2878 & 2966 & 2408 & 2447 & 3057 & 3559 & 2914 \\
Set Packing (1000) & 5362 & \textbf{3626} & \textbf{3682} & \textbf{3777} & \textbf{4096} & \textbf{3557} & \textbf{3640} & \textbf{3765} & \textbf{4036} \\
Combinatorial Auctions (200) & 1493 & \textbf{1304} & 1504 & \textbf{1482} & \textbf{1482} & \textbf{1384} & 2046 & 1952 & 1952 \\
Independent Set (500) & 1701 & \textbf{854} & \textbf{1364} & \textbf{1427} & \textbf{1441} & \textbf{856} & \textbf{1281} & \textbf{1301} & \textbf{1318} \\
Weighted Matching (1000) & 3545 & \textbf{2726} & \textbf{3184} & \textbf{3184} & \textbf{3184} & \textbf{2253} & \textbf{2541} & \textbf{2541} & \textbf{2541} \\
Multi-dim Knapsack - Small (100) & 1014 & \textbf{916} & \textbf{997} & \textbf{996} & \textbf{995} & \textbf{889} & \textbf{961} & \textbf{960} & \textbf{960} \\
Multi-dim Knapsack - Med (100) & 16196 & \textbf{13455} & \textbf{15508} & \textbf{15190} & \textbf{15174} & \textbf{12423} & \textbf{13703} & \textbf{13588} & \textbf{13585} \\
Multi-dim Knapsack - Large (100) & 18171 & \textbf{13140} & \textbf{15825} & \textbf{15697} & \textbf{15697} & \textbf{11102} & \textbf{12217} & \textbf{12195} & \textbf{12195} \\ \hline
Geo Mean & 4149 & \textbf{3498} & \textbf{4050} & \textbf{4017} & \textbf{3951} & \textbf{3442} & \textbf{3981} & \textbf{3974} & \textbf{3922} \\ \hline
\end{tabular}
}}
\end{center}
\end{table}

{The results in Section~\ref{sec:MIPLIB} indicate that a higher exponent for the larger gain may improve performance in the absence of overestimation. Therefore, we also test the Pruning Focused rule that uses the Def-SB~(3,7) score whenever no branching variable leads to pruning. We refer to this as Pruning Focused~(3,7) in Tables~\ref{tab:greedy pruning} and ~\ref{tab:greedy pruning pct}. This rule leads to a remarkable decrease in tree sizes for many problems such as set-packing, independent set and weighted matching. On the other hand, other problems, such as set-covering, combinatorial auctions and weighted set coverage, see a substantial increase in tree size. Overall, the mean tree sizes improve over the Pruning Focused rule for all primal gaps. With optimal primal information, the mean tree sizes are 17\% less than Def-SB and better than Eff-SB~(3,7). But these gains again fall to 4-5\% with suboptimal bounds which are much lower than the same for Eff-SB~(3,7). Therefore, this rule is less robust.}

\begin{table}
\begin{center}
		
		\caption{Comparison of shifted geometric mean of tree sizes for the  Def-SB and Pruning Focused rule. Tree sizes are presented for Def-SB and the percentage reduction in tree size relative to Def-SB are presented for others.}\label{tab:greedy pruning pct}
\resizebox{1.\textwidth}{!}{
\setlength{\tabcolsep}{0.7em} 
  {\renewcommand{\arraystretch}{1.1}
\begin{tabular}{|l|r|rrrr|rrrr|}
\hline
Rule & \multicolumn{1}{c|}{Def-SB} & \multicolumn{4}{c|}{Pruning focussed} & \multicolumn{4}{c|}{Pruning focussed (3,7)} \\ \hline
Primalgap & \multicolumn{1}{c|}{} & \multicolumn{1}{c}{0\%} & \multicolumn{1}{c}{2\%} & \multicolumn{1}{c}{5\%} & \multicolumn{1}{c|}{10\%} & 0\% & 2\% & 5\% & 10\% \\ \hline
\multicolumn{1}{|c|}{Problem} &  &  &  &  &  &  &  &  &  \\ \hline
Big-Bucket Lot-Sizing (30) & 3099 & \textbf{10\%} & -1\% & -1\% & \textbf{3\%} & \textbf{14\%} & \textbf{4\%} & \textbf{6\%} & \textbf{8\%} \\
Capacitated Lot-Sizing (30) & 5302 & \textbf{11\%} & -9\% & -3\% & \textbf{1\%} & \textbf{8\%} & -7\% & -2\% & -1\% \\
Fixed Charge Flow (150) & 7658 & \textbf{17\%} & -3\% & \textbf{1\%} & \textbf{6\%} & \textbf{20\%} & \textbf{5\%} & \textbf{9\%} & \textbf{11\%} \\
Capacited Facility Location (100) & 3733 & \textbf{4\%} & \textbf{2\%} & \textbf{2\%} & \textbf{2\%} & \textbf{6\%} & \textbf{4\%} & \textbf{4\%} & \textbf{4\%} \\
Portfolio Optimization (100) & 19965 & -5\% & -13\% & -9\% & -9\% & -9\% & -18\% & -13\% & -13\% \\
Weighted Set Coverage (200) & 1274 & -5\% & -7\% & -5\% & -5\% & -25\% & -39\% & -36\% & -36\% \\
Set Covering (3000) & 2232 & -1\% & -29\% & -33\% & -8\% & -10\% & -37\% & -59\% & -31\% \\
Set Packing (1000) & 5362 & \textbf{32\%} & \textbf{31\%} & \textbf{30\%} & \textbf{24\%} & \textbf{34\%} & \textbf{32\%} & \textbf{30\%} & \textbf{25\%} \\
Combinatorial Auctions (200) & 1493 & \textbf{13\%} & -1\% & \textbf{1\%} & \textbf{1\%} & \textbf{7\%} & -37\% & -31\% & -31\% \\
Independent Set (500) & 1701 & \textbf{50\%} & \textbf{20\%} & \textbf{16\%} & \textbf{15\%} & \textbf{50\%} & \textbf{25\%} & \textbf{24\%} & \textbf{23\%} \\
Weighted Matching (1000) & 3545 & \textbf{23\%} & \textbf{10\%} & \textbf{10\%} & \textbf{10\%} & \textbf{36\%} & \textbf{28\%} & \textbf{28\%} & \textbf{28\%} \\
Multi-dim Knapsack - Small (100) & 1014 & \textbf{10\%} & \textbf{2\%} & \textbf{2\%} & \textbf{2\%} & \textbf{12\%} & \textbf{5\%} & \textbf{5\%} & \textbf{5\%} \\
Multi-dim Knapsack - Med (100) & 16196 & \textbf{17\%} & \textbf{4\%} & \textbf{6\%} & \textbf{6\%} & \textbf{23\%} & \textbf{15\%} & \textbf{16\%} & \textbf{16\%} \\
Multi-dim Knapsack - Large (100) & 18171 & \textbf{28\%} & \textbf{13\%} & \textbf{14\%} & \textbf{14\%} & \textbf{39\%} & \textbf{33\%} & \textbf{33\%} & \textbf{33\%} \\ \hline
Geo Mean & 4149 & \textbf{15.7\%} & \textbf{2.4\%} & \textbf{3.2\%} & \textbf{4.8\%} & \textbf{17.0\%} & \textbf{4.1\%} & \textbf{4.2\%} & \textbf{5.5\%} \\ \hline
\end{tabular}
}}
\end{center}
\end{table}

\section{Importance of Identifying Branching Constraints Necessary for Infeasibility} \label{sec:necc_branching_constrs}

\begin{figure}
    \centering
    \includegraphics[width=0.5\linewidth]{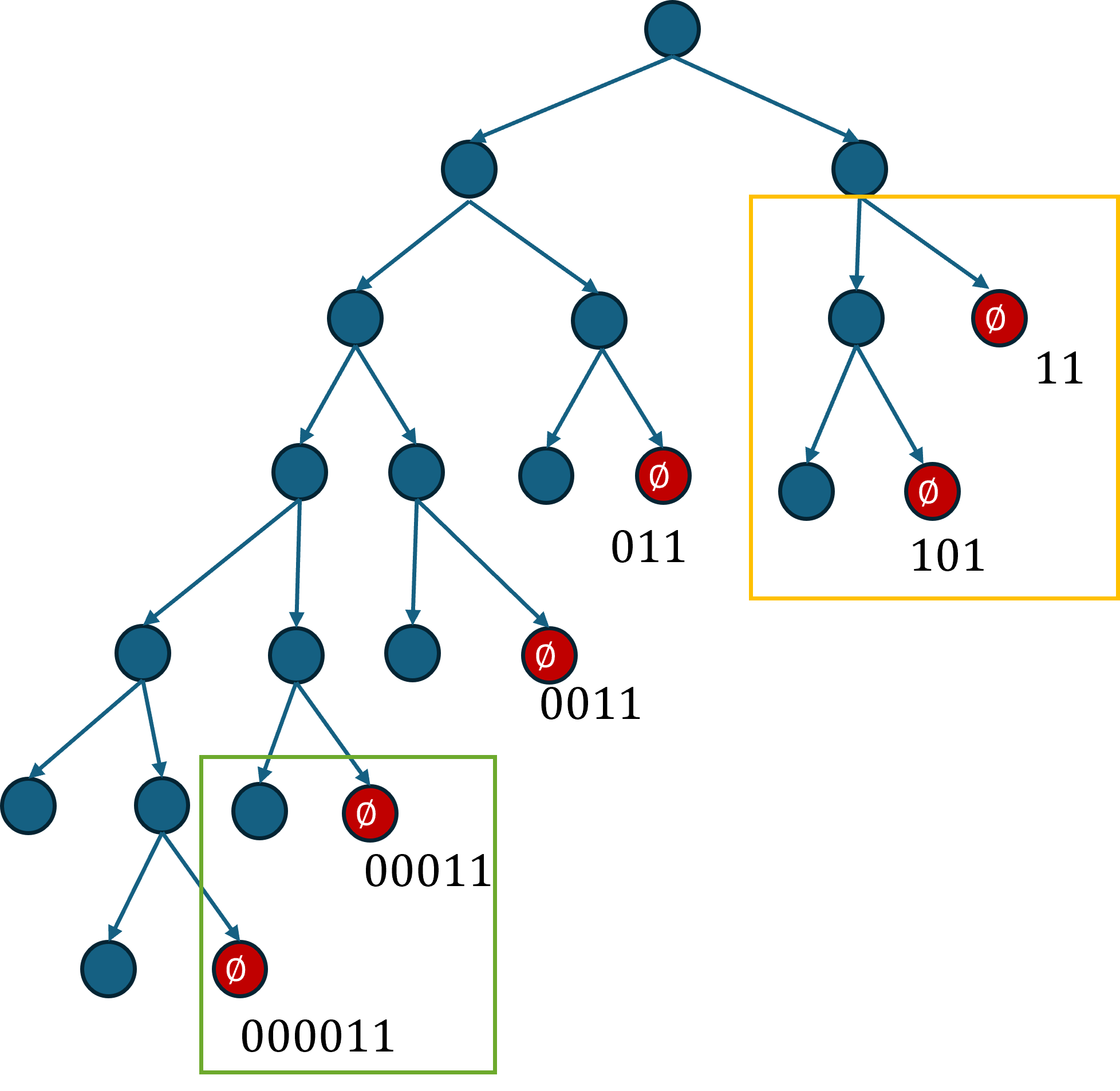}
    \caption{A possible partial branch-and-bound tree for the knapsack problem with constraint $\sum_i 2x_i \leq 3$.}
    \label{fig:packingtree}
\end{figure}

Consider a knapsack problem with constraint $\sum_i 2x_i \leq 3$. A possible partial branch-and-bound tree for this problem is illustrated in Figure~\ref{fig:packingtree}, where the infeasible nodes are marked in red and the binary string under the node denotes the path from the root to the node. 
In this example, at any node, the size of the subtree on the branch fixing a variable to 0 is larger than that on the branch fixing a variable to 1. 
This is due to the fact that infeasibility is more frequent on the branch fixing a variable to 1. 
However, considering all branching constraints at infeasible nodes, we cannot infer this asymmetrical trend. Depending upon the node's path, some infeasible leaves (highlighted by the green box) have a majority of branching constraints fixing variables to 0, whereas at some other leaves (highlighted by the yellow box), a majority of branching constraints fix variables to 1. Therefore, it is crucial to identify minimally infeasible subsystems of branching constraints. 

\section{Cardinality Constrained Problems} \label{sec:cardinality_constr}
% \chapter{Cardinality Constrained Problems} \label{sec:cardinality_constr}
Here, we focus on problems with a cardinality constraint enforcing an upper bound on the number of variables that can be non-zero. As a special case of packing constraint, it allows us to simplify the analysis and construct a specialized model for branching. Moreover, this class of problems is also important due to its relevance in several real-world problems. 

Consider an instance with a cardinality constraint and a node in its branch-and-bound tree where $n$ binary variables are free. Further, let the cardinality constraint along with the local branching constraints at the node dictate that at most $k$ of these free variables may take value 1. To model the asymmetry in the tree size on the two branches, we consider a simplified estimator $T(n, k)$ that ignores all other instance data and models the size of the subtree rooted at this node as a function of $n$ and $k$ only. Now, consider branching on any variable. On the branch fixing the variable to 0, we have $n-1$ free variables of which $k$ can be positive. Whereas on the branch fixing the variable to 1, $k-1$ of the $n-1$ free variable can be positive. This yields the following recursive relationship, 
\begin{equation}
    T(n, k) = T(n-1, k) + T(n-1, k-1) + 1
\end{equation}
Due to its natural interpretation of choosing $k$ out of $n$ items, this equation is approximately satisfied by $T(n, k) = {\binom{n}{k}}$. Thus, the relative weights of the two branches, $f_0, f_1$, modelled by the size of the respective subtrees as a fraction of the total size are, 
\begin{align*}
    f_0 = \frac{T(n-1, k)}{T(n, k)} = \frac{\binom{n-1}{k}}{\binom{n}{k}} = \frac{n-k}{n}, \\
    f_1 =  \frac{T(n-1, k-1)}{T(n, k)} = \frac{\binom{n-1}{k-1}}{\binom{n}{k}} = \frac{k}{n},
\end{align*}

We update the generalized score in Eq.\eqref{eq:asymmetrical_score} by choosing $a_0 = f_0/4$ and $a_1 = f_1/4$ and refer to this as the \textit{efficacious cardinality score} (Eff-Card). 

We then run branch-and-bound with this score on the two problems in our set that have cardinality constraints - Portfolio Optimization (CCP) and Weighted set coverage. We compare Eff-Card with LA-SB~(3,7,0.15) as it has a superior performance on these problems than Eff-SB~(3,7). 
Table \ref{tab:effcard} presents the geometric mean of the tree sizes of Eff-Card in comparison with LA-SB~(3,7,0.15) and Def-SB, and Table \ref{tab:eff_card_pct} presents the same results in terms of percentage reduction with respect to Def-SB. 
Eff-Card results in smaller tree sizes than both of the other scores and this decrease in tree size is particularly remarkable for suboptimal primal bounds, indicating improved robustness of the score.

\begin{table}
\begin{center}
		
		\caption{Evaluating Eff-Card based on shifted geometric mean of tree sizes rounded to integers. Smaller tree sizes between LA-SB(3,7,0.15) and Eff-Card are highlighted in bold.}\label{tab:effcard}
\resizebox{1.\textwidth}{!}{
\setlength{\tabcolsep}{0.5em} 
  {\renewcommand{\arraystretch}{1.1}
\begin{tabular}{|l|r|rrrr|rrrr|}
\hline
Rule                          & \multicolumn{1}{c|}{Def-SB} & \multicolumn{4}{c|}{LA-SB (3,7,0.15)}                                                                   & \multicolumn{4}{c|}{Eff-Card}                                                                           \\ \hline
Primal Gap                    & \multicolumn{1}{c|}{}       & \multicolumn{1}{c}{0\%} & \multicolumn{1}{c}{2\%} & \multicolumn{1}{c}{5\%} & \multicolumn{1}{c|}{10\%} & \multicolumn{1}{c}{0\%} & \multicolumn{1}{c}{2\%} & \multicolumn{1}{c}{5\%} & \multicolumn{1}{c|}{10\%} \\ \hline
\multicolumn{1}{|c|}{Problem} & \multicolumn{1}{r|}{}       & \multicolumn{1}{r}{}    & \multicolumn{1}{r}{}    & \multicolumn{1}{r}{}    & \multicolumn{1}{r|}{}     & \multicolumn{1}{r}{}    & \multicolumn{1}{r}{}    & \multicolumn{1}{r}{}    & \multicolumn{1}{r|}{}     \\ \hline
Portfolio Optimization (100)  & 19965                       & 18257                   & 18491                   & 18516                   & 18516                     & \textbf{17775}          & \textbf{17962}          & \textbf{17977}          & \textbf{17977}            \\
Weighted Set Coverage (200)   & 1274                        & 1309                    & 1415                    & 1415                    & 1415                      & \textbf{1235}           & \textbf{1246}           & \textbf{1246}           & \textbf{1246}             \\ \hline
Geo Mean                      & 5151                        & 4986                    & 5207                    & 5212                    & 5212                      & \textbf{4785}           & \textbf{4831}           & \textbf{4833}           & \textbf{4833}             \\ \hline
\end{tabular}
}}
% \end{center}
% \end{table}
\vspace{24pt}

% \begin{table}
% \begin{center}		
		\caption{
        Comparison of the shifted geometric mean of tree sizes to evaluate Eff-Card. Tree sizes are presented for Def-SB, and the percentage reduction in tree size relative to Def-SB is presented for others.}\label{tab:eff_card_pct}
\resizebox{1.\textwidth}{!}{
\setlength{\tabcolsep}{0.5em} 
  {\renewcommand{\arraystretch}{1.2}
% Please add the following required packages to your document preamble:
% \usepackage[table,xcdraw]{xcolor}
% Beamer presentation requires \usepackage{colortbl} instead of \usepackage[table,xcdraw]{xcolor}
\begin{tabular}{|l|r|rrrr|rrrr|}
\hline
Rule                          & \multicolumn{1}{c|}{Def-SB} & \multicolumn{4}{c|}{LA-SB (3,7,0.15)}                                                                   & \multicolumn{4}{c|}{Eff-Card}                                                                           \\ \hline
Primal Gap                    & \multicolumn{1}{c|}{}       & \multicolumn{1}{c}{0\%} & \multicolumn{1}{c}{2\%} & \multicolumn{1}{c}{5\%} & \multicolumn{1}{c|}{10\%} & \multicolumn{1}{c}{0\%} & \multicolumn{1}{c}{2\%} & \multicolumn{1}{c}{5\%} & \multicolumn{1}{c|}{10\%} \\ \hline
\multicolumn{1}{|c|}{Problem} & \multicolumn{1}{r|}{}       & \multicolumn{1}{r}{}    & \multicolumn{1}{r}{}    & \multicolumn{1}{r}{}    & \multicolumn{1}{r|}{}     & \multicolumn{1}{r}{}    & \multicolumn{1}{r}{}    & \multicolumn{1}{r}{}    & \multicolumn{1}{r|}{}     \\\hline
Portfolio Optimization (100)  & 19965                       & 9\%                     & 7\%                     & 7\%                     & 7\%                       & \textbf{11\%}           & \textbf{10\%}           & \textbf{10\%}           & \textbf{10\%}             \\
Weighted Set Coverage (200)   & 1274                        & -3\%                    & -11\%                   & -11\%                   & -11\%                     & \textbf{3\%}            & \textbf{2\%}            & \textbf{2\%}            & \textbf{2\%}              \\ \hline
Geo Mean                      & 5151                        & 3\%                     & -1\%                    & -1\%                    & -1\%                      & \textbf{7\%}            & \textbf{6\%}            & \textbf{6\%}            & \textbf{6\%}              \\ \hline
\end{tabular}
}}
\end{center}
\end{table}

\section{Description of Random Instances} \label{sec:random_instances}
In this section, we discuss the problems considered for computational experiments and the details of generating randomized instances.\\

\noindent
\textbf{Multi-dimensional Knapsacks} Consider the following multi-dimensional knapsack problem,
\begin{eqnarray*}
	&\textup{max}& {c^{\top}}{x} \\
	&\textup{s.t.} &a_i^{\top} x \leq b_i  \ \forall \ i \in \{1, \hdots, m\} \\
	&& x \in \{0, 1\}^n.
\end{eqnarray*}
We chose $n=100$ and $m=50$.
The constraint matrix is sparse and every $a_{ij}$ is independent and identically distributed.
Each element $a_{ij}$ is $0$ with probability $0.95$. Otherwise, it is an integer randomly selected from the set $\{1, 2, \dots, 200\}$ with uniform probability. Each component of the objective function is independent and uniformly picked from a range of integers $c_i \in \{1, 2, \dots, 200\}$.
The capacity parameter $b_i$ is a fraction of the sum of weights of the constraint, rounded down to an integer value. This fraction is $25\%$ for problems labelled small, $50\%$ for medium and $75\%$ for large. \\

% \noindent
% \textbf{Covering IPs} Similar to Packing IPs, the problem is formulated as follows, 
% \begin{eqnarray*}
% 	&\textup{min}& {c^{\top}}{x} \\
% 	&\textup{s.t.} &a_i^{\top} x \geq b_i  \ \forall \ i \in \{1, \hdots, m\} \\
% 	&& x \in \{0, 1\}^n.
% \end{eqnarray*}
% All parameters are chosen in the same manner as Packing IP.
% This fraction of $b_i$ to sum of weights in the constraint is $75\%$ for problems labelled strict, $25\%$  for problems labelled loose and $50\%$ otherwise. \\

\noindent
\textbf{Constrained Lot-sizing}
We consider the classical lot-sizing problem of determining production volumes while minimizing production cost, fixed cost of production and inventory holding cost across the planning horizon. In addition, the maximum quantity produced in any time period is constrained. The MILP model for the lot-sizing problem with $n$ time periods is as follows,
\begin{eqnarray}
	& \textup{min} & \sum_{i=1}^{n} \left( c_i \, x_i + f_i \, y_i \right ) +   \sum_{i=1}^{n -1} h_i \, s_i   \\
	&\textup{s.t.}& {x_1}{ = s_1 + d_1}{} \label{lotsize:flow0} \\
	&& s_{i-1} + x_{i} = d_{i} + s_{i}  \ \forall \ i \in \{2, \dots, n\} \\
	&& s_{n-1} + x_n = d_n \label{lotsize: flown}\\
	&& x_i \le \left(\sum_{j=i}^n d_j \right) \, y_i \ \forall \ i \in \{1, \hdots, n\} \\
        && x_i \le u_i \, y_i \ \forall \ i \in \{1, \hdots, n\} \\
	&& x \in  \mathbb{R}^n_+, \ s \in  \mathbb{R}^{n -1}_+, \  y \in  \{0, 1\}^n.
\end{eqnarray}
where variable $x_i$ is the quantity produced in period $i$ and $y_i$ is a binary variable with value $1$  if production occurs in period $i$ and $0$ otherwise. Lastly, variable $s_i$ is the inventory at the end of period $i$. In our computational experiments, we generate problems with $n=30$. 

Unit cost of production for each period, $c_i$, is independently and uniformly sampled from $\{1, \hdots, 10\}$, fixed cost of production, $f_i$, from $\{300, \hdots, 600\}$ and unit inventory holding cost $h_i$ from $\{1, \hdots, 10\}$. Similarly, the demand for each time-period, $d_i$ is independent and uniformly distributed in $\{50, \hdots, 100\}$ and capacity $u_i$ from $\{150, \hdots, 250\}$. \\

\noindent
\textbf{Big-bucket Lot-sizing}
{We consider a variant of the lot-sizing problem referred to in literature as the big-bucket lot-sizing problem. Here, multiple products are produced using shared resources~\cite{quadt2008capacitated}.} We do not consider unit cost of production here. On the other hand, set-up time and processing time are considered to be constrained. 
%The problem is modelled as follows, 
%%\paragraph*{Parameters \\ \\}
The following are the parameters of the corresponding MILP model,\\ \\
\begin{tabular}{l l}
	$P$ & Number of products, $\mathcal{P} = \{1, \hdots, P\}$ \\
	$T$ & Number of time periods, $\mathcal{T} = \{1, \hdots, T\}$ \\
	$f_i^p$ & Fixed cost of producing product $p$ in period $i$ \\
	$h_i^p$ & Inventory holding cost of product $p$ in period $i$ \\
	$t_i^{s, p}$ & Set up time of product $p$ in period $i$ \\
	$t_i^{u, p}$ & Processing time per unit of product $p$ in period $i$ \\
	$C_i$ & Time available in period $i$ \\
	$z^p$ & Initial inventory of product $p$ at the beginning of planning horizon. \\
\end{tabular}
\newline \newline The variables used to model the problem are the following,\\ \\
%%\paragraph*{Variables \\ \\}
\begin{tabular}{l l}
	$x_i^p$ & Quantity product $p$ produced in period $i$ \\
	$s_i^p$ & Quantity of product $p$ stored as inventory at the end of period $i$ \\
	$y_i^p$ & Binary variable indicating if product $p$ was produced in period $i$ ($y_i^p = 1$ if $x_i^p > 0$)\\
\end{tabular}
\newline \newline The MILP model used for the big-bucket lot-sizing problem is described below.
%%\paragraph*{Model - Big bucket lot-sizing}
\begin{eqnarray*}
	&\textup{min} &  \sum_{i \in \mathcal{T}} \sum_{p \in \mathcal{P}} \left(f^p_i \, y^p_i + h^p_i \, s^p_i \right) \label{bblotsize: obj} \\
	& \textup{s.t.} & s^p_0 = z^p, \ p \in \mathcal{P}  \label{bblotsize: s_init}\\
	&& s^p_T = 0, \ p \in \mathcal{P} \label{bblotsize: s_fin} \\
	&& s^p_{i-1} + x^p_{i}  = d^p_{i} + s^p_{i},  i \in \mathcal{T},  p \in \mathcal{P} \label{bblotsize: flow} \\
	&& x^p_i  \leq \left(\sum_{j=i}^n d^p_j \right) \, y^p_i, i \in \mathcal{T}, p \in \mathcal{P} \label{bblotsize: x-y} \\
	&& \sum_{p \in \mathcal{P}} \left( t^{s,p}_i \, y^p_i + t^{u,p}_i \, x^p_i \right) \leq C_i,  i \in \mathcal{T} \label{bblotsize: cap} \\
	&& x^p_i \in  \mathbb{R}_+, s^p_i \in  \mathbb{R}_+, y^p_i \in  \{0, 1\},\quad i \in \mathcal{T}, p \in \mathcal{P}.
\end{eqnarray*}
Our experiments we consider 10 time periods ($T=10$) and 3 products ($P=3$). 
Set up time for a product in each period, $t^{s, p}_i$ is independently sampled from $\{200, \hdots, 500\}$ with equal probability, unit processing time $t^{u, p}_i$ from $\{1, \hdots, 10\}$, time limitation $C_i$ from $\{1500 , \hdots, 3000 \}$, fixed cost of production, $f_i^p$, from $\{300, \hdots, 600\}$ and unit inventory holding cost $h_i^p$ from $\{1, \hdots, 10\}$.  
Similarly, the demand for each time-period, $d_i^p$ is independent and uniformly distributed in $\{0, \hdots, 100\}$ and the initial inventory for each product $z^p$ is sampled from $\{0, \hdots, 200\}$.\\

\noindent
\textbf{Maximum Weighted Matching}
The matching problem on graphs $G = (V, E)$ involves selecting a subset of edges $F \subset E$ such that at most one edge in $F$ is incident on any vertex in $V$.
The objective is to maximize the sum of weights on the edges. It is formulated as follows, 
\begin{eqnarray*}
	&\textup{max}& \sum_{e \in E} w_e x_e \\
	&\textup{s.t.}& \sum_{e \in E: u \in e} x_e \le 1, \qquad \forall u \in V \\
	&& x \in \{0,1\}^{|E|}
\end{eqnarray*}
The instances in this work have 300 nodes and 1000 edges. A random geometric graph is used to construct these instances. The position of the vertices is generated independently from a 2-dimensional uniform distribution. Pairs of nodes are sorted in the ascending order of the distance between them and the closest 1000 pairs are included as edges in the graph. The objective coefficients are independently picked from a uniform distribution in $[0, 1]$.\\

\noindent
\textbf{Maximum Weighted Coverage}
Consider a collection of sets $\mathcal{S}$ in a universe $\mathcal{U}$, that is $\mathcal{S} \subseteq 2^{\mathcal{U}}$. The goal of the maximum coverage problem is to pick $k$ sets from $\mathcal{S}$ such that the sum of weights of the elements covered by their union is maximized. The IP formulation of this problem is as follows, 
\begin{eqnarray*}
	&\textup{max}& \sum_{e \in \mathcal{U}} w_e y_e \\
	&\textup{s.t.}& y_e \le \sum_{S \in \mathcal{S}: e \in S} x_S \\
        && \sum_{S \in \mathcal{S}} x_S \le k \\
        && y \in [0, 1]^{|\mathcal{U}|} \\
	&& x \in \{0,1\}^{|\mathcal{S}|}
\end{eqnarray*}
The instances in this work have $|\mathcal{U}| = 1000$, $|\mathcal{S}| = 200$ and $k= 12$. Each set in the collection $\mathcal{S}$ is constructed by randomly including every element with probability $0.03$. The weights in the objective are independently and uniformly picked from $[0, 1]$.\\

\noindent
\textbf{Portfolio Optimization - Chance Constraint Programming}
We consider the probabilistically-constrained portfolio optimization problem for $n$ asset types, approximated by the sample approximation approach~\cite{pagnoncelli2009computational}, where the constraint on the overall return may be violated for at most $k$ out of the $m$ samples. The MILP formulation of this problem is as follows:
\begin{eqnarray*}
	&\textup{min} &\sum_{i = 1}^n x_i \\
	&\textup{s.t.} & a_i^T x + r z_i\ge r, \quad \forall i = 1, \hdots, m \\
	&& \sum_{i=1}^m z_i \le k \\
	&& x \in \mathbb{R}^n_+, \, z \in \{0,1\}^m. 
\end{eqnarray*}
We sample scenarios from the distribution presented in~\cite{qiu2014covering}, which is shown to be computationally difficult to solve. Each component of the constraint matrix, $a_{ij}$ is independently sampled from a uniform distribution in $[0.8, 1.5]$ and $r$ is equal to $1.1$. For our experiments, we set $n=20$, $m=100$ and $k=10$.\\

\noindent
\textbf{Fixed Charge Network Flow}
We consider the network flow problem where multiple commodities, need to be routed through a graph $G=(V, E)$. Let $k$ be the total number of commodities. Each commodity $p \in [k]$ has a single source $s_p$, a destination $t_p$ and a demand $d_p$ specifying the quantity that must be routed. The graph's edges are ordered pairs of nodes $(u, v)$ and allow bidirectional flow. There is a limited capacity $c_p$ on the total traffic flowing through the edge. Moreover, there is a fixed cost $f_e$ for using an edge and a cost $u_e$ per unit demand of any commodity flowing through it. The goal is to identify a feasible flow that satisfies the demand for all commodities while minimizing the total cost. The IP formulation of this problem is as follows, 
\begin{eqnarray}
	&\textup{min}_{x, y, z} & \sum_{e\in E} f_e y_e \, + \sum_{p \in [k], e\in E} d_p u_e z^p_e  \label{fcnf:obj}\\
	&\textup{s.t.}& \sum_{e \in \delta^+(v)} x^p_e - \sum_{e \in \delta^-(v)} x^p_e = \textbf{1}_{v = s_p} - \textbf{1}_{v = t_p}, \, \text{ for all } p \in [k], v \in V \label{fcnf:flow_conserv} \\
        & & \sum_{p \in [k]} d^p z^p_e \leq c_e y_e, \, \text{ for all } e \in E \label{fcnf:cap}\\
        & &  z^p_e \geq x^p_e, \, \text{ for all } p \in [k], \, e \in E \label{fcnf:abs1}\\
        & &  z^p_e \geq - x^p_e, \, \text{ for all } p \in [k], \, e \in E \label{fcnf:abs2}\\
	&& x^p_e \in [-1, 1] \, \text{ for all } p \in [k], \, e \in E  \\
 && z^p_e \in [0, 1] \, \text{ for all } p \in [k], \, e \in E  \\
 && y_e \in 
 \{0, 1\} \, \text{ for all } e \in E  \\
\end{eqnarray}
where variables $x^p_e$ represent the fraction of demand for $p$ flowing on $e$ with its sign representing the direction of the flow. The variable $z^p_e$ represents the absolute value of $x^p_e$ and $y_e$ is a binary variable indicating if an edge is used or not. The first component of Eq.(\ref{fcnf:obj}) models the fixed cost of using edges, and the second component models the total unit cost of the flow. The constraints Eq.(\ref{fcnf:flow_conserv}) ensure flow conservation of all commodities at all nodes and Eq.(\ref{fcnf:cap}) restrict the traffic through an edge to its capacity. Finally, Eq.(\ref{fcnf:abs1}, \ref{fcnf:abs2}) model the relationship between $x^p_e$ and $z^p_e$.

In our experiments, a random geometric graph with 50 nodes and 150 edges is constructed similarly as described in Maximum weighted matching. The problem has 3 commodities and the origin-destination pair for each is picked independently with uniform probability over the set of all ordered pairs of nodes. The demands are also independently and uniformly picked from the interval $\{100, \hdots, 300\}$. For each edge, the unit costs, fixed costs and capacities are independently picked from the intervals $\{3, \hdots, 10\}$, $\{3D, \hdots, 8D\}$ and $\{90, \hdots, 240\}$ respectively with uniform probability, where $D$ is the total demand across all commodities. 
\\

\noindent
\textbf{Set Packing and Set Covering}
We generate set-packing and set-covering problems with the random distribution used in \cite{yang2022learning}. 
We present it here for the sake of completion.
The set-covering problem is represented as, 
\begin{eqnarray*}
	&\textup{min}& {c^{\top}}{x} \\
	&\textup{s.t.} &a_i^{\top} x \geq 1  \ \forall \ i \in \{1, \hdots, m\} \\
	&& x \in \{0, 1\}^n.
\end{eqnarray*}
where the coefficients of $a_i$ are in $\{0, 1\}$.
Each objective coefficient $c_j$ is drawn from $\{1, \dots 100\}$ with uniform probability. The expected number of non-zeros in each constraint $n_i$ is selected from $\{ \frac{2n}{25} + 1, \hdots, \frac{3n}{25} -1 \}$. The coefficients $a_{ij}$ are then 1 with probability $\frac{n_i}{n}$.
The instances in our experiments are generated with $n=300$ and $m=3000$.
The set-packing problem is generated similarly with $n=200$ and $m=1000$.
\\

\noindent
\textbf{Independent Set on Graphs}
The unweighted maximum independent set problem involves selecting the largest subset of vertices of a graph $G$ such that no pair of vertices in the subset share an edge in $G$.
We use the instance generator of Ecole~\cite{prouvost2020ecole} to generate a clique cover formulation considered in~\cite{bergman2016decision, gasse2019exact}. Let $\mathcal{C}$ denote a set of cliques whose union covers all edges of a graph $G = (V, E)$. Then the maximum independent set problem on $G$ can be formulated as, 
\begin{eqnarray*}
	&\textup{max}& \sum_{v \in V} x_v \\
	&\textup{s.t.} & \sum_{v \in S} x_v \leq 1  \ \forall \ S \in \mathcal{C} \\
	&& x_v \in \{0, 1\}^{|V|}.
\end{eqnarray*}
Random graphs with 500 nodes are generated using Barabási–Albert model with affinity 4 and a greedy clique partitioning heuristic is used to obtain a clique partition.  \\

\noindent
\textbf{Capacitated Facility Location}
We consider the classical capacitated facility location problem~\cite{cornuejols1991comparison} with $m$ facilities and $n$ customer. The problem consists of selecting a subset of facilities to be open while ensuring that the open facilities meet the demand $d_i$ for every customer $i$. Each facility $j$ has a production capacity $s_j$ and a fixed cost $f_j$ is incurred for opening it. Additionally, a cost $c_{ij}$ is associated with transporting demand from facility $j$ to customer $i$. The objective is to minimize the total cost and is formulated as follows,
\begin{eqnarray*}
	&\textup{min}& \sum_{i=1}^n \sum_{j=1}^m c_{ij}x_{ij} + \sum_{j=1}^m f_j y_j \\
	&\textup{s.t.} & \sum_{j=1}^m x_{ij} = 1  \ \forall \ i \in \{1, \hdots, n\} \\
    && \sum_{i=1}^n d_i x_{ij} \leq s_j y_j  \ \forall \ j \in \{1, \hdots, m\} \\
    && x \in [0, 1]^{m \times n} \\
	&& y \in \{0, 1\}^m.
\end{eqnarray*}
Instances are generated with Ecole~\cite{prouvost2020ecole} that use the random distribution from~\cite{gasse2019exact} with 100 customers, and 100 facilities. The ratio of total capacities to total demand is selected as 5.
\\

\noindent
\textbf{Combinatorial auctions}
The combinatorial auctions problem presented in ~\cite{leyton2000towards, gasse2019exact} is considered and instances are generated using Ecole~\cite{prouvost2020ecole} following the ``arbitrary relationships'' scheme found in~\cite{leyton2000towards}. 
Participants bid on a collection of items and the problem consists of identifying the winning bids that maximize total price but do not repeat the constituting items. Let $\mathcal{B}$ be the set of all bids and $\mathcal{S}$ be the set of all items. Each bid $b$ has a price $p_b$ and includes items $S_b \subseteq \mathcal{S}$. The problem can then be formulated as the following set-packing problem, 
\begin{eqnarray*}
	&\textup{max}& \sum_{b \in \mathcal{B}} p_b x_b \\
	&\textup{s.t.} & \sum_{b : i \in S_b} x_{b} \leq 1  \ \forall \ i \in \mathcal{S} \\
	&& x_b \in \{0, 1\} \ \forall b \in \mathcal{B}.
\end{eqnarray*}
Instances in this work have 200 items, 1000 bids and default values of all other parameters in the Ecole generator.

\end{document}